\documentstyle[epsf,12pt]{article}
\newtheorem{sect}{}[section]
\newtheorem{subsect}{}[sect]

\newcommand{\sqt}{\mbox{ $\Box \hspace{-9.1pt} \raisebox{1.1pt}{$\times$} $} }

\newcommand{\qsqsj}[6]{\mbox{$ \left[ \begin{array}{ccc}
                                    #1 & #2 & #3 \\ #4 & #5 & #6
                                    \end{array} \right]_q   $ } }

\def\ox{\otimes}

\begin{document}

\title{Structures and Diagrammatics  of Four Dimensional
Topological Lattice Field Theories }

\author{
J. Scott Carter \\
Department of Mathematics  \\
University of South Alabama \\
Mobile, Alabama 36688 
\and
Louis H. Kauffman \\ 
Department of Mathematics \\
University of Illinois at Chicago \\
851 S. Morgan St. \\
Chicago, Illinois 60607-7045 
\and
Masahico Saito \\
Department of Mathematics \\ University of South Florida
\\Tampa, FL 33620 
}

\maketitle

\begin{abstract}
Crane and Frenkel proposed a 
state sum invariant   
for triangulated 4-manifolds.
They defined and used new algebraic structures called Hopf categories
for their construction. 
Crane and Yetter 
studied Hopf categories and gave some examples
using  
group cocycles
that are associated to the Drinfeld double of a finite group. 

In this paper we define a state sum invariant of triangulated $4$-manifolds
using Crane-Yetter cocycles as Boltzmann weights.
Our invariant 
generalizes  
the  3-dimensional invariants 
defined by
Dijkgraaf and Witten 
and the invariants that are defined via Hopf algebras.  
We 
present diagrammatic methods for the study of
such invariants 
that illustrate connections between Hopf categories 
and moves to triangulations. 

\end{abstract}

\newpage

\tableofcontents

\section{Introduction}

Witten's formulation \cite{Witten} of 
an intrinsic definition of the Jones polynomial \cite{Jones}
based on
physical models 
lead
 to 
the more rigorous 
mathematical definitions via representations of quantum groups
that were 
given by Reshetikhin,       
Turaev, and Viro \cite{RT}, \cite{TV}.
These 
quantum invariants
are speculated to 
 generalize 
to higher dimensions.
Such 
putative invariants 
have  
their origins in a theory of  
quantum gravity \cite{JBOx} 
and higher categories \cite{JB:form}.
In relation to the current work, the following progress has been made.

Quantum spin 
networks were 
generalized 
by Crane-Yetter \cite{CY} 
to give 
$4$-manifold 
invariants that were based on 
Ooguri's proposal \cite{Oog}. 
The invariants can be used to compute the signature as shown in
\cite{JRob,CKY1,CKY2}.  
Birmingham-Rakowski \cite{BR} generalized the
 Dijkgraaf-Witten \cite{DW}
 invariant of 3-manifolds,
defined by group 3-cocycles,  to 
triangulated 4-manifolds using 
pairs of cocycles. 
Crane and Frenkel \cite{CF} 
constructed 
Hopf categories to define 4-manifold invariants, and 
they gave examples 
using canonical bases of quantum groups. 
In \cite{CYex} Crane and Yetter used cocycles to construct 
Hopf categories. 

In this paper we provide direct relations between 
the 
cocycle 
conditions of \cite{CYex} and Pachner moves of $4$-manifolds,
thus 
constructing 
a generalization of the  Dijkgraaf-Witten
invariants to dimension 4. The relations are established diagrammatically,
providing  
connections 
 between Hopf categorical structures and 
triangulations via dual graphs and their movies. 

The current paper is self-contained, but the reader might enjoy our
 introduction to the subject given in \cite{CKS}, 
where many of our  ideas and motivations are introduced
 in a more leisurely fashion. 
For the diagrammatic foundation of the invariants in dimension 3 see 
\cite{K&P}, \cite{KL}, and \cite{K&D}. 
For the algebraic approach see \cite{CFS}.
 Finally, there is a relation to higher dimensional knot theory as found in \cite{CS:book}.

Let us continue our motivational remarks. 
In $3$ dimensions, 
planar diagrams played a key role 
in 
the definitions of 
both 
knot invariants and manifold invariants.
Such diagrams are convenient 
since 
they help one  
grasp 
the  
categorical 
and algebraic structures needed for defining invariants.
One of the difficulties in generalizing to dimension $4$ or higher
is the lack of such visualizations and diagrammatic 
machinery. 
The purpose of this paper is to 
provide basic diagrammatic tools to study $4$-manifold 
triangulations, and to use such formulations to define invariants.

In particular, we  formulate the Crane-Frenkel approach
in terms of cocycles as initial data and prove the 
invariance under Pachner moves
in a diagrammatic way.
We introduce spin networks for  the study of
such invariants. 
We hope that the present work 
serves 
as a basic tool  
in exploring the possibilities in higher dimensions.

There is a close relationship bewteen certain physical models in
statistical mechanics and quantum field theory and the formulation of
``quantum'' invariants of knots, links and three-manifolds. We hope that this
relationship continues into dimension four. In particular, one can hope
that a four dimensional toplogical field theory (such as we study here)
would be related to the calculation of amplitudes in quantum gravity. The
naive reason for this hope is quite simple: An amplitude for quantum
gravity must sum over the possible metrics on the four-dimensional
spacetime. Averaging over metrics should in priniciple produce basic
numbers that are metric-independent. In other words a valid process of
averaging over metrics should produce topological invariants of the
underlying four-space. Of course on the physical side it will be neccessary
to extricate the topological part from the complexities of the model. On
the pure mathematical side it will be neccessary to see the relevance of
the mathematics. Nevertheless this hope for an application to quantum
gravity is one of the forces that drives our project.

The organization of the paper is as follows.
In Section~\ref{23}, we review state sum invariants for triangulated
manifolds in dimensions $2$ and $3$. We emphasize
diagrammatic relations 
between triangulations and algebraic structures.
At the end of the section, we summarize the idea of categorification
in  relation to the construction of higher dimensional invariants.
In Section~\ref{Pac4sec} we present
diagrams
of Pachner moves in 
dimension $4$. We also introduce singular moves in dimension $4$, and
prove that singular moves together with $3$-dimensional Pachner moves
imply $4$-dimensional Pachner moves. These Lemmas will be used to prove the 
well-definedness of our invariants.
In Section~\ref{T&Dsec} we give generalization of spin networks to
 dimension $4$. Triangulations are represented by movies of graphs,
and 
these graph movies 
will be used to give a direct relation 
between Hopf category structures and triangulations.
Cocycle conditions defined by \cite{CYex} will be reviewed 
in Section~\ref{data}. Symmetry of cocycles are defined.
In Section~\ref{bzman} the state sum invariants will be defined,
and will be proved to be well-defined in Section~\ref{invsec}.
Our proofs are diagrammatic. 
They provide the
basic machinery necessary to define other invariants
defined via Hopf categories.
The axioms of Hopf categories are related to moves on 
triangulations of 4-manifolds in a manner 
similar to the relationship between Hopf 
algebras and moves on 3-manifolds.  Section~\ref{Hopfsec} reviews
 the definition of Hopf category and make the connection with the
 rest of the paper.

\section{Quantum 2- and 3- manifold invariants
}
\label{23}

In this section, we review topological lattice field theories 
in dimension $3$ and explain how they are generalized from 
those in 
dimension $2$. 
First we review dimension $2$ following \cite{FHK,ChFS} 
where semi-simple algebras are used. An alternative
approach is given  by Lawrence in \cite{Ruth} 
in which the algebra is assumed to be Frobinius. 
Next the Turaev-Viro \cite{TV} theory is reviewed following 
\cite{CFS} and \cite{KL}.
 Invariants of 3-manifolds derived from Hopf algebras are presented
 following \cite{ChFS}. 
Alternative approaches are found in Kuperburg~\cite{Greg} and
 Kauffman-Radford~\cite{KR}.
Some of the summary appeared in \cite{CKS}.
We summarize Wakui's defintion \cite{Wakui} of 
the Dijkgraff-Witten invariants \cite{DW}, 
but here we show invariance using the Pachner Theorem.
This section closes with a conceptual scheme for generalizing to dimension 4.

\begin{sect}{\bf Topological lattice field theories in dimension $2$.\/}
\addcontentsline{toc}{subsection}{Topological lattice field theories in dimension $2$}
{\rm
Let $A$ denote a finite dimensional associative 
algebra over the complex numbers ${\bf C}$. 
Let $\{ \phi_i | i=1, \cdots, n \}$ denote an ordered basis for $A$, and 
for $x,y,z \in \{ 1, \ldots, n\}$,
let 
$C_{xy}^z$ denote {\it the structure constant} of the
algebra $A$.  Thus 
the multiplication between basis elements is given by the formula:
$$\phi_x \cdot \phi_y = \sum _z C_{xy}^z \phi_z.$$
Apply the associativity law, $(ab)c=a(bc)$, to the basis elements as follows:
\begin{eqnarray*}
(\phi_a \phi_b)\phi_c & = & (\sum_j C_{ab}^j \phi_j)\phi_c
 =   \sum_{j, d} C_{ab}^j C_{jc}^d \phi_d \\ 
\phi_a (\phi_b \phi_c) & = & \phi_a (\sum_i C_{bc}^i \phi_i)
 =   \sum_{i, d} C_{ai}^d C_{bc}^i \phi_d .
\end{eqnarray*}

In this way, we obtain the equation 
$$ \sum_j C_{ab}^j C_{jc}^d =     
\sum_i  C_{ai}^d C_{bc}^i $$      
whose geometrical interpretation will be presented shortly.

For $x,y \in \{1, 2, 
\ldots, n (=dim A) \}$, define 
$$g_{xy}= \sum_{u,v} C^v_{ux} C^u_{vy}.$$ 
 Then this is invertible
precisely when the algebra $A$ is semi-simple \cite{FHK},
 and  the matrix inverse 
$g^{xy}$ of $g_{xy}$ defines a bilinear form on the algebra $A$.
The geometric interpretation of this bilinear form and that of the 
associativity identity will allow us to define from a semi-simple 
associative algebra, an invariant of 2-dimensional manifolds.

We follow the definition 
given in  
\cite{FHK}.
Let $T$ be a triangulation of a closed 2-dimensional
manifold $F$.
Let ${\cal N}= \{ 1, 2, \cdots, n \} $.
This is called the set of {\it spins}. 
Let ${\cal ET}= \{ (e,f) | e \subset f \}$ 
be the set 
of all the pairs of edges, $e$, and faces, $f$, 
such that $e$ is an edge of $f$. The set ${\cal ET}$ is a partial flag.
A {\it labeling} is a map $L: {\cal ET} \rightarrow {\cal N}$.
Thus a labeling is an assignment of spins to all the edges
with respect to faces. 
Given a labeling, we assign weightings to 
faces and edges as follows:
Suppose that we are given functions
$C$ and $g$,                                 
$C: {\cal N}^3 \rightarrow {\bf C}$, $C(x,y,z)= C_{xyz}$,
and $g: {\cal N}^2 \rightarrow {\bf C}$, $g(x,y)=g^{xy}$. 
If a face has three edges labeled with spins $x,y,z$,  
then assign the complex number $C_{xyz}$ to the face.
It is assumed that 
the function $C$ posesses a cyclic symmetry; so $C_{xyz}$ $=C_{yzx}$ 
$=C_{zxy}$.
If an edge is shared by two faces, and the edge with respect 
to these faces receives spins $x$ and $y$, then assign 
the complex number $g^{xy}$ to the edge. 
Then define a {\it partition function} $\Psi(T) $
by $\sum_{L} \prod  C_{xyz} g^{uv}$
where the sum is taken over all the labelings and the product 
is taken over all the 
elements of ${\cal ET}$.
Values of $C_{xyz}$ and $g^{uv}$
are
given  in terms of the structure constants for the 
algebra and the bilinear form $g_{uv}$ in the fourth paragraph forward. 
First, we discuss topological aspects that motivate their definition.

In order for the partition function to be topologically invariant,
it cannot depend on the choice of triangulation.
There are two steps in constructing such an invariant quantity.
First,
we work topologically.
There is a  set of local moves to triangulations that suffices  
to relate any two triangulations 
of a given manifold.
These moves were discovered by Pachner \cite{Pac} in the general case of 
$n$-manifolds, and they  generalize 
a classical theorem of Alexander 
\cite{LX}. 
Therefore for the 
partition function 
to be independent of 
the choice of triangulation, it is sufficient to prove that 
the weighting assigned to triangles and edges
satisfies equations 
that  
correspond
to these local moves. 
The second step, therefore, 
is algebraic.
We seek 
functions 
$C$ and $g$ that satisfy these equations. We will indicate that the 
structure constants of an associative algebra can be used for the 
function $C$ and that the bilinear form on $A$ can be used to define 
the function $g$, 
as the notation suggests. 

Let us consider the topological aspects.
The Pachner moves in dimension 2 are depicted in 
Figure~\ref{2dpac}.
 The move on the left 
of Figure~\ref{2dpac} is called
the $(2\rightleftharpoons 2)$-move; 
that on the right is called the $(1 \rightleftharpoons 3)$-move.
The names of the moves indicate the number of triangles that are involved.

\begin{figure}
\begin{center}
\mbox{
\epsfxsize=4in
\epsfbox{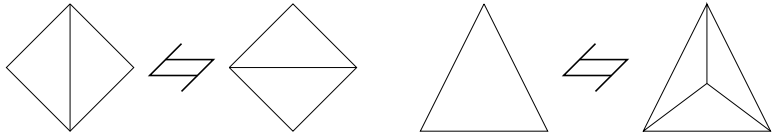}
}
\end{center}
\caption{2D Pachner moves}
\label{2dpac}
\end{figure}

\begin{figure}
\begin{center}
\mbox{
\epsfxsize=3in
\epsfbox{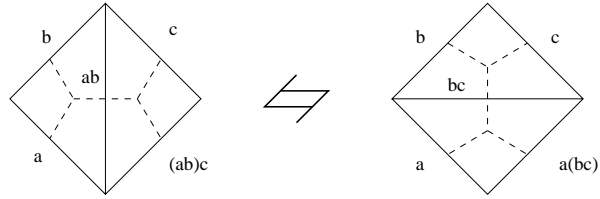}
}
\end{center}
\caption{Associativity and a 2-dimensional Pachner  move}
\label{2ass}
\end{figure}

We now 
intepret 
associativity and the bilinear form
in a semi-simple algebra over ${\bf C}$
in terms of 
the  Pachner moves.
Specifically, 
the $(2 \rightleftharpoons 2)$-Pachner moves is related to
the associativity law
$(ab)c= a(bc)$.
The relationship is depicted in Figure~\ref{2ass}.
The dual graphs,
indicated in the
Figure by dotted segments, are
sometimes useful for visualizing
the relations between triangulations and the algebraic structure.
The diagram given in Figure~\ref{bubble} illustrates the geometrical 
interpretation of the bilinear form $g_{xy}= \sum_{u,v} C^v_{xu} 
C^u_{yv}.$ 
In the figure, two triangles share two edges in the left picture,
representing the local weighting $\sum_{u,v} C^v_{xu} C^u_{yv}$,
and the right represents a single edge corresponding to  $g_{xy}$.
Finally, this relationship together with the associativity 
identity can be used to show that the partition function is invariant 
under the $(1 \rightleftharpoons 3 )$-Pachner move. The essence of the 
proof is indicated in Figure~\ref{bubbled}.

\begin{figure}
\begin{center}
\mbox{
\epsfxsize=3in
\epsfbox{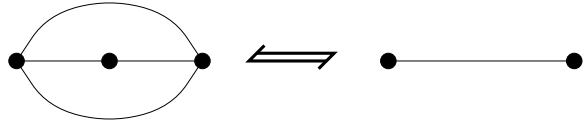}
}
\end{center}
\caption{The semi-simplicity axiom and degenerate triangulations}
\label{bubble}
\end{figure}

\begin{figure}
\begin{center}
\mbox{
\epsfxsize=3in
\epsfbox{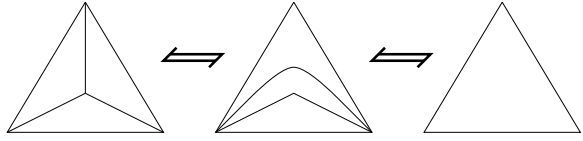}
}
\end{center}
\caption{Semi-simplicity, associativity, and the (3,1)-move}
\label{bubbled}
\end{figure}

Having illustrated the  
algebra axioms diagrammatically,
 we turn to 
show how the structure constants and the bilinear form
of associative semi-simple algebras solve the equations corresponding 
to the Pachner moves. Given, structural constants $C^z_{xy}$ and a 
non-degenerate bilinear form $g^{uz}$ with inverse $g_{uz}$,
define $C_{xyu}$ by the equation (using Einstein summation convention
of summing over repeated indices),                     
$$ C_{xyu} \equiv g_{uz}C_{xy}^z .$$ 
Then since
$$ \sum_j C_{ab}^j C_{jc}^d =   
\sum_i  C_{ai}^d C_{bc}^i    $$  
the partition function defined in this way is invariant under 
the $(2 \rightleftharpoons 2)$-move.
Furthermore, we have (again, under summation convention) 
$$C^a_{de}C^j_{ab}C^d_{jc} = C^a_{de}C^d_{ai}C^i_{bc} = g_{ie}C^i_{bc} 
= C_{ebc},$$
and so the partition function is invariant under the 
$(1 \rightleftharpoons 3)$-move. In this way, a semi-simple 
finite dimensional algebra 
defines an invariant of surfaces. On the other hand,          
given a partition 
function one can define a semi-simple algebra with these structure 
constants and that bilinear form. In \cite{FHK}, this is stated as 
Theorem 3:

\noindent
{\it The set of all TLFTs is in one-to-one correspondence with the set 
of finite dimensional semi-simple associative algebras.}

Observe that the $(1 \rightleftharpoons 3)$-move follows from the 
$(2 \rightleftharpoons 2)$-move and a non-degeneracy condition. In the 
sequel, we will see similar phenonema in dimensions 
 3 and 4.


In general, the idea of defining a partition function to produce a 
manifold invariant is (1) to assign spins to simplices of a 
triangulation, and (2) to find weightings that satisfy equations 
corresponding to Pachner moves. This approach, of course, depends on 
finding such solutions to (often extremely overdetermined) equations. 
Such solutions come from certain algebraic structures. Thus one hopes 
to extract appropriate algebraic structures from the Pachner moves in 
each dimension.  
This is the motivating philosophy of quantum topology. 

In the following sections we review such invariants in 
dimension $3$ in more detail to explain such relations
between triangulations and algebras.

}\end{sect}

\begin{sect}{\bf Pachner moves in dimension  $3$.\/} 
\label{pac23}
\addcontentsline{toc}{subsection}{Pachner moves in dimension  $3$\/}
{\rm 
In this section
we  review the Pachner moves \cite{Pac}
of triangulations of manifolds in dimension $3$.
The Pachner moves in $n$-dimensions form 
a set of moves
on triangulations such that any two different triangulations
of a manifold can be related by a sequence of moves from this set.
Thus, two triangulations represent the same manifold if and only if
one is obtained from the other by a finite sequence of such moves.
In Figure~\ref{3dpac} the  
3-dimensional Pachner moves are depicted, these are called
the 
$(2 \rightleftharpoons 3)$-move 
and  the 
$(1 \rightleftharpoons 4)$-move.

\begin{figure}
\begin{center}
\mbox{
\epsfxsize=4in
\epsfbox{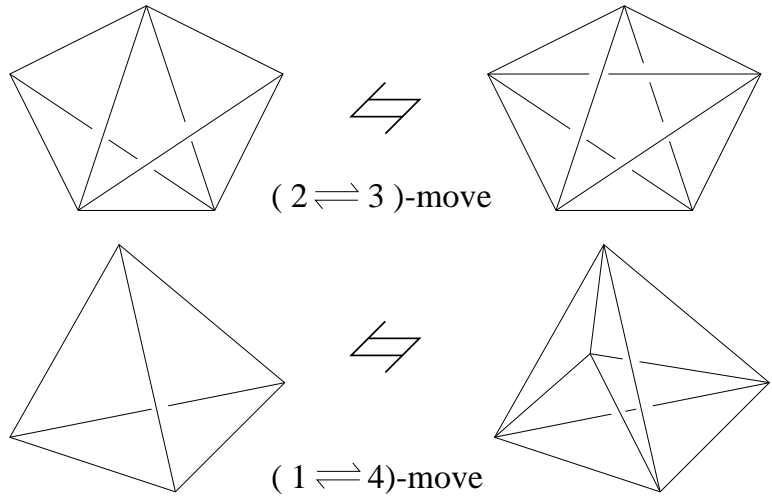}
}
\end{center}
\caption{3-dimensional  Pachner moves}
\label{3dpac}
\end{figure}

\begin{figure}
\begin{center}
\mbox{
\epsfxsize=4in
\epsfbox{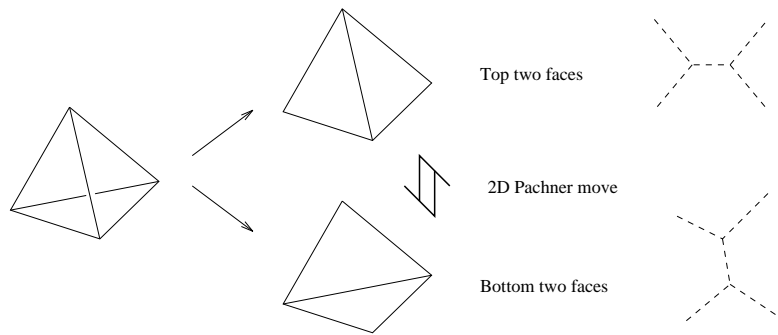}
}
\end{center}
\caption{Movie of a tetrahedron and a  2-dimensional  Pachner move}
\label{3ass}
\end{figure}

\begin{figure}
\begin{center}
\mbox{
\epsfxsize=5in
\epsfbox{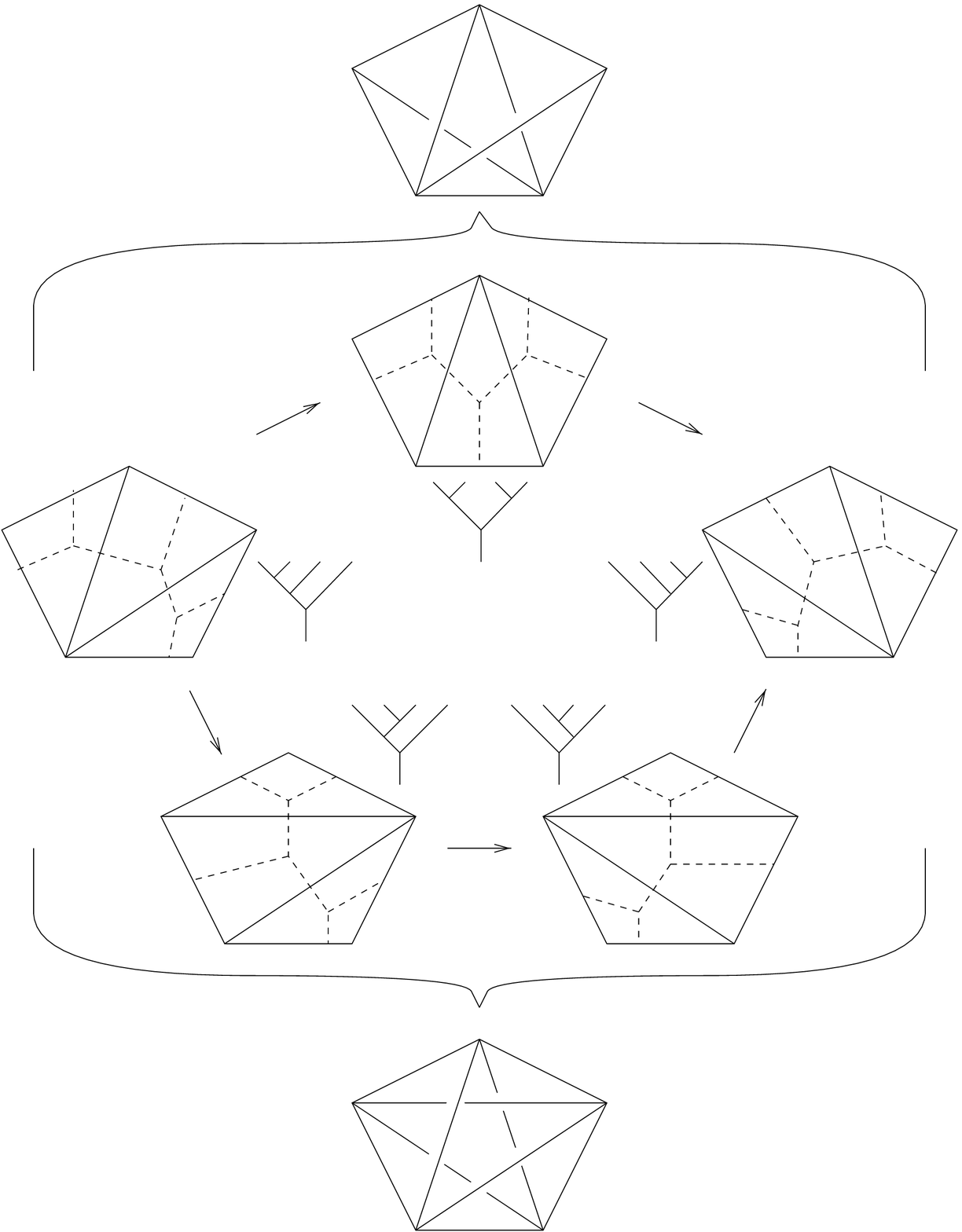}
}
\end{center}
\caption{The pentagon, trees, and a 3-dimensional Pachner move}
\label{penta}
\end{figure}

Notice that  the  2-dimensional Pachner moves
relate the  faces of a tetrahedron.
Specifically, the $(2 \rightleftharpoons 2)$-move consists of
two pairs of
triangles and they together form a tetrahedron
(Figure~\ref{3ass}). 
Meanwhile, the 
$(1 \rightleftharpoons 3)$-move 
relates three triangular faces of a tetrahedron  to the 
remaining face. The three triangles form the 
central projection of a tetrahedron. 
Analogous facts are true
for the
3-dimensional Pachner moves as well; let us explain.
One side of each move is a union of $3$-faces
of the boundary of a $4$-simplex and the other side of the
move is the rest of the $3$-faces, and they together form
the boundary of a $4$-simplex. For example, the 
$(1 \rightleftharpoons 
4)$-move  indicates two 3-balls on the 
boundary of a 4-simplex as they 
appear in a central projection of the simplex.

In Figure~\ref{3ass}, the relation between faces of a 
tetrahedron and their dual graphs is depicted.
The middle picture shows pairs of front and back faces
of a tetrahedron 
on the center left. 
Note that these pairs represent 
the 
$(2 \rightleftharpoons 2)$
Pachner move (as indicated by the vertical double arrow in the middle). 
Thus the  $(2 \rightleftharpoons 2)$ Pachner move corresponds to
a tetrahedron,  
a  $1$-dimensionally 
 higher simplex. 
On the right of the figure, the change on dual graphs is depicted
by dotted lines.
In Figure~\ref{penta}, a similar correspondence is depicted for 
the  $(2 \rightleftharpoons 3)$ Pachner move. 
Here faces of unions of tetrahedra are depicted from left to right,
in two different ways that correspond to the Pachner move.
These are the faces taken from the union of tetrahedra depicted
in top and bottom of the figure, respectively. 
The dual graphs are also depicted, which are the graphs used for
the Biedenharn-Elliott identity of the $6j$-symbols.
This direct diagrammatic correspondence 
is pointed out in \cite{CFS}.

}\end{sect}

\begin{sect}{\bf  Turaev-Viro invariants. \/}\addcontentsline{toc}{subsection}{Turaev-Viro invariants }
{\rm
One way to view the Turaev-Viro invariants \cite{TV,KL}
is to ``categorify'' the 
TLFT in dimension 2. In this process, the semi-simple algebra is 
replaced by a semi-simple monoidal category --- 
namely the category of representations of 
$U_q(sl(2))$ where $q$ is a primitive $4r\/$th root of unity.
First we review the definition of the Turaev-Viro invariants,
and then explain the viewpoint mentioned above.

A triangulation of a $3$-manifold is given.  
A 
coloring, $f$,  
is {\it admissible} if 
whenever 
edges with colors $a,b,j$ bound a 
triangle, then the triple $(a,b,j)$ is 
a
{\it $q$-admissible triple} in the 
sense that  
\begin{enumerate}
\item
$a+b +j$ is an integer, 
\item $a+b-j$, \ $b+j-a,$ and $a+j-b$  are all $\ge 0$
\item $a+b+j \le r-2.$
\end{enumerate} 

If edges with labels $a,b,c,j,k,n$ are the edges of a tetrahedron such 
that each of $(a,n,k)$, $(b,c,n),$  $(a,b,j)$, and $(c,j,k)$ is
a $q$-admissible triple, then the tetrahedron, $T$, receives a weight of 
$T_f= \qsqsj{a}{b}{n}{c}{k}{j}.$ If any of these is not admissible, then 
the weight associated to a tetrahedron is, by definition, 0.

For a fixed coloring $f$ of the  edges of the triangulation of a 3-
manifold $M$, the value 
$$|M|_f = \Delta^{-t}  \prod  \Delta_{f(E)} \prod T_f$$
is associated where $t$ is the number of vertices in the 
triangulation, the first product is taken over all the edges in the 
triangulation, the second product is over all the tetrahedra, the 
factor $\Delta$ is a normalization factor (= const.) 
and $\Delta_{f(E)}$ is a 
certain quantum integer associated to the color of the edge $E$. 
To obtain an invariant of the manifold one forms the sum
$$|M| = \sum_{f} |M|_f$$ where the sum is taken over all colorings.
Further details can be found in \cite{TV,KL} or \cite{CFS}.

Several points should be made here. First, the sum is finite because 
the set of possible colors is finite. Second, the quantity $|M|$ is 
a topological invariant because the $6j$-symbols satisfy the 
Beidenharn-Elliott identity and an orthogonality condition. The 
orthogonality is a sort of non-degeneracy condition on the 
$6j$-symbol. In \cite{KL,CFS} it is shown how to use orthogonality and 
Beidenharn-Elliott
(together with an identity among certain quantum integers) 
 to show invariance under the $(1 \rightleftharpoons 
4)$ move.
 Third, the $6j$-symbol is a measure of non-associativity as 
we now explain.

The
situation 
at hand 
can be seen as a categorification.
In $ 2$-dimensions  associativity $(ab)c=a(bc)$ played a key role. 
In $3$-dimensions the $6j$-symbols are defined by comparing two different 
bracketting 
$(V^a \ox V^b) \ox V^c$ and $V^a \ox (V^b \ox V^c)$
of representations $V^a$, $V^b$, and $V^c$.
Here the algebra elements became vector spaces as we went up 
 dimensions by one, and the symbol measuring the difference
in associativity satisfies the next order associativity,
corresponding to the Pachner move. 

Given representations $V^a$, $V^b$,  $V^c$,  we can form their triple 
tensor product and look in this product for a copy of the 
representation $V^k$. If there is such a copy, it can be obtained by 
regarding $V^k$ as a submodule of $V^a \ox V^n$ where $V^n$ is a 
submodule of $V^b \ox V^c$, or it can be obtained as a submodule 
of $V^j \ox V^c$ where $V^j$ is a submodule of $V^a \ox V^b$. From 
these two considerations we obtain two bases for the set of 
$U_q(sl(2))$ maps 
$V^k \rightarrow V^a \ox V^b \ox V^c$. The $6j$-symbol is the change of basis 
matrix between these two.

Considering such inclusions into four tensor products, 
we obtain the Biedenharn-Elliott identity. Each such inclusion is 
represented by a tree diagram. Then the Biedenharn-Elliott
identity is derived from the tree diagrams depicted in Fig. \ref{penta}.


}\end{sect}

\begin{sect}{\bf Invariants defined from Hopf algebras.\/ }
\addcontentsline{toc}{subsection}{Invariants defined from Hopf algebras}

{\rm 
 In this section we review invariants defined by 
Chung-Fukuma-Shapere \cite{ChFS} and Kuperberg \cite{Greg} 
(we follow the description in \cite{ChFS}).
We note that the invariants obtained in this section are
also very closely related to the invariants defined and
studied by Hennings, Kauffman, Radford and Otsuki
(see \cite{KR} for example). 
For background material on Hopf Algebras see \cite{sweed} or \cite{Mont}, for example.

\begin{subsect}{\bf Definition (Bialgebras).\/}  {\rm 
A {\it bialgebra} over a field $k$ is a quintuple
$(A, m, \eta, \Delta, \epsilon)$ such that

\begin{enumerate} 
\item
 $(A, m, \eta)$ is an algebra
where $m : A \ox A \rightarrow A$ is the multiplication
and $\eta  : k \rightarrow A$ is the  unit.
({\it i.e.}, these  
are $k$-linear maps such that
$m(1 \ox m)=m(m \ox 1)$,
$m(1 \ox \eta)= 1 = m(\eta \ox 1)$).

\item
 $\Delta : A \rightarrow A \ox A $ is an algebra
homomorphism (called the  {\it comultiplication})
satisfying $(id \ox \Delta ) \Delta = ( \Delta \ox id) \Delta $,

\item
 $\epsilon : A \rightarrow k$ is an algebra homomorphism
called the {\it counit},
satisfying $ (\epsilon \ox id)\Delta = id = (id \ox \epsilon)\Delta $.

\end{enumerate}
} \end{subsect}

\begin{subsect}{\bf Definition (Hopf algebras).\/ } {\rm
An {\it antipode} is a map $s : A \rightarrow A$
such that
$m \circ (s \ox 1) \circ \Delta = \eta \circ \epsilon
         = m \circ (1 \ox s) \circ \Delta $.

A {\it Hopf  algebra } is a bialgebra with an antipode.
} \end{subsect}

\begin{figure}
\begin{center}
\mbox{
\epsfxsize=3in
\epsfbox{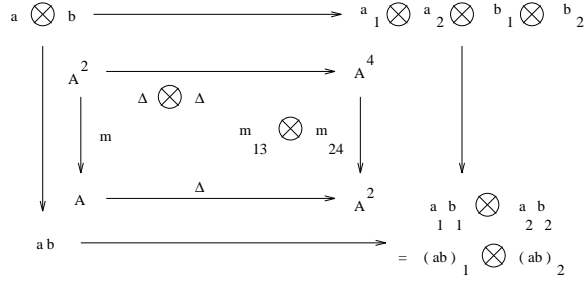}
}
\end{center}
\caption{The relation between multiplication and comultiplication}
\label{coh}
\end{figure}

The image of the comultiplication is often written as
$\Delta(a)= a_1 \ox a_2$ 
for $a \in A$. The image in
fact is a linear combination of such tensors but the
coefficients and the summation are abbreviated; 
this is the so-called 
Sweedler notation \cite{sweed}.
The most important property 
from the present point of view   
is the compatibility condition
between the multiplication and the comultiplication
({\it i.e.}, the condition that the comultiplication is an algebra homomorphism), 
and we
include 
the commuting diagram for this relation in Figure~\ref{coh}.
 The condition is written more specifically
$\Delta \circ m = (m \ox m) \circ P_{23} \circ (\Delta \ox \Delta)$
where $P_{23}$ denotes the permutation of the second and the third 
factor: $P_{23}(x \ox y \ox z \ox w) = (x \ox z \ox y \ox w )$.
In the Sweedler  
notation, 
it is also written as 
$(ab)_1 \ox (ab)_2 = a_1 b_1 \ox a_2b_2$.

The definition of invariants defined in \cite{ChFS} is
similar to the 
$2$-dimensional 
case. 
Given a triangulation $T$ of a $3$-manifold $M$,
give spins to  edges with respect to  faces (triangles).
The weights then are assigned to edges and  to faces.      
The structure constants $C_{xyz}$ (resp. $\Delta_{xyz}$)
 of multiplication (resp. comultiplication)
are assigned as weights to faces (resp. edges).
If an edge is shared by more than three faces, then a composition
of comultiplications are used. For example for four faces sharing
an edge, use the structure constant for $(\Delta \ox 1) \Delta$.
The coassociativity ensures that the other choice 
$(1 \ox \Delta)\Delta$ gives the same constant $\Delta_{v_1,v_2,v_3,v_4}$.
Thus the partition function takes the form 
$\Psi(T)= \sum_L \prod C_{xyz} \Delta_{v_1, \cdots, v_n}$.
This formula exhibits the form of the partition
function for this model, but is not technically complete. 
The full formula
uses the antipode in the Hopf algebra 
to take care of relative orientations
in the labellings of the simplicial complex.  
See \cite{ChFS} for the details.

In \cite{ChFS} the well definedness 
was proved by using
singular triangulations  
--- these  generalize 
triangulations by                                     
allowing certain cells as building blocks.
In this case the move called the {\it cone move}
for a singular triangulation plays an essential role.
This move is depicted in Figure~\ref{3dcone} with a dual graph to 
illustrate the relationship to the compatibility condition.

\begin{figure}
\begin{center}
\mbox{
\epsfxsize=4in
\epsfbox{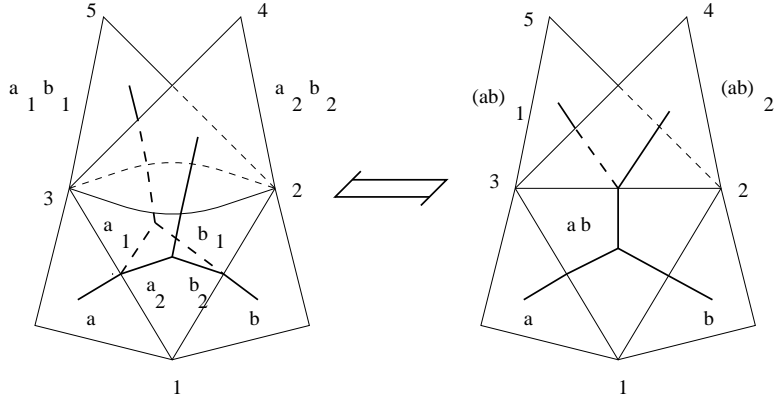}
}
\end{center}
\caption{The  cone move in dimension $3$}
\label{3dcone}
\end{figure}

Let us now explain the relation of this move to the compatibility condition
verbally.
In the left hand side of Fig. \ref{3dcone} 
there are distinct and parallel
triangular faces sharing the edge $(12)$ and $(13)$; 
these triangles  have different
edges connecting the vertices $2$ and $3$.
One of these is shared by the face $(234)$ while the other is shared 
by the face $(235)$.

The parallel faces (123) and (123)' are collapsed to a single face to
obtain the right hand side of Fig.~\ref{3dcone}.
Now there is a single face with edges $(12)$, $(23)$, and 
$(31)$, and the edge $(23)$ is shared by three faces 
$(123)$, $(234)$, and $(235)$.

The thick segments indicate part of the dual graph.
Each segment is labeled by Hopf algebra elements.
Reading from bottom to top, one sees that 
the graphs represent maps 
from $A \ox A$ to itself. 
The left-hand-side of the figure 
represents
$$(m \ox m) \circ (1 \ox P \ox 1) \circ (\Delta \ox \Delta) (a \ox b)
   $$ 
$$= (m \ox m) \circ (1 \ox P \ox 1)  (\Delta a \ox \Delta b ) $$
$$ = (m \ox m) ( a_1 \ox b_1 \ox a_2 \ox b_2 )
= (a_1 b_1 ) \ox ( a_2  b_2 ) $$
while the right-hand-side represents 
$$ \Delta \circ m (a \ox b)
= \Delta (ab) 
= (ab)_1 \ox (ab)_2 
$$
and these are equal by the consistency condition between
multiplication and comultiplication.
This shows that the Hopf algebra structure gives 
solutions to the equation
corresponding  
to the cone move.

That the partition function in this case does not depend on the 
choice of triangulation   
is proved by showing that 
the Pachner moves follow from the cone move and other 
singular moves. 
 Figure~\ref{chfs} explains why the $(2 \rightleftharpoons 3)$-move 
follows from singular moves (this figure                    
is basically the same as a Figure
in \cite{ChFS}).

Let us explain the figure. 
The first polyhedron is the right-hand-side 
of the $(2 \rightleftharpoons 3)$-move.
There are three internal faces and three tetrahedra. 
Perform the cone move along edge $(25)$ thereby duplicating face 
$(125)$. Internally, we have face $(125)$ glued to face $(235)$ along 
edge $(25)$ and face $(125)'$ glued to face $(245)$ along edge $(25)'$.
These faces are depicted in the second polyhedron.
By associativity these faces can be replaced by four faces 
parallel to four faces on the boundary;
$(123)$, $(135)$, $(124)$, $(145)$. This is 
the configuration in the 
third polyhedron.
Then there are two $3$-cells bounded by these parallel faces.
Collapse these cells and push the internal faces onto the boundary
(this is done by singular moves).
The result is the fourth polytope which now is a single polytope
without any internal faces.
This is the middle stage in the sense that we have proved that 
the 
right-hand-side
of the $(2 \rightleftharpoons 3)$-move 
is in fact equivalent to this polytope. 

Now introduce a pair of internal faces parallel to the faces 
$(135)$ and $(145)$ to get the fifth polytope 
(the left bottom one).
Perform  associativity again to change it to 
a pair of faces  
$(134)$ and $(345)$ to get the 
sixth
polytope.
Perform a cone move along the pair of faces with vertices 
$(345).$ These faces share edges 
$(35)$ 
and $(45);$ edge 
$(34)$ 
is 
duplicated.) 
The last picture which is the left-hand-side of 
the $(2 \rightleftharpoons 3)$-move.

In summary, 
we perform cone moves and collapse  
some 
$3$-cells to the boundary and prove that both sides of
the Pachner move is in fact equivalent to 
the polyhedral $3$-cell without internal faces.
We give a generalization of this Theorem to dimension $4$
in 
Lemmas~\ref{pac33}, \ref{pac24}, and \ref{pac15}.
 
\begin{figure}
\begin{center}
\mbox{
\epsfxsize=5in
\epsfbox{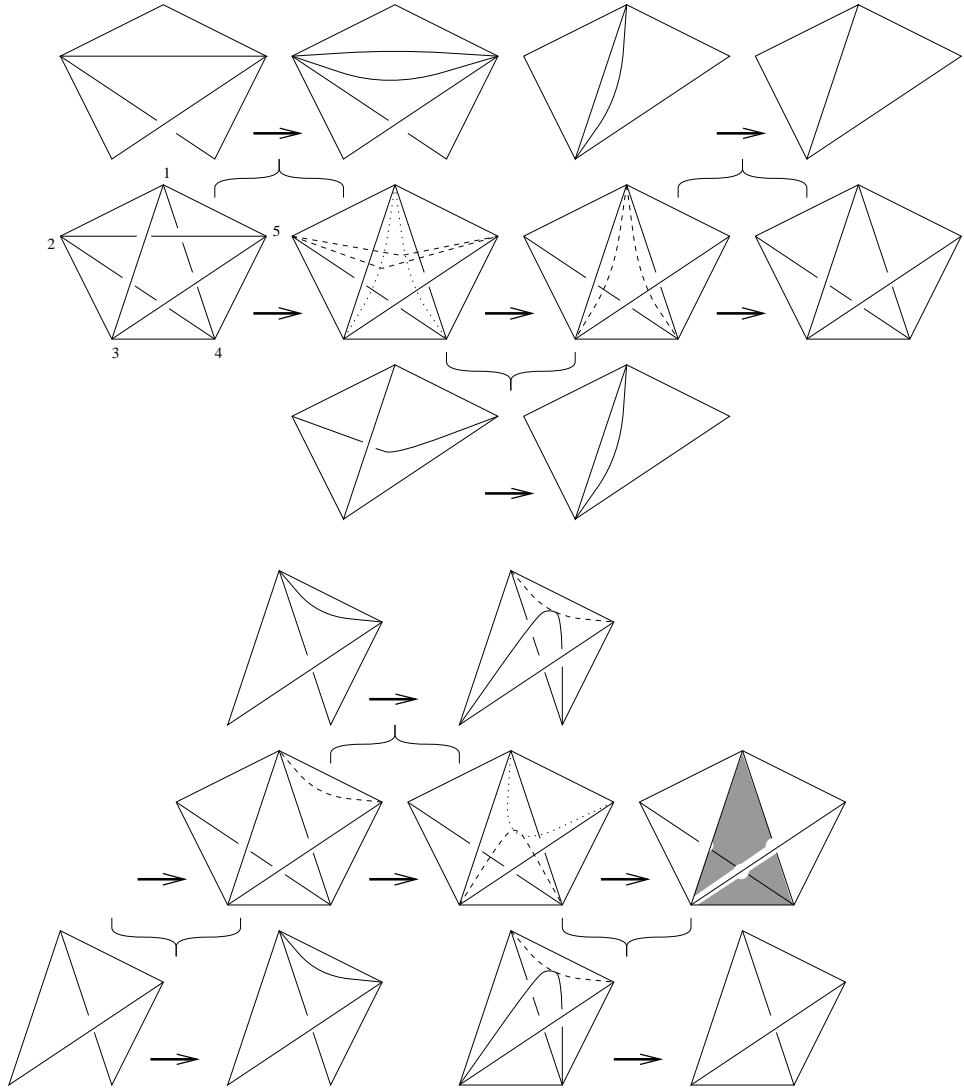}
}
\end{center}
\caption{A Pachner move follows from cone moves}
\label{chfs}
\end{figure}
 
}\end{sect}

\begin{figure}
\begin{center}
\mbox{
\epsfxsize=4in
\epsfbox{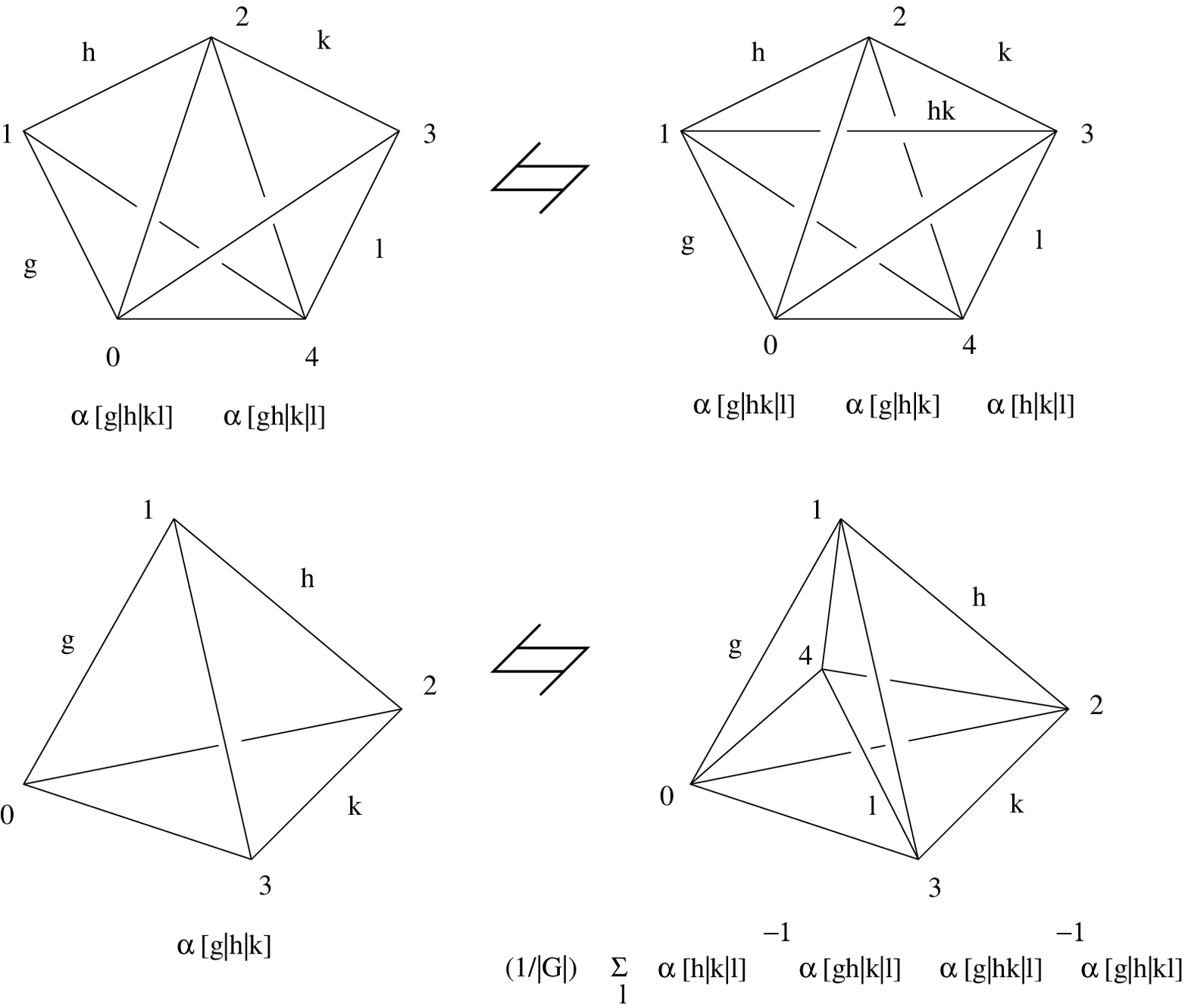}
}
\end{center}
\caption{Pachner moves and the $3$-cocycle condition}
\label{3cocycle}
\end{figure}

\begin{sect}{\bf Dijkgraaf-Witten invariants. \/}\addcontentsline{toc}{subsection}{Dijkgraaf-Witten invariants }
{\rm 
We review the Dijkgraaf-Witten invariants for 3-dimensional manifolds.
In \cite{DW} Dijkgraaf and Witten  gave a combinatorial
definition for Chern-Simons invariants with finite gauge
groups using $3$-cocycles of the group cohomology.
We follow Wakui's description \cite{Wakui}
except we use the Pachner moves. See \cite{Wakui}
for more detailed treatments.

Let $T$ be a triangulation of an oriented closed
$3$-manifold $M$,
with $a$ vertices and $n$ tetrahedra.
Give an ordering to the set of vertices. Let $G$ be a finite group.
Let $\phi : $ \mbox{ $\{$ oriented edges $\} \rightarrow G$ } 
be a map such 
that

(1) for any triangle with vertices $v_0, v_1, v_2$ of $T$,
$\phi(\langle v_0, v_2 \rangle )=\phi( \langle v_1, v_2 \rangle )
 \phi( \langle v_0, v_1 \rangle )$,
where $ \langle v_i, v_j \rangle $ denotes the oriented edge,
and

(2) $\phi(-e) = \phi(e)^{-1}$.

Let $\alpha :G \times G \times G \rightarrow A$,
$(g,h,k) \mapsto \alpha [g|h|k] \in A$, be a $3$-cocycle
with a multiplicative abelian group $A$.
The $3$-cocycle condition is
$$\alpha [h|k|l] \alpha [gh|k|l]^{-1} \alpha [g|hk|l]
\alpha [g|h|kl]^{-1} \alpha [g|h|k] =1. $$
Then the Dijkgraaf-Witten invariant is defined by
$$ Z_M = \frac{1}{|G|^a} \sum_{\phi} \Pi_{i=1}^{n} W( \sigma, \phi )
^{\epsilon_i}. $$
Here $a$ denotes the number of the vertices of the
given triangulation,                     
$W(\sigma, \phi )= \alpha [g|h|k] $ where
$\phi( \langle v_0, v_1 \rangle )=g$, $\phi( \langle v_1, v_2 \rangle  )=h$,
 $\phi( \langle v_2, v_3 \rangle )=k$,
for the tetrahedron $\sigma = |v_0 v_1 v_2 v_3|$ with
the ordering $v_0<v_1<v_2<v_3$,
and  
$\epsilon=\pm 1 $ according to 
whether
or not the orientation
of $\sigma$ with respect to the vertex ordering matches the
orientation of $M$.

 Then one  checks the invariance of this state sum 
under Pachner moves, see Figure~\ref{3cocycle}.

} \end{sect}

\begin{sect}{\bf Summary: Going up 
dimensions.\/}
\addcontentsline{toc}{subsection}{Going up  dimensions}
{\rm
As we reviewed the invariants in dimensions $2$ and $3$,
there are two ways to go up the dimension from $2$ to $3$.
One way is to consider the algebras formed by representations
of quantum groups as in Turaev-Viro invariants.
In this case algebra elements are regarded as vector spaces
(representations) and algebras are replaced by
such categories. This process is called  categorification.
The second way is to include a comultiplication in addition to
the multiplication of an algebra as in Hopf algebra invariants
and consider bialgebras (in fact Hopf algebras) instead of algebras.

Crane and Frenkel defined invariants in dimension $4$
using these ideas. We reach at the idea of the 
algebraic  
structure
called {\it Hopf categories} either by 
(1) categorifying Hopf algebras, or (2) including 
comultiplications to categories of representations.
The following chart represent this idea.

\[
\begin{array}{cccc}
2D & & \fbox{Associative algebras} &  \\[5mm]
 & \parbox{1in}{Categorification} \hfill \swarrow & & \searrow 
\hfill \parbox{1in}{Adding a comultiplication} \\[5mm]
3D & \fbox{Turaev-Viro invariants} & &
\fbox{Hopf algebra invariants} \\[5mm]
 &  \parbox{1in}{Adding a comultiplication} \hfill \searrow & & 
\swarrow \hfill \parbox{1in}{Categorification} \\[5mm]
 4D & & \fbox{Crane-Frenkel invariants}
\end{array}
\]

In the following sections we follow this idea to define 
invariants in dimension $4$.

We also point out here that the theories reviewed above 
have remarkable features in that they have 
direct relations between algebraic structures and 
triangulations via diagrams (trivalent planar graphs, or 
{\it spin networks}).
On the one hand such diagrams appear as dual complexes
through movie 
descriptions of duals of triangulations, and on the other hand
they appear as diagrammatic representations of maps in algebras. 
In the following sections we explore such relations 
and actually utilize diagrams to prove well-defined-ness
of the invariants proposed by Crane and Frenkel.
}\end{sect}

\section{Pachner Moves in dimension $4$ } \label{Pac4sec}  

In  Section~\ref{pac23} we reviewed the Pachner moves for triangulations
in dimensions $2$ and $3$ and their relations to 
associativity of algebras.
In this section, we describe Pachner moves in dimension $4$.
Relations of these moves to the Stasheff polytope 
was discussed in \cite{CKS}.

In general, an $n$-dimensional Pachner move of type $(i
\rightleftharpoons j)$,
where $i+j=n+2$, is obtained by decomposing the (spherical)
boundary of an
$(n+1)$-simplex into the union of two $n$-balls such that one of the
balls is the union of $i \  \ $ $n$-simplices, the other ball is
the union of $j \ \ $ $n$-simplices, and the intersection of these
balls is an $(n-1)$-sphere. By labeling the vertices of the
$(n+1)$-simplex these moves are easily expressed. For example,
the table below indicates the lower dimensional Pachner moves:

\vspace{1cm}

\begin{tabular}{|l|c|c|} \hline
$n=1$ & $(1 \rightleftharpoons 2)$ & $(01) \rightleftharpoons (02) \cup
(12)$ \\
 \hline
$n=2$ & $(1 \rightleftharpoons 3)$ &  $(012) \rightleftharpoons (013) \cup
(023)
 \cup (123)$
\\
    & $(2 \rightleftharpoons 2)$ &  $(012) \cup (023) \rightleftharpoons
(013) \
\cup  (123)$ \\
\hline
$n=3$ & $(1 \rightleftharpoons 4)$ &  $(0123) \rightleftharpoons (0134)
\cup (02
34) \cup
(1234)$ \\
    & $(2 \rightleftharpoons 3)$ &  $(0123) \cup (1234) \rightleftharpoons
(0124
) \cup
(0134)  \cup (0234)$ \\ \hline
$n=4$ & $(1 \rightleftharpoons 5)$ &  $(01234) \rightleftharpoons (01235)
\cup (
01245)
\cup (01345) \cup (02345) \cup
(12345)$  \\
    & $(2 \rightleftharpoons 4)$ &  $(01234) \cup (01235)
\rightleftharpoons
(12345) \cup  (01245)  \cup (01345)  \cup (02345)$ \\
    & $(3 \rightleftharpoons 3)$ &  $(01234) \cup (01245) \cup (02345)
\rightleftharpoons
      (01235) \cup  (01345)  \cup (12345)$  \\ \hline
\end{tabular}

\vspace{1cm}

The relationship between the general Pachner move and the higher order 
associativity relations are explained in \cite{CKS}.
Next we turn to a more explicit description of the $4$-dimensional 
Pachner Moves.

\begin{sect}{\bf  4-dimensional Pachner moves.\/}
\addcontentsline{toc}{subsection}{ 4-dimensional Pachner moves}
{\rm
In this section we explain the $4$-dimensional
Pachner moves.
One side of  a 4-dimesional  Pachner move is the  
union of $4$-faces
of a $5$-simplex (homeomorphic to a $4$-ball), and
the other side of the move is the union of the rest of
$4$-faces.

\begin{figure}
\begin{center}
\mbox{
\epsfxsize=5in
\epsfbox{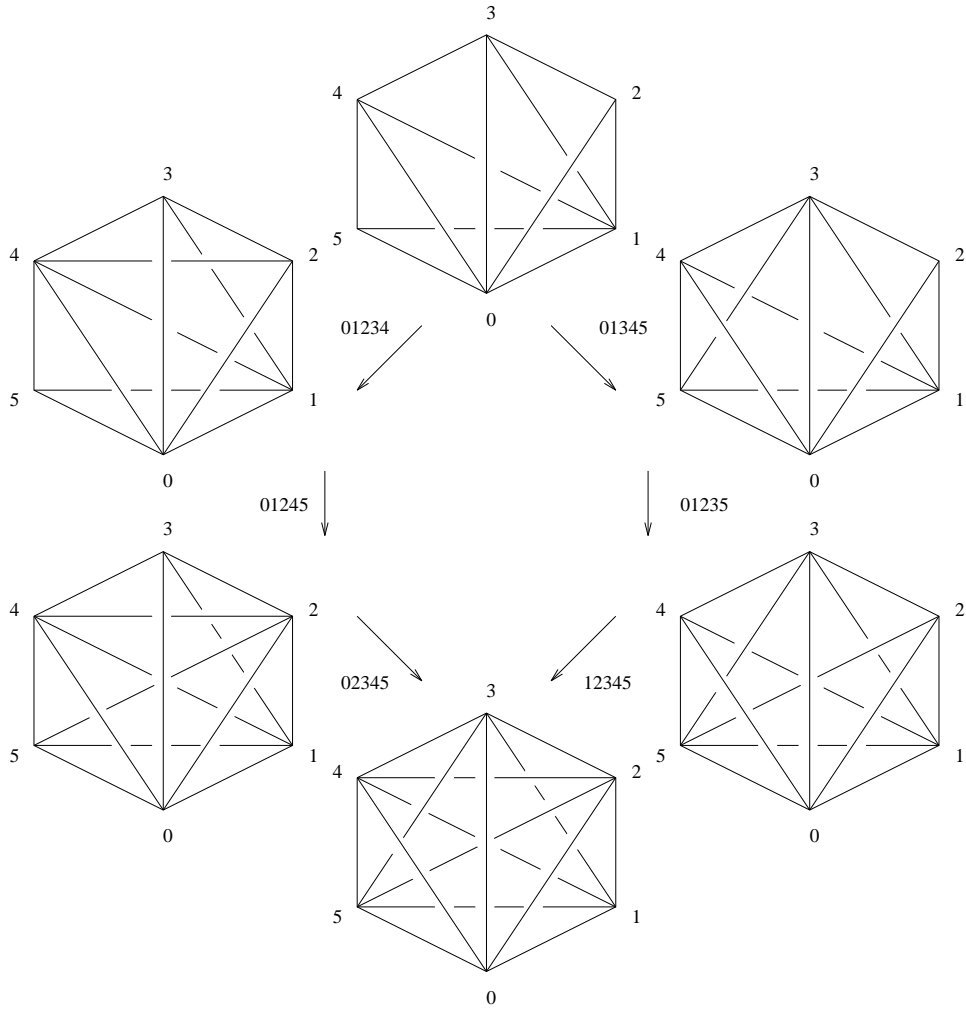}
}
\end{center}
\caption{ 4-dimensional Pachner move I}
\label{pachner1}
\end{figure}

\begin{figure}
\begin{center}
\mbox{
\epsfxsize=5in
\epsfbox{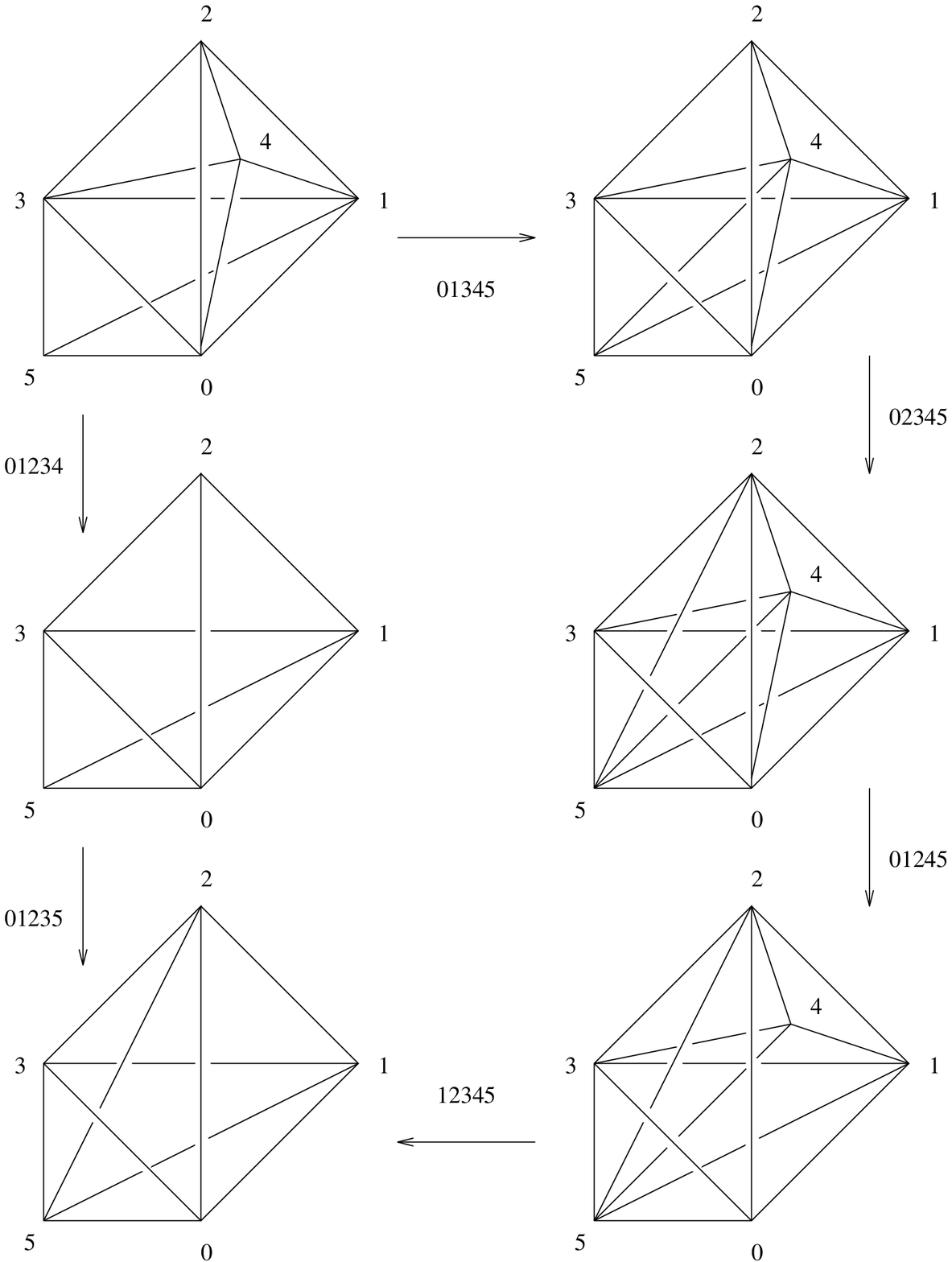}
}
\end{center}
\caption{ 4-dimensional Pachner move II}
\label{pachner2}
\end{figure}

\begin{figure}
\begin{center}
\mbox{
\epsfxsize=5in
\epsfbox{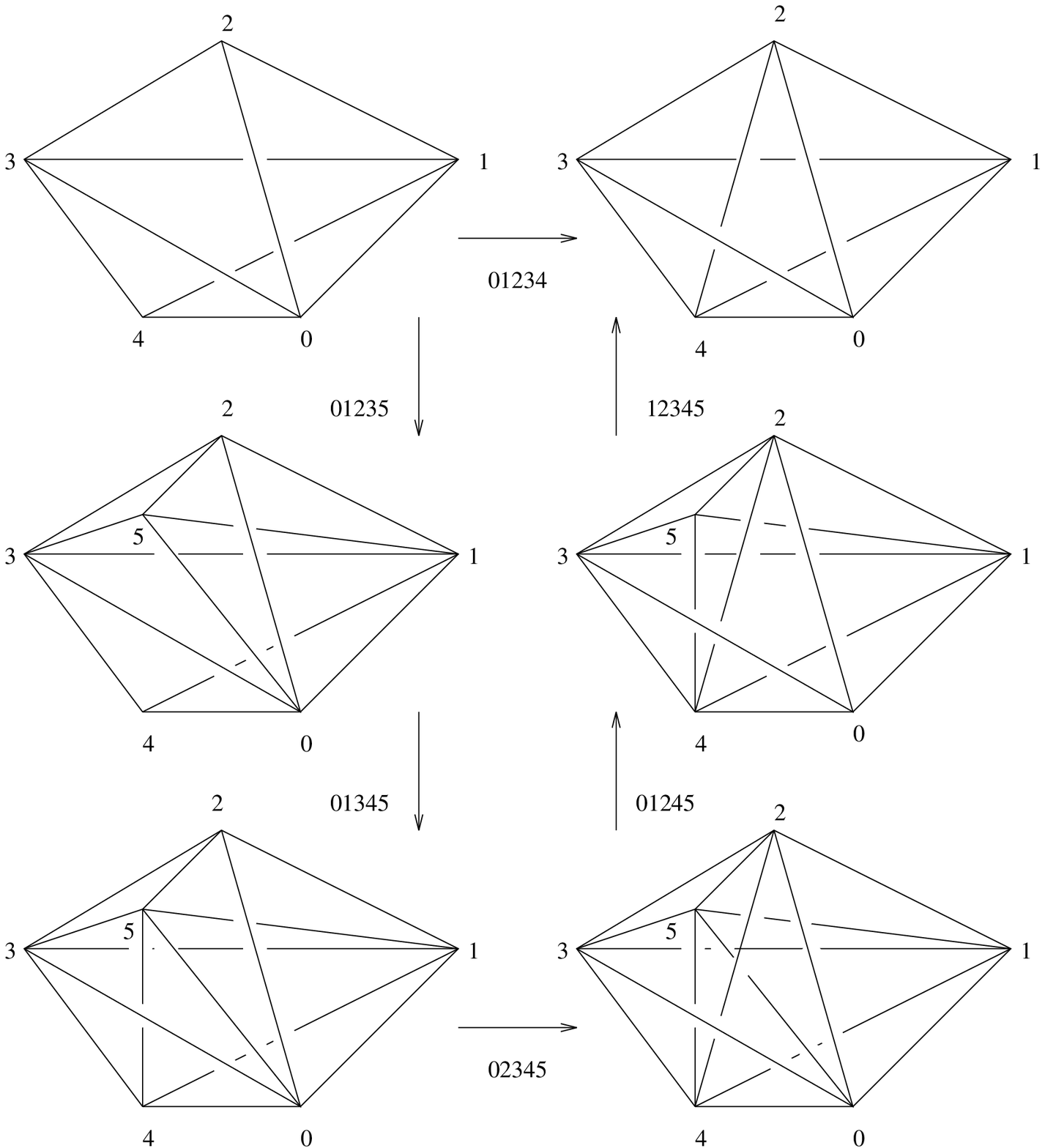}
}
\end{center}
\caption{ 4-dimensional Pachner move III}
\label{pachner3}
\end{figure}

In Figures~\ref{pachner1} \ref{pachner2}, and \ref{pachner3}
the $(3 \rightleftharpoons 3)$-move, 
$(2 \rightleftharpoons 4)$-move,
and $(1 \rightleftharpoons 5)$-move are  depicted, respectively.
Recall here that each 3-dimensional Pachner move represents a $4$-simplex.
Therefore the 3-dimensional Pachner move depicted in
the top left of Figure~\ref{pachner1}, 
the move represented by an arrow
labeled $(01234)$, represents the $4$-simplex with vertices
$0$, $1$, $2$, $3$ and $4$.
Then the left-hand side of Fig.~\ref{pachner1} represents
the union of three $4$-simplices
$(01234)\cup (01245) \cup (02345)$.
Similarly, the right-hand side of Fig.~\ref{pachner1} represents the
union of the three $4$-simplices
$(01345)\cup (01235) \cup (12345)$.}\end{sect}

\begin{sect}{\bf Singular moves.\/} \label{mostex}
\addcontentsline{toc}{subsection}{Singular moves}
{\rm 
In dimension 4, the Pachner moves can be decomposed as singular moves and 
lower dimensional   
moves. Here we define a 4-dimensional 
singular moves (called cone, pillow, taco moves) 
and show how the Pachner
moves follow. This material was discussed in \cite{CKS}.

}\end{sect}

\begin{figure}
\begin{center}
\mbox{
\epsfxsize= 6in
\epsfbox{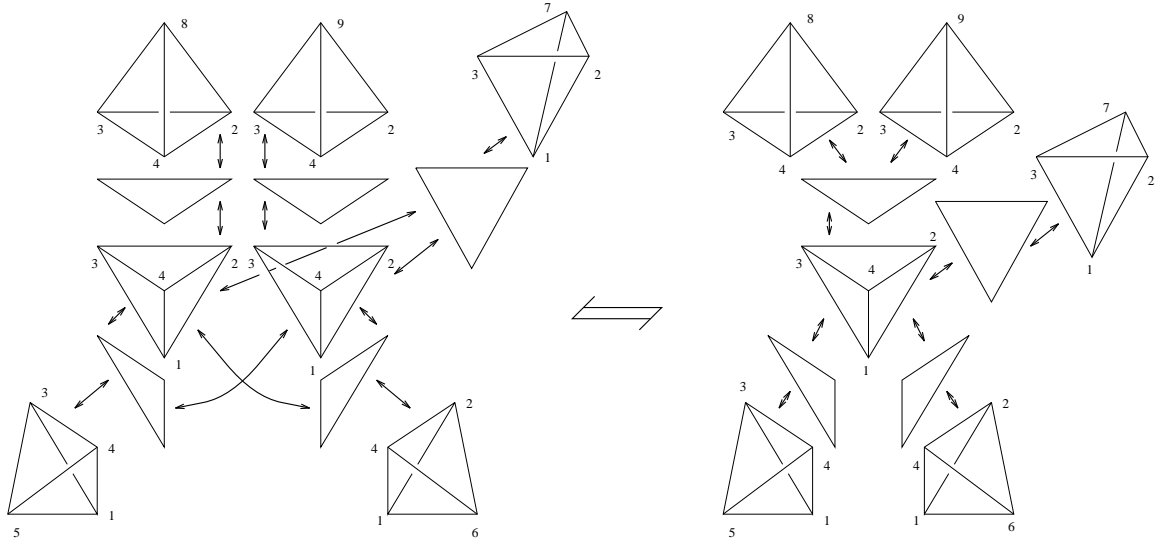}
}
\end{center}
\caption{ 4-dimensional cone move}
\label{4dcone}
\end{figure}

\begin{subsect} {\bf Definition (cone move).} {\rm
The  {\it cone move} for $CW$-complexes 
for $4$-manifolds is
defined as follows.

Suppose there is a pair of tetrahedra
$(1234)_1$ and $(1234)_2$ 
that share 
the same faces
$(123)$, $(124)$ and $(134)$, but have different
faces 
 $(234)_1$ and $(234)_2$, such that
(1) $(234)_1$ and $(234)_2$ bound a 3-ball $B$ in
the 4-manifold,  (2) the union of $B$ , $(1234)_1$
and $(1234)_2$ is diffeomorphic to the $3$-sphere
bounding a $4$-ball $W$ in the 4-manifold.

The situation is depicted in Figure~\ref{4dcone} 
which we now explain.
The 
left-hand-side 
of the Figure has two copies of tetrahedra
with vertices $1$, $2$, $3$, and $4$.
They share the same faces $(123)$, $(124)$, and $(134)$
but have two different faces with vertices $2$, $3$, and $4$.

\begin{center}
\begin{tabular}{||c|c|c||} \hline  \hline
Triangle & is a face of & tetrahedron \\ \hline \hline
$(234)_1$ & $\subset$ & $(2348)$ \\ \hline
$(123)_2$ & $\subset$ & $(2349)$ \\ \hline
$(123)$    & $\subset$ & \begin{tabular}{c} $(1234)_1$ \\ $(1234)_2$ \\ $(1237)$ \end{tabular} \\ \hline
 $(124)$   & $\subset$ & \begin{tabular}{c}  $(1234)_1$ \\ $(1234)_2$ \\$(1246)$
\end{tabular} \\ \hline 
$(134)$   & $\subset$ & \begin{tabular}{c}  $(1234)_1$ \\ $(1234)_2$ \\$(1345)$ 
\end{tabular} \\ \hline \hline
\end{tabular} \end{center}

Collapse these two tetrahedra to a single tetrahedra to get
the 
right-hand-side 
of the Figure. Now we have a single tetrahedron
with vertices $1$, $2$, $3$, and $4$.
The face $(234)$ now is shared by three tetrahedra
$(1234)$, $(2348)$, and $(2349)$
while three faces $(123)$, $(124)$, and $(134)$ are
shared by two tetrahedra.
} \end{subsect}

\begin{subsect} {\bf Definition (taco move).\/} {\rm 
Suppose we have a $CW$-complex 
such that there is a pair of  tetrahedra 
$(1234)_1$ and $(1234)_2$ that share two faces 
$(123)$ and $(124)$ but have different faces 
$(134)_1$, $(134)_2$ 
and $(234)_1$, $(234)_2$ (of $(1234)_1$, $(1234)_2$ respectively).
Suppose further that $(134)_1$, $(134)_2$, $(234)_1$, and $(234)_2$
together bound a $3$-cell $B$ and $(1234)_1$, $(1234)_2$, and $B$
bounds a $4$-cell.
Then collapse this $4$-cell to get a single tetrahedron $(1234)$.
As a result $(134)_1$ (resp. $(234)_1$) and $(134)_2$ (resp. $(234)_2$)
are identified. 
This move is called the {\it taco} move.
} \end{subsect}

\begin{subsect} {\bf Definition (pillow move).} {\rm
Suppose we have a $CW$-complex 
such that there is a pair of tetrahedra sharing all four faces
cobounding a $4$-cell.
Then collapse these tetrahedra to a single tetrahedron.
This move is called the {\it pillow move}.
} \end{subsect}

\begin{figure}
\begin{center}
\mbox{
\epsfxsize= 3in
\epsfbox{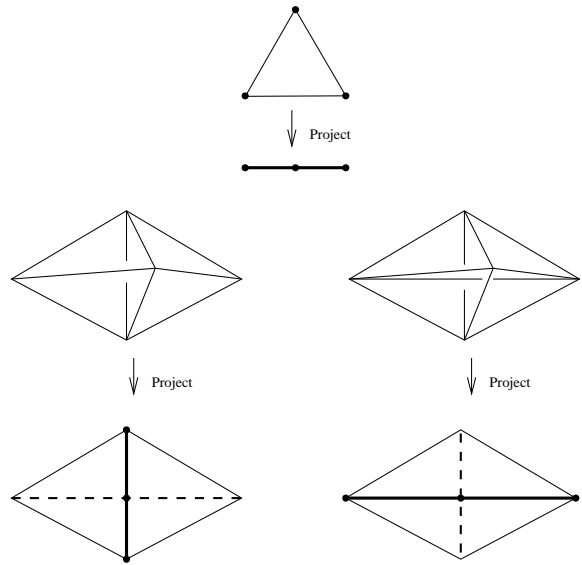}
}
\end{center}
\caption{Projecting the (2,3)-move}
\label{project32}
\end{figure}

\begin{figure}
\begin{center}
\mbox{
\epsfxsize= 4in
\epsfbox{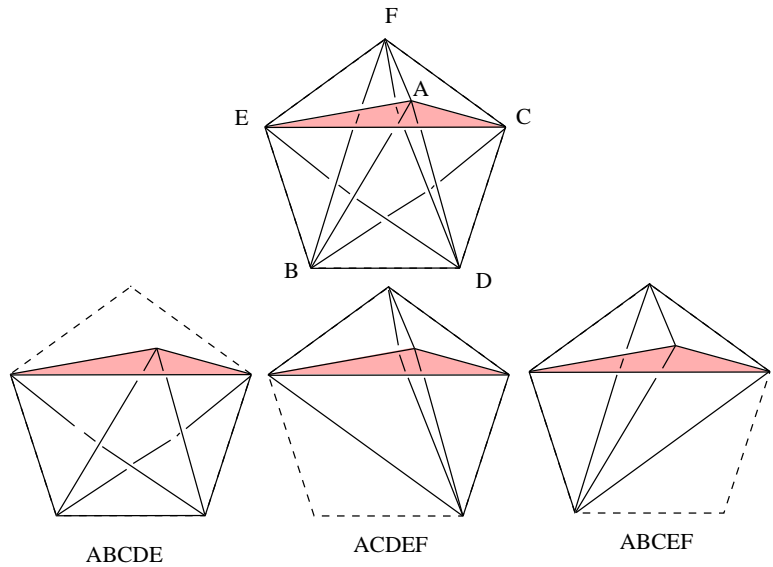}
}
\end{center}
\caption{The left-hand-side of the (3,3)-move}
\label{tacolem1}
\end{figure}

\begin{figure}
\begin{center}
\mbox{
\epsfxsize= 4in
\epsfbox{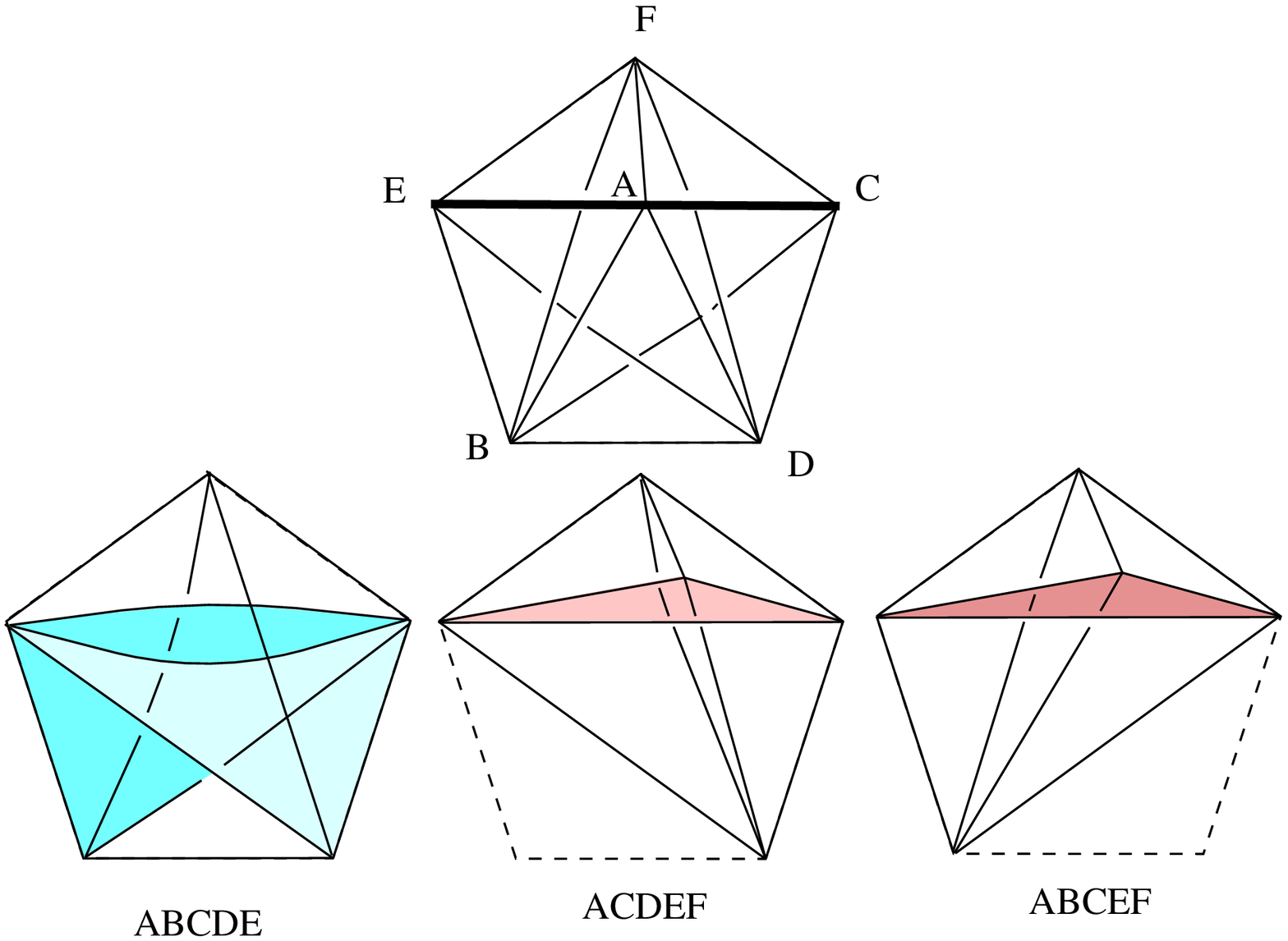}
}
\end{center}
\caption{Splitting the tetrahedron $(ACEF)$ via a cone move}
\label{tacolem2}
\end{figure}

\begin{figure}
\begin{center}
\mbox{
\epsfxsize= 4in
\epsfbox{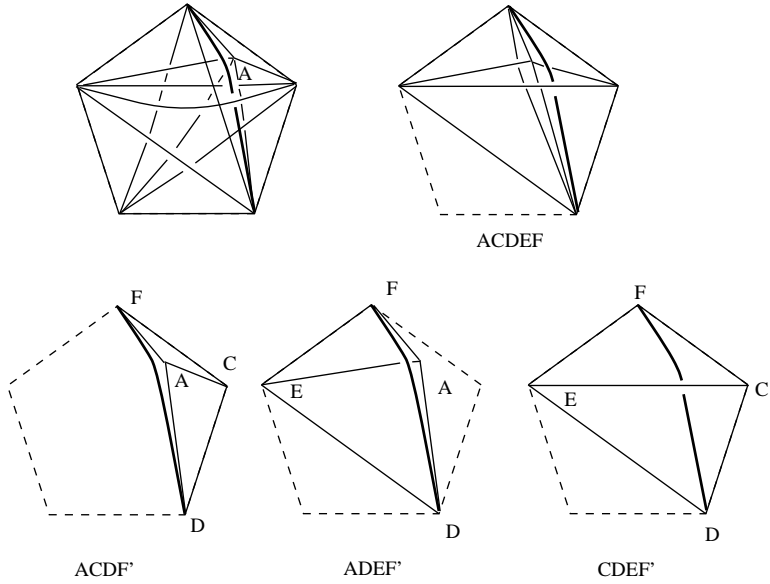}
}
\end{center}
\caption{Performing a $(2,3)$-move to $(ACDE) \cup (ACEF)$}
\label{tacolem3}
\end{figure}

\begin{figure}
\begin{center}
\mbox{
\epsfxsize= 4in
\epsfbox{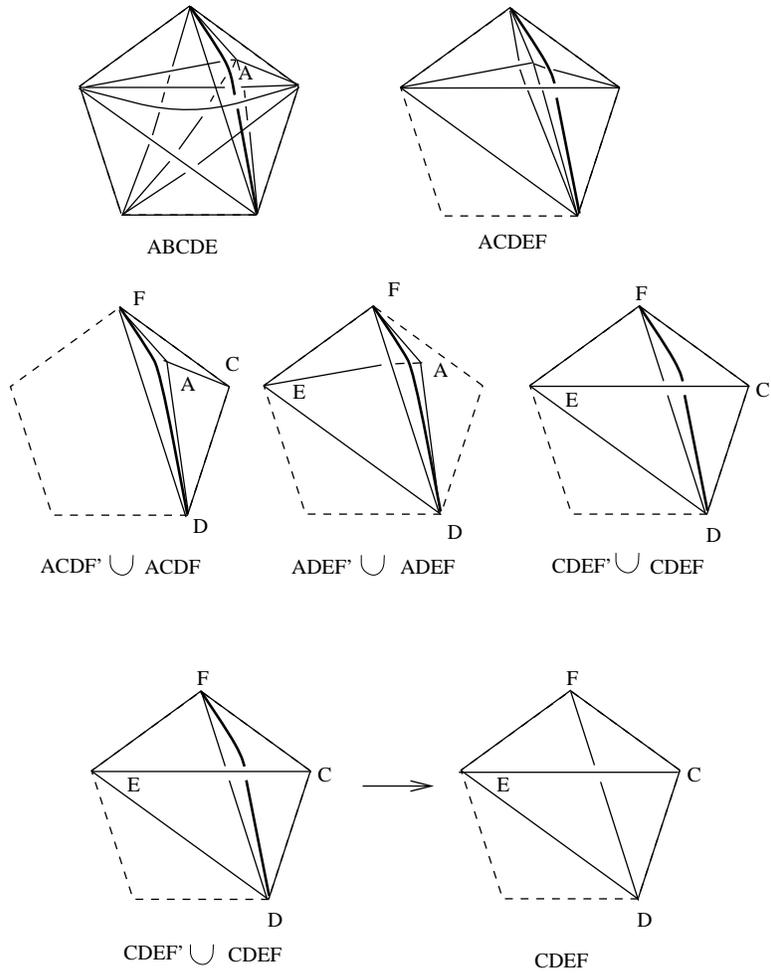}
}
\end{center}
\caption{Performing a taco move to the pair $(CDEF)$ and $(CDEF)'$
}
\label{tacolem4}
\end{figure}

\begin{figure}
\begin{center}
\mbox{
\epsfxsize= 4in
\epsfbox{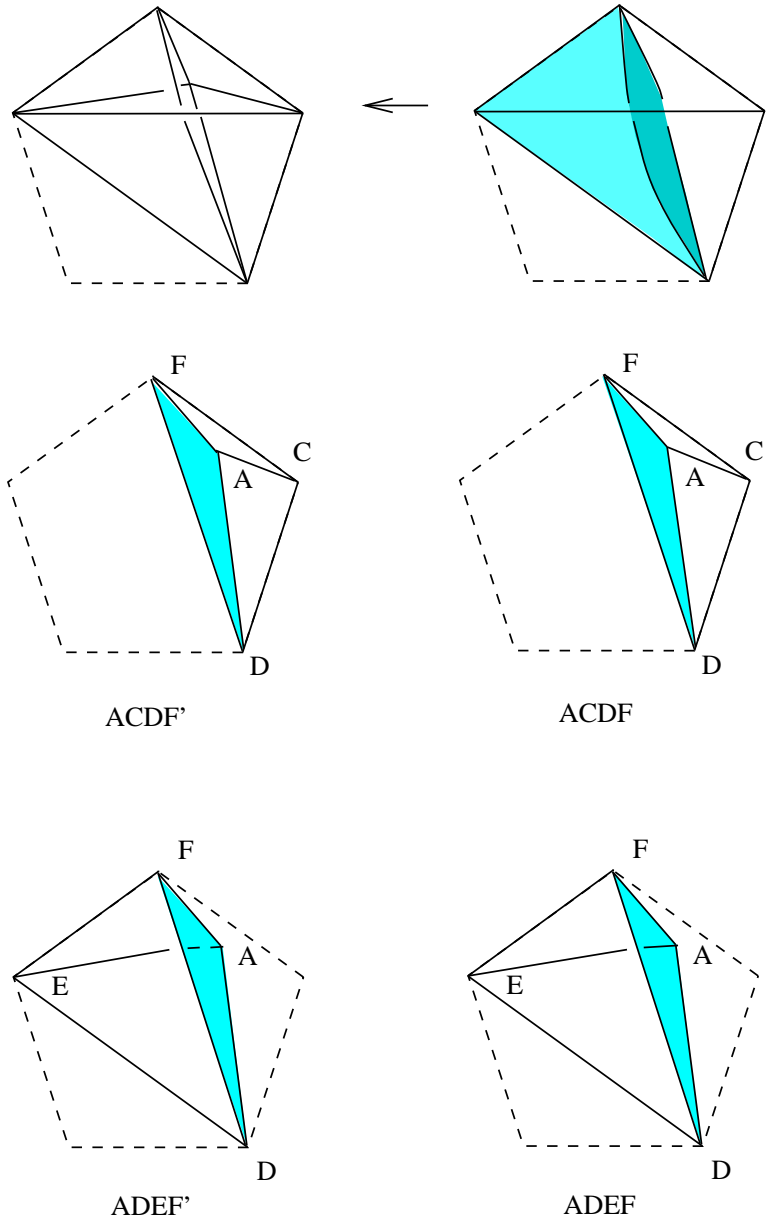}
}
\end{center}
\caption{Performing a cone move to the pair $(ADEF)$ and $(ADEF)'$
}
\label{tacolem5}
\end{figure}

\begin{figure}
\begin{center}
\mbox{
\epsfxsize= 4in
\epsfbox{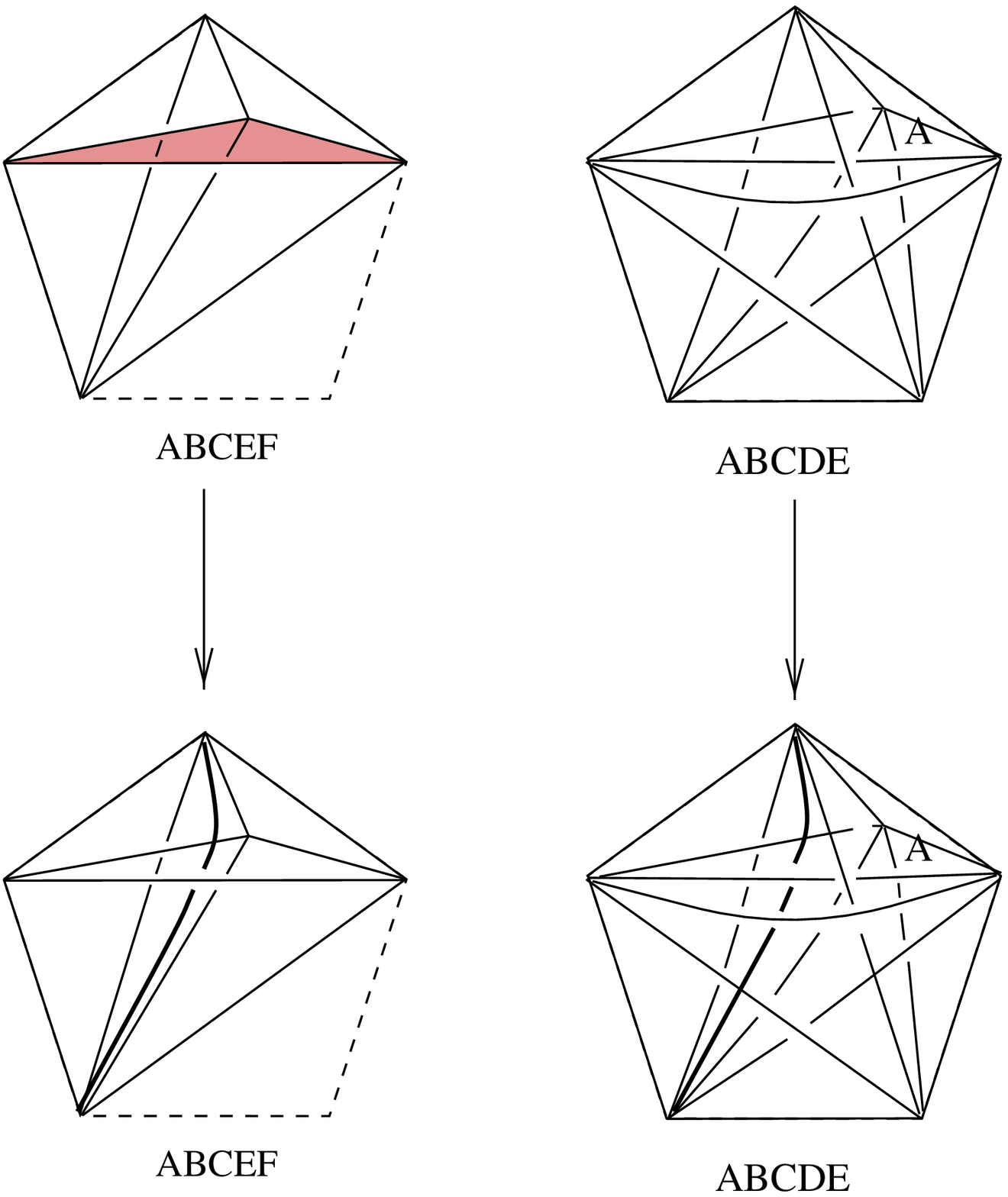}
}
\end{center}
\caption{Performing a $(2,3)$-move to the pair $(ABCE)$ and $(ACEF)_1$
}
\label{tacolem6}
\end{figure}

\begin{figure}
\begin{center}
\mbox{
\epsfxsize= 4in
\epsfbox{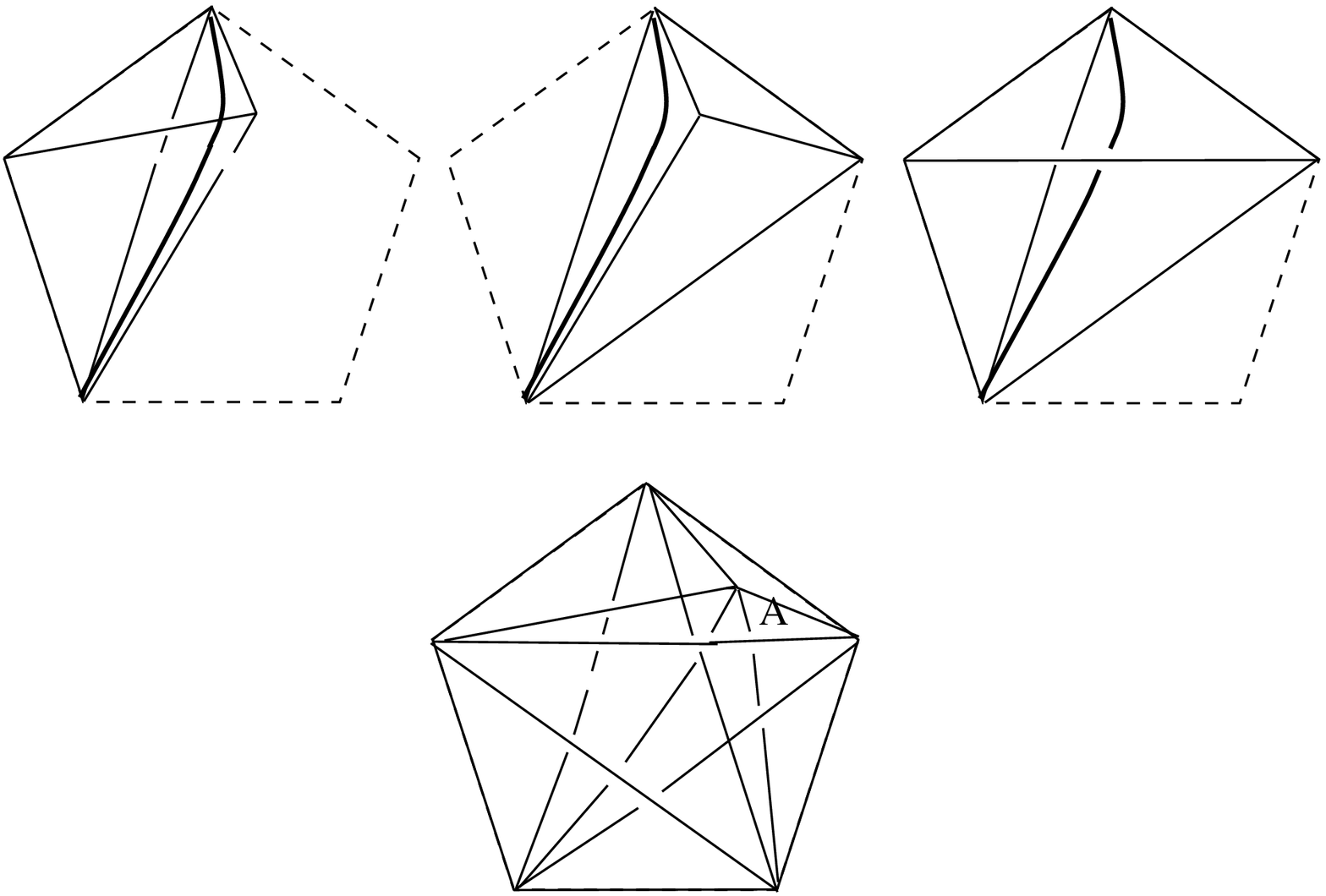}
}
\end{center}
\caption{Collapsing to a single polytope}
\label{tacolem7}
\end{figure}

\begin{figure}
\begin{center}
\mbox{
\epsfxsize= 4in
\epsfbox{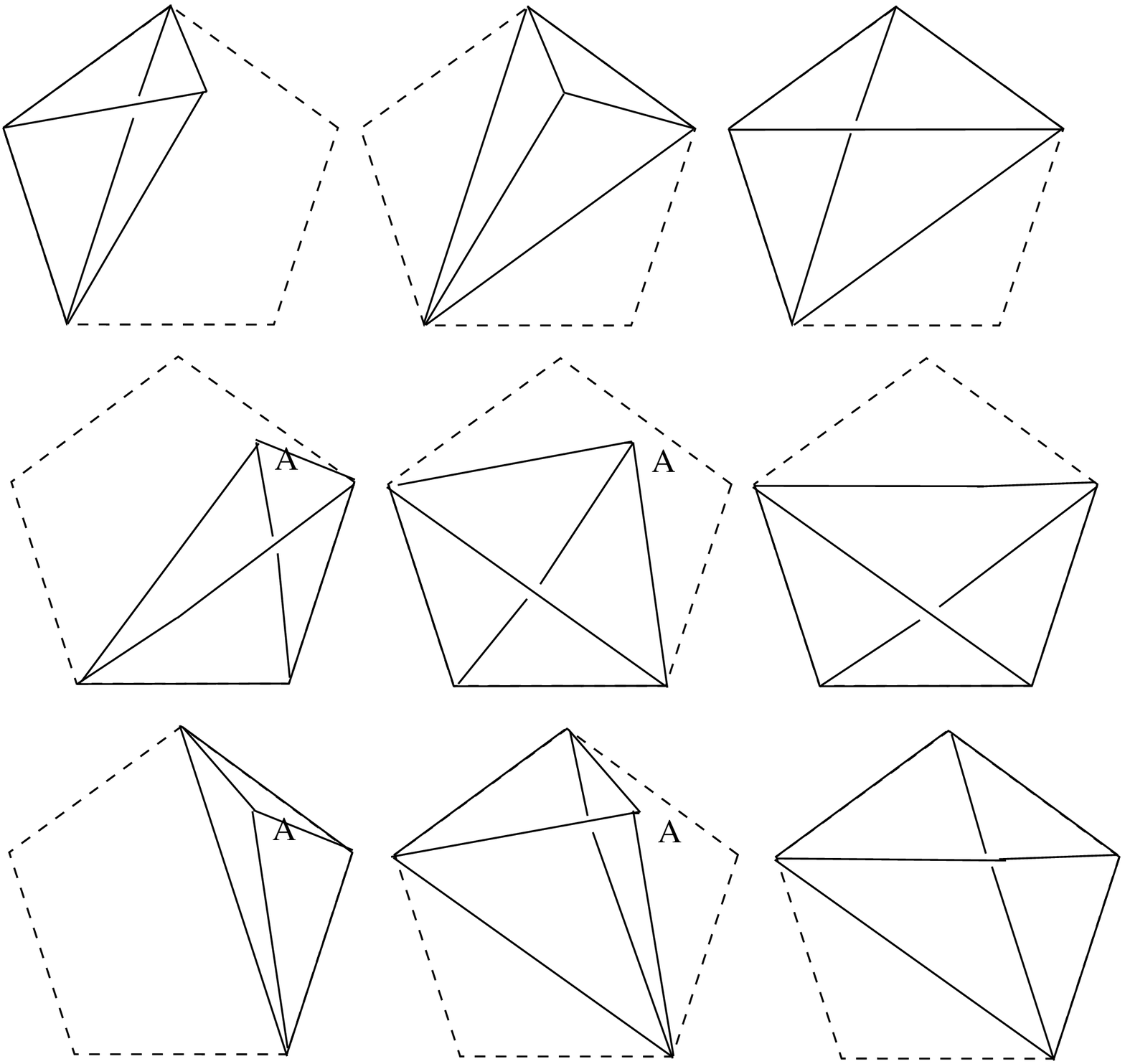}
}
\end{center}
\caption{The tetrahedral faces of the middle stage}
\label{tacolem8}
\end{figure}

\begin{subsect} {\bf Lemma.} \label{pac33} 
The 
$(3 \rightleftharpoons 3)$
Pachner move is described as a sequence of
cone moves, pillow moves, taco moves,              
 and   $3$-dimensional Pachner moves.
\end{subsect}
  {\it Proof.}
The proof can be facilitated by following the Figures~\ref{project32} through
\ref{tacolem8}. Figure~\ref{project32} is a preliminary sketch that
 indicates in dimension 3 the methods of the subsequent figures.
 It illustrates that the  $( 2 \rightleftharpoons 3)$-move in dimension 3 can be interpreted
 in terms of the  $( 2 \rightleftharpoons 2)$-move via a non-generic projection.
 The thick vertical line on the left-hand-side of the figure
 is the projection of the triangle along which the two tetrahedra are glued.
 The thick horizontal line on the right is the projection of one of the three
 triangles that are introduced on the right-hand-side of the move.
 The other two triangles project to fill the lower right quadrilateral. 
The dotted lines indicate that some edges in the figure will project to 
these lines. Some information is lost during the projection process,
 but at worst, the projected figures serve as a schematic diagram
 of the actual situation.

In Fig.~\ref{tacolem1} the union of the three 4-simplices 
$(ABCDE)$, $(ACDEF)$, and $(ABCEF)$ is illustrated; 
these share the triangle
$(ACE)$ which is shaded in figure.
 The union forms the left-hand-side of the 
of the $(3 \rightleftharpoons 3)$-move. 
Let $P$ denote this union. 
 In the top of Fig.~\ref{tacolem2}, the triangle $(ACE)$
 has been projected to the thick line  $(EAC)$. At the bottom of Fig.~\ref{tacolem2}, the 4-simplex $(ACEF)$ has been split into simplices
 $(ACEF)_1$ and $(ACEF)_2$ by a cone move. The cone move is illustrated
 in this projection, and the schematic resembles the cone move
 in dimension 3 as is seen on the bottom left of the figure. 
Thus $(ACEF)_1$ and $(ACEF)_2$ share the same faces
$(ACF)$, $(AEF)$ and $(CEF)$ but have different
faces $(ACE)_1$ and $(ACE)_2$.
The face $(ACE)_1$ is shared with $(ABCE)$
and the face $(ACE)_2$ is shared with $(ACDE)$ respectively.

After the splitting,     
 $P$ consists of three $4$-polytopes,
$\tau^1_j$, $j= 1,2,3$.
Here the polytope $\tau^1_1$ is bounded by
tetrahedra $(ABCD)$, $(ABDE)$,
$(ABCE)$, $(ACDE)$,
$(BCDE)$, $(ACEF)_1$, and $(ACEF)_2$.
The polytope $\tau^1_2$ is bounded by
tetrahedra $(ABCE)$, $(ABEF)$,
$(ACEF)_1$,
$(ABCF)$, and $(BCEF)$.
The polytope $\tau^1_3$ is bounded by
tetrahedra $(ACDE)$, $(ACEF)_2$ 
$(ACDF)$, $(ADEF)$, and $(CDEF)$.
The polytope $\tau^1_1$ corresponds to $(ABCDE)$ and
 it is illustrated on the left bottom of Fig.~\ref{tacolem2}
 (labeled $(ABCDE)$ to indicate the correspondence).
 On the bottom right  of the figure, we see 
  the polytope $\tau^1_2$ labeled $(ABCEF)$. In
 the bottom center of the figure the polytope $\tau^1_3$
 labeled $(ACDEF)$ to indicate its antecedent.
 Our first work will be on $\tau^1_1$ and $\tau^1_3$.

Next perform a   Pachner move to
the pair of tetrahedra $(ACDE) \cup (ACEF)_2$
sharing the face $(ACE)_2$.
Note that these two tetrahedra are shared by
$\tau^1_1$ and $\tau^1_3$ so that the   Pachner
move we perform does not affect $\tau^1_2$.
Thus we get three $4$-cells $\tau^2_j$, $j= 1,2,3$,
where $\tau^2_2 =  \tau^1_2$, and 
$\tau^2_1$ is bounded by
$(ABCD)$, $(ABDE)$, $(ABCE)$, $(BCDE)$,
 $(ACEF)_1$, $(ACDF)'$, $(ADEF)'$, and $(CDEF)'$.
Here $(ACDF)'$, $(ADEF)'$, and $(CDEF)'$ denote new tetrahedra
obtained as a result of performing a   Pachner move to
$(ACDE) \cup (ACEF)_2$.
Then the last polytope $\tau^2_3$ is bounded by
$(ACDF)'$, $(ADEF)'$, and $(CDEF)'$ that are explained above,
and $(ACDF)$, $(ADEF)$, $(CDEF)$ that used to be faces of $\tau^1_3$.

The $ ( 2 \rightleftharpoons 3)$-move to $(ACDE) \cup (ACEF)$ is illustrated
 in Fig.~\ref{tacolem3}. In the upper left the the 4-cell $\tau^2_1$ 
is shown while $\tau^2_3$ is shown on the upper right.
 In the lower part of the figure the three new tetrahdera
 $(ACDF)'$ $(ADEF)'$, and $(CDEF)'$ are illustrated.

Then we can collapse $\tau^2_3$ to the tetrahedra
 $(ACDF)$, $(ADEF)$, $(CDEF)$ as in the following 3 paragraphs and 
2 tables.

 The polytope 
$\tau^2_3$ is a $4$-cell bounded by  $(ACDF)$, $(ADEF)$, $(CDEF)$,
$(ACDF)'$, $(ADEF)'$, and $(CDEF)'$.
The incidence relations  for these tetrahdra are indicated in the next table. Also see the top two rows of Fig.~\ref{tacolem4}. 

\begin{center}
\begin{tabular}{||c|c|c||} \hline  \hline
Triangles & are faces of & tetrahedra \\ \hline \hline
 $\left\{ \begin{array}{c} (ACD) \\ (ACF) \end{array} \right\}$& $\subset$ & $(ACDF)\cup (ACDF)'$ \\ \hline
 $\left\{ \begin{array}{c} (ADE) \\ (AEF) \end{array} \right\}$& $\subset$ & $(ADEF)\cup (ADEF)'$ \\ \hline
 $\left\{ \begin{array}{c} (CDE) \\ (CEF) \end{array} \right\}$& $\subset$ & $(CDEF)\cup (CDEF)'$ \\ \hline
 $(ADF)$& $\subset$ & $(ACDF)\cup (ADEF)$ \\ \hline
$(DEF)$& $\subset$ & $(ADEF)\cup (CDEF)$ \\ \hline
$(CDF)$& $\subset$ & $(ACDF)\cup (CDEF)$ \\ \hline 
$(ADF)'$& $\subset$ & $(ACDF)'\cup (ADEF)'$ \\ \hline 
$(DEF)'$& $\subset$ & $(ADEF)'\cup (CDEF)'$ \\ \hline 
$(CDF)'$& $\subset$ & $(ACDF)'\cup (CDEF)'$ \\ \hline 
\hline
\end{tabular}
\end{center}

Then perform the taco move to the pair $(CDEF)$ and $(CDEF)'$
that share two faces $(CDE)$ and $(CEF)$. 
This move is illustrated at the bottom of Fig.~\ref{tacolem4}. 
Then the faces $(CDF)$ and $(CDF)'$, $(DEF)$ and $(DEF)'$
are identified after the move respectively.
The result is a $4$-cell bounded by $(ACDF)$, $(ADEF)$,
$(ACDF)'$, and  $(ADEF)'$.
(Precisely speaking these tetrahdera share new faces so that 
we should use the different labels, but adding a new layer of
 labels here will cause more confusion than leveing the old labels intact).
The incidence relations among the triangles and the tetrahedra are summarized
 in the next table.

\begin{center}
\begin{tabular}{||c|c|c||} \hline  \hline
Triangles & are faces of & tetrahedra \\ \hline \hline
 $\left\{ \begin{array}{c} (ACD) \\ (ACF) \\(CDF) \end{array} \right\}$& $\subset$ & $(ACDF)\cup (ACDF)'$ \\ \hline
 $\left\{ \begin{array}{c} (ADE) \\ (AEF) \\ (DEF) \end{array} \right\}$& $\subset$ & $(ADEF)\cup (ADEF)'$ \\ \hline
 $(ADF)$& $\subset$ & $(ACDF)\cup (ADEF)$ \\ \hline
$(ADF)'$& $\subset$ & $(ACDF)'\cup (ADEF)'$ \\ \hline 
\hline
\end{tabular}
\end{center}

The cone move to $(ADEF)$ and $(ADEF)'$ 
(which is illustrated schematically in Fig.~\ref{tacolem5})
followed by the pillow move to $(ACDF)$ and $(ACDF)'$
collapses $\tau^2_3$ to $(ACDF) \cup (ADEF) \cup (CDEF)$
as claimed.

Thus we get two polytopes $\tau^2_1$ and $\tau^2_2$.
Next perform a   Pachner move to
$(ABCE) \cup (ACEF)_1$ which shares $(ACE)_1$.
As a result we get three new tetrahedra
$(ABEF)' \cup (ABCF)' \cup (BCEF)'$.
The $ ( 2 \rightleftharpoons 3)$-move is illustrated in Fig.~\ref{tacolem6}; the labels on the polytopes indicate their antecendents.

Thus we obtain $\tau^3_1$ bounded by
$(ABCD)$, $(ABDE)$,
$(BCDE)$,  $(ACDF)$, $(ADEF)$,  $(CDEF)$,
$(ABEF)'$, $ (ABCF)'$, and $(BCEF)'$,
and $\tau^3_2$ bounded by
$(ABEF)$, $(ABCF)$, $(BCEF)$, and
$(ABEF)' \cup (ABCF)' \cup (BCEF)'$.

Hence we now can collapse $\tau^3_2$ to
the 
tetrahedra $(ABEF)$, $(ABCF)$, and $(BCEF)$
in the same manner as we did to $\tau^2_3$.
The collapsing is indicated in Fig.~\ref{tacolem7}.
 The result is a single polytope $\tau^4$ 
resulted from
$\tau^3_1$ which has the same boundary tetrahedra as
those of the left hand side of the  4-dimensional Pachner move.
Figure~\ref{tacolem7} indicates the resulting polytope at the bottom of the figure. In Fig.~\ref{tacolem8} the 3-dimensional boundary is illustrated.
Notice the following: (1) triangle (ACE) is no longer present; (2)
among the nine tetrahedra illustrated, neither triangle   $(ACE)$ nor triangle $(BDF)$ appears; (3) these are all of the tetrahedral faces of the 
5 simplex that contain neither $(ACE)$ nor $(BDF)$.
Thus we can apply the same method starting with $(BDF)$ 
to get to this polytope. 
This proves that $(3 \rightleftharpoons 3)$-move 
is described as
a sequence of singular moves (cone, taco, and pillow moves) and
  Pachner moves.

$\Box$

\begin{subsect} {\bf Lemma.} \label{pac24}
The $(2 \rightleftharpoons 4)$-move 
is described as a sequence of
cone moves, pillow moves, taco moves,  and $3$-dimensional  Pachner moves.
\end{subsect}
{\it Proof.}
We use the following labeling for the $(2 \rightleftharpoons 4)$-move
in this proof:

$$ (ABCDE) \cup (ABCEF) 
\rightleftharpoons (ABCDF) \cup (ABDEF) \cup (ACDEF) \cup (BCDEF) .$$

Perform 
a 
$(3 \rightleftharpoons 3)$-move (which was proved to be 
a sequence of the singular moves in the preceeding Lemma) to
$(ABCDF) \cup (ABDEF) \cup (ACDEF)$ 
to get $(ABCDE)' \cup (ABCEF)' \cup (ACDEF)'$.
Then the polytope now consists of 
 $(ABCDE)'$, $(ABCEF)'$, $(ACDEF)'$, and $(ACDEF)$. 

Perform a  Pachner move to the tetrahedra
$(ACDF) \cup (ADEF) \cup (CDEF)$, 
that are shared by $(ACDEF)'$ and $(ACDEF)$,
to get $(ACDE)' \cup (ACEF)'$.

This changes $(ACDEF)' \cup (ACDEF)$ to a $4$-cell 
bounded by $(ACDE)$, $(ACDE)'$, $(ACEF)$, and $(ACEF)'$.
The cone move followed by the pillow move collapses this polytope 
yielding $ (ABCDE) \cup (ABCEF)$, the left-hand-side of  the 
$(2 \rightleftharpoons 4)$-move.
$\Box$

\begin{subsect} {\bf Lemma.} \label{pac15} 
The $(1 \rightleftharpoons 5)$-move 
is described as a sequence of
cone moves, pillow moves, taco moves, and $3$-dimensional  Pachner moves.
\end{subsect}
{\it Proof.}
We use the following labelings:

$$ (ABCDE) \rightleftharpoons
(ABCDF) \cup (ABCEF) \cup (ABDEF) \cup (ACDEF) \cup (BCDEF).$$

Perform the  
$(3 \rightleftharpoons 3)$-move to 
$(ABCDF) \cup (ABDEF) \cup (BCDEF)$ to get 
$(ABCDE) \cup (ABCEF)' \cup (ACDEF)'$. 

The $4$-simplices
$(ACDEF)$ and $(ACDEF)'$ share 
all their tetrahedral faces except
$(ADEF)$ (and $(ADEF)'$).
Perform a $(1 \rightleftharpoons 3)$-move 
to 
each of these 
shared tetrahedra 
to get $4$-cells bounded by copies of $(ADEF)$ sharing all the $2$-faces.
Thus the pillow moves will collapse $(ACDEF)$ and $(ACDEF)'$.
The same argument collapses $ (ABCEF) \cup (ABCEF)'$
to get the left-hand-side of the  $(1 \rightleftharpoons 5)$-move.
$\Box$

\begin{subsect} {\bf Remark.\/} 
{\rm 
In \cite{CF} Crane and Frenkel proposed  
constructions of $4$-manifold quantum invariants using 
Hopf categories. 
Hopf categories generalize the definition of Hopf algebra 
to  a 
categorical setting
in the same way that modular categories generalize modules.
One of the conditions in their definition is
called the coherence cube which generalizes the compatibility 
condition of Hopf algebras between mutiplication
and comutiplication (See Section~\ref{Hopfsec}). 
They showed that this condition corresponds to 
the cone move. 
Thus 
Lemmas in this section 
can be used to prove the well-definedness 
of invariants they proposed by showing that 
their definition 
is invariant under 
Pachner moves.
} \end{subsect}

\section{Triangulations and Diagrams} \label{T&Dsec}

In dimension 3, quantum spin networks are used on the one hand
to provide calculations of identities among representations of
 quantum groups \cite{CFS}.
 On the other hand they are cross sections of the dual complex of 
a triangulated 
3-manifold (see Section~\ref{pac23}).  

In this section,  we use similar graphs 
to
relate 
them  
to the dual complex of triangulated 4-manifold.
 We begin the discussion on the local nature of triangulated 4-manifolds
 near 2-dimensional faces.

\begin{sect} {\bf Graphs, 2-complexes,
 and triangulations.\/} \label{graphdef}
\addcontentsline{toc}{subsection}{Graphs, 2-complexes, and triangulations}
 {\rm 
Let $\Phi$ be a triangulation of an oriented closed $4$-manifold $M$.
In this section we associate graphs to triangulations and their duals.
} \end{sect}

\begin{subsect} {\bf Definition.} {\rm
The dual complex $\Phi ^*$ of $\Phi$ is defined as follows.
Pick a vertex $v$ of $\Phi ^*$ in the interior of each
$4$-simplex of $\Phi$. Connect two vertices $v_1 $ and $v_2$
of $\Phi ^*$ 
if and only if the 
corresponding $4$-simplices of $\Phi$ share
a $3$-face. Thus each edge of $\Phi ^*$ is dual to a 
tetrahedron 
of
$\Phi$.
Edges $e_1, \cdots, e_k$ of $\Phi ^*$ 
bound 
a face $f$ 
if and only if 
the
corresponding tetrahedra share a $2$-face of $\Phi$.
A set of $2$-faces $f_1, \cdots, f_k$ of $\Phi ^*$
bounds a $3$-face (a polyhedron) 
if and only if 
the corresponding
faces of $\Phi $ share an edge of $\Phi$.
Finally a set of $3$-faces of $\Phi^*$ bounds a $4$-face 
if and only if
the corresponding edges of $\Phi $ share a vertex.
Thus
$\Phi^*$ gives a CW-complex structure 
to 
the $4$-manifold.
} \end{subsect}

\begin{subsect} {\bf Definition.} {\rm
Let $\Phi$ be a triangulation of a $4$-manifold $M$,
and 
let $\Phi^*$ be the 
dual complex.
Each $3$-face of $\Phi^*$ is a polytope.
Choose a triangulation of each  $3$-face into tetrahedra
so that it defines a triangulation of the $3$-skeleton of $\Phi^*$.
We require that such a triangulation does not have 
interior vertices in the $2$-faces of $\Phi^*$.
Thus the restriction on each $n$-agonal 
$2$-face consists of $(n-2)$
triangles.
Such a choice of triangulation is called a
{\it $3$-face triangulation} (a {\it triangulation} for short)
of $\Phi ^*$.
A  $3$-face triangulation is denoted by $\Phi^!$.
} \end{subsect}

\begin{subsect}{\bf Definition (Carrier Surface).\/} {\rm In each tetrahedron
 of the triangulation $\Phi$,  we embed the dual spine 
to the tetrahedron.
 The intersection of the dual spine with a triangular face 
is a graph consisting of a 3-valent vertex 
with edges intersecting the edges of the tetrahedron. 
There is a vertex 
in the center of the 2-complex at which 
four edges (corresponding to the faces of the tetrahedron) and six faces
(corresponding to the edges) intersect. 
 The union (taken over all tetrahedra in the triangulated 4-manifold)
 of these 2-complexes form a 2-complex, $C$,  
that we call the {\it carrier surface}. 
Let us examine the incidence relations of the carrier surface along faces
 and edges of the triangulation.

Consider a 2-face, $f$, of the triangulation $\Phi$.
 Suppose that $n$ tetrahedra are incident along this triangle $f$. 
 Then the dual face $f^*$
is an $n$-gon. The 4-manifold in a neighborhood of the face $f$
 looks like the 
Cartesian 
product $f \times f^*$. 
The carrier surface in this neighborhood then appears as $Y \times X_n$
 where $X_n$ is the 1-complex that consists of the cone on $n$-vertices
({\it i.e.}, the $n$-valent vertex),           
 and $Y$ is the graph that underlies that alphabet character 
(a neighborhood of a trivalent vertex).
 For example $X_2$ is an interval, $X_3 = Y$, $X_4 =X$, {\it etc.}
We can think of $X_n$ being embedded in $f^*$ with the edges of $X_n$
 intersecting the centers of the edges of $f^*$ and the vertex of $X_n$
 lying at the ``center''  of $f^*$ 
({\it i.e.} we may assume that $f^*$ is a regular polygon).

 Consider an edge, $e$, of $\Phi$, 
and the 3-cell, $e^*$, that is dual to $e$.
 The faces of $e^*$ are $n$-gons, $f^*$,
 that are dual to the triangular faces, $f$, 
which 
are incident to $e$.
 The carrier surface intersects a face $f^*$ in the graph $X_n$.
The carrier surface intersects $e^*$ in a 2-complex that is the cone
 on the union of the $X_n$ where the union is taken over all
 the faces of $e^*$.

The situation is depicted in Fig.~\ref{face}  in which three tetrahedra
 intersect along a triangular face.
 On the right hand side of the figure, we illustrate a {\it graph movie}.
 The two graphs that are drawn there represent the intersection of the carrier
 complex with the boundary of $f \times f^*$. 
In a neighborhood of this face the carrier complex looks like $Y\times Y$.
 In this and subsequent figures, the vertices that are labeled with 
open circles correspond to the dual faces $f^*$. In this figure, 
three such circled vertices appear since the dual face appears on each
 of the duals to the three edges.

}\end{subsect}

\begin{sect}{\bf Faces and diagrams.\/} 
\label{f&d}
\addcontentsline{toc}{subsection}{Faces and diagrams}
{\rm
Suppose that the face $(012)$ of a triangulation of 
a $4$-manifold is shared by three tetrahedra $T_i$, $i=1,2,3$.
Take a neighborhood $N$ of 
 the  
face $(012)$ 
in the $3$-skeleton of the triangulation 
 such that 
$N \cap T_i$ is diffeomorphic to $(012) \times I$
for each $i=1,2,3$.

In Fig.~\ref{face} the
projection of a neighborhood $N$ of the face
$(012)$  
is depicted in $3$-space. 
Denote by $0'$, $0''$, $0'''$ the vertices obtained from the vertex $0$
by pushing it into $T_i$, 
$i=0,1,2,$   
respectively
(they are depicted in Fig.~\ref{face}). 
Similar 
notation is 
used for the other vertices.

The {\it graph movie} 
 for                     
 $N$ is constructed as follows. 
Regard 
$N$ 
as a 
$3$-dimensional
polyhedral complex consisting of the 
following faces:
$(0'1'2')$, $(0''1''2'')$, $(0'''1'''2''')$,
$(011'0')$, $(010''1'')$, $(010'''1''')$,
$(122'1')$, $(122''1'')$, $(122'''1''')$,
$(200'2')$, $(200''2'')$, $(200'''2''')$.
Then  
trivalent vertices are assigned to 
the middle points of the triangular faces 
$(0'1'2')$, $(0''1''2'')$, $(0'''1'''2''')$,
and the middle points of the edges 
$(01)$, $(12)$, $(20)$.
These are connected by segments as indicated in the figure 
where this 1-complex is depicted in two parts.
The middle point in the interior of $N$ is the cone point
of this $1$-dimensional complex. 
Within $N$ we have an embedding of 
the Cartesian  product 
$Y \times Y$ where $Y$ represents 
the obvious graph with one trivalent vertex. 
The graphs on the right of Fig.~\ref{face} 
represent portions of the boundary of 
$Y\times Y$. The space $Y\times Y$ is indicated in Fig.~\ref{YxY} in which 
the subspace $ \circ \times Y$, 
where $\circ$ denotes a vertex, is indicated as a fat vertex times $Y$. 
The labels on the Figure will be explained 
in Section ~\ref{deflabel} and Fig.~\ref{phi}.

Below (Sections~\ref{dc&t} and \ref{Hopfsec})
 we will
 relate these spin networks to cocycle conditions 
in a specific Hopf category. In this way, we will obtain a direct connection
 among these structures. 

\begin{figure}
\begin{center}
\mbox{
\epsfxsize=4in
\epsfbox{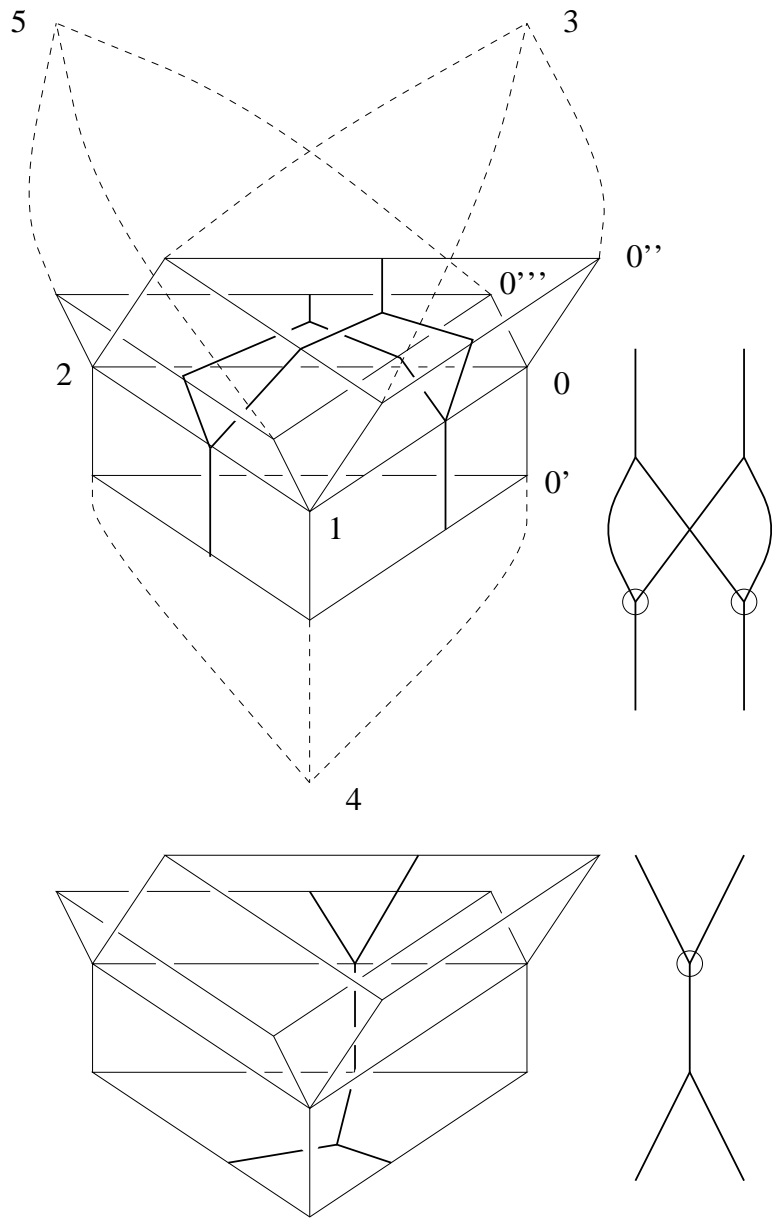}
}
\end{center}
\caption{Graphs and  triangulations around a face}
\label{face}
\end{figure}

\begin{figure}
\begin{center}
\mbox{
\epsfxsize=6in 
\epsfbox{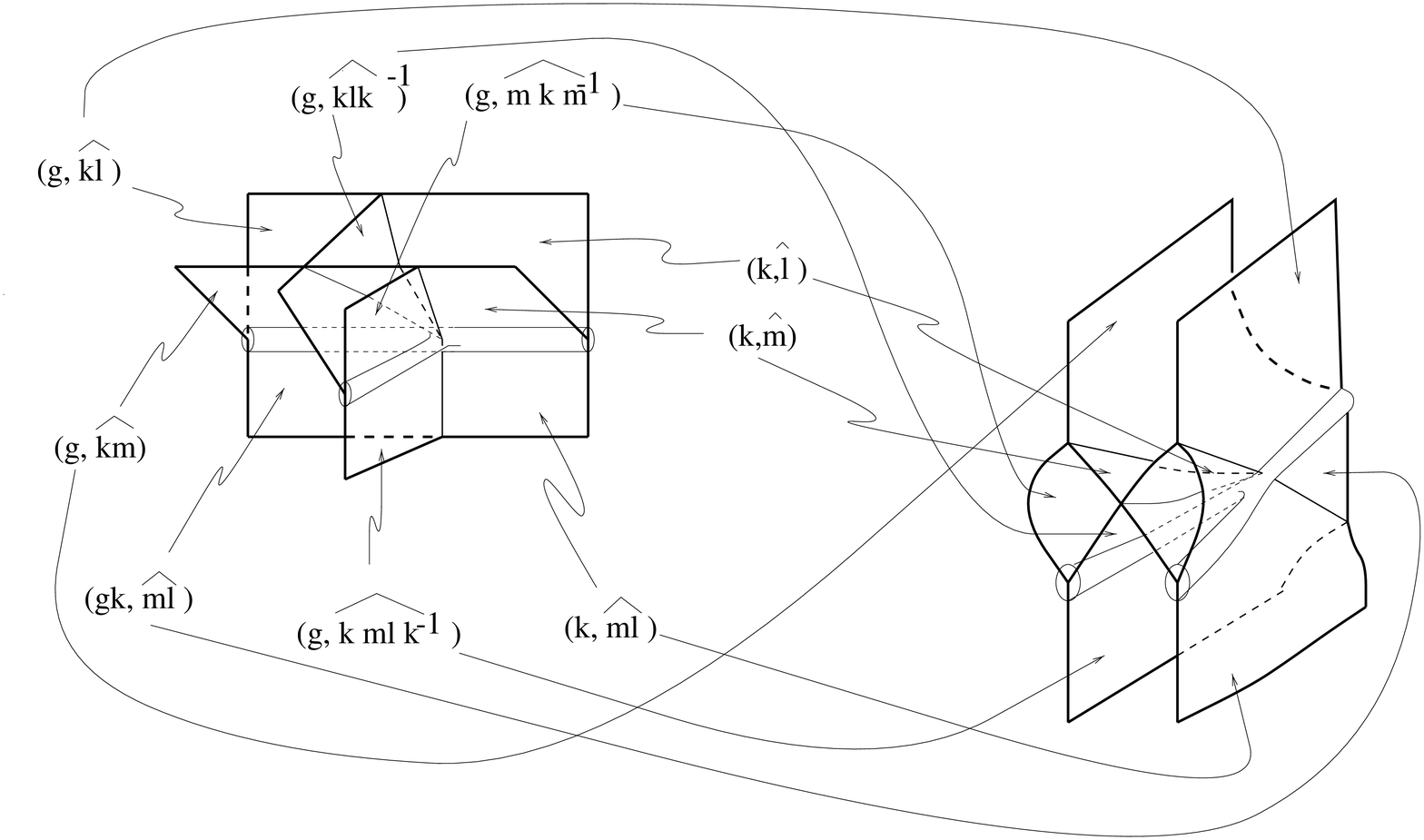}
}
\end{center}
\caption{The product space YxY}
\label{YxY}
\end{figure}

}\end{sect}

\begin{subsect}{\bf Defintion.\/}
{\rm We perturb the carrier surface  to construct a 2-dimensional complex that has the following properties:
\begin{enumerate}

\item
The vertices of the complex all have valence 4 or valence 6.

\item 
Exactly three sheets meet along an edge;

\item
The set of edges can be partitioned into two subsets;
 we color the edges accordingly. 

\item 
A valence 4 vertex  has 4 edges of the same 
color    
incident to it;

\item
A valence 6 vertex has 3 edges of each color incident to it. 

\item
Thus, the 2-complex has a tripartite graph as its 1-complex
 and a bipartition on the set of edges.  
 
\end{enumerate}

Such a 2-complex will be called a {\it perturbed carrier}.
} \end{subsect}

\begin{figure}
\begin{center}
\mbox{
\epsfxsize=4in
\epsfbox{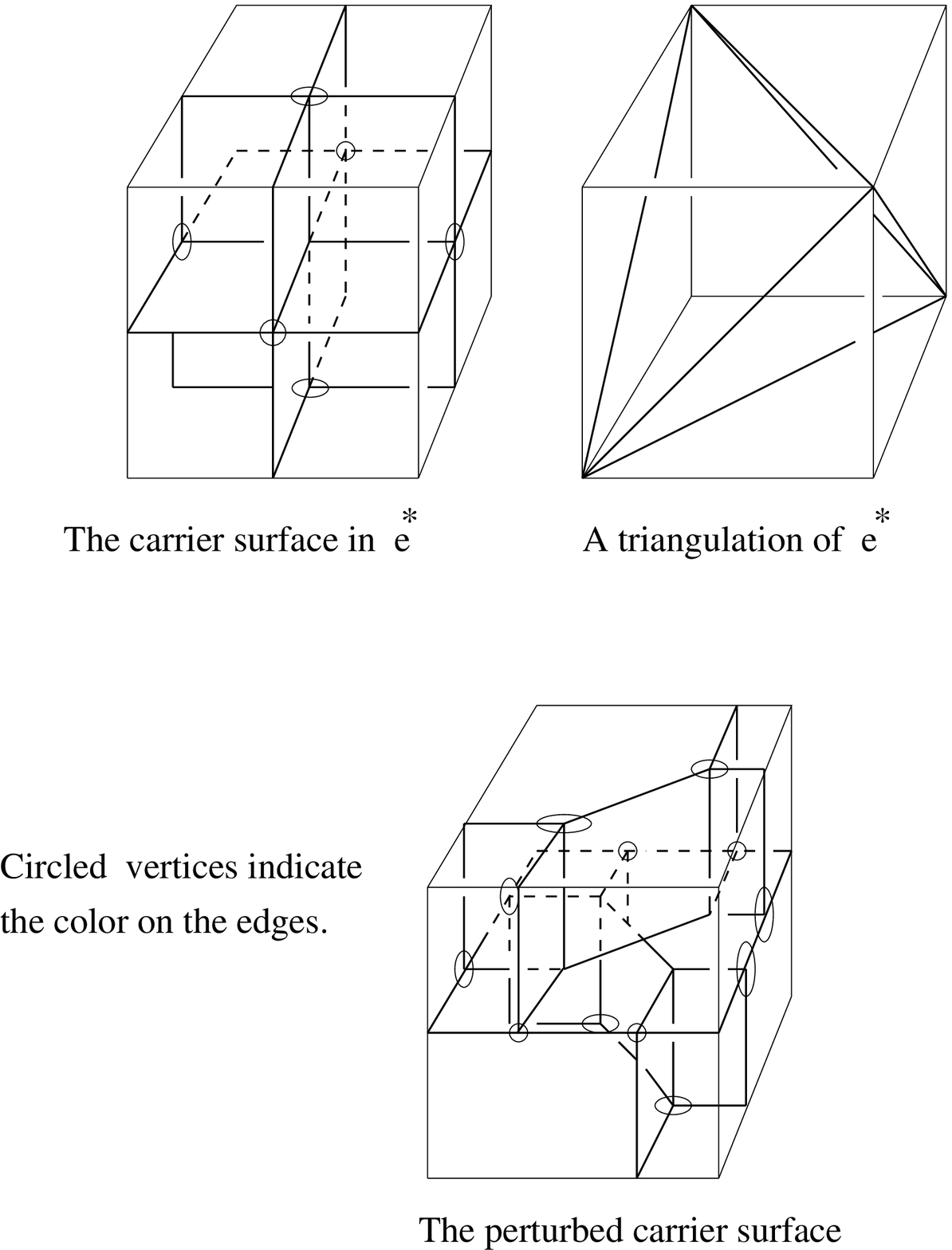}
}
\end{center}
\caption{A neighborhood of an edge whose dual is a cube}
\label{pertcube}
\end{figure}

\begin{subsect}{\bf Lemma.} \label{corresplemma}
 A perturbed carrier can be 
constructed from          
 the carrier surface by means of a 3-face triangulation. 
\end{subsect}
 {\it Proof.} Consider a 3-face triangulation; recall this is a
triangulation of the dual 3-cells of the triangulation $\Phi$,
 and a 3-face is the dual to an edge $e^*$. A $n$-agonal face 
of $e^*$ is divided into $(n-2)$ triangles. The graph $X_n$ in 
the $n$-agonal face is replaced with the dual to the triangulation.
 In $e^*$, the cone on the union of the $X_n\/$s is replaced by
 the union of the duals to the tetrahedra in the triangulation. 
These are the surfaces with 6 faces, 1 vertex, and 4 edges;
they glue together in $e^*$ to form the subcomplex
in which  
all of the vertices 
have one color.
An example is illustrated in Fig.~\ref{pertcube} in which 
the dual of an edge is a cube.

The vertices that have two different colored edges 
incident to them are found on the triangular faces of the 3-face
 triangulation. Three of the edges are coming from the dual face,
 the other three edges are coming from the dual complex of the original
 tetrahedra. 
The local structure at the 6-valent vertices was explained in detail above.
 This completes the proof. $\Box$

\begin{subsect} {\bf Definition.\/} {\rm
A {\it graph movie} is a sequence of graphs that 
appear as cross 
sections 
of a portion of the perturbed carrier
when a height function is chosen, such that the stills of movies are graphs 
having trivalent (circled and uncircled) 
vertices  
and between two stills,
the movie changes in one of the following ways:
} \end{subsect}

\begin{figure}
\begin{center}
\mbox{
\epsfxsize=4in               
\epsfbox{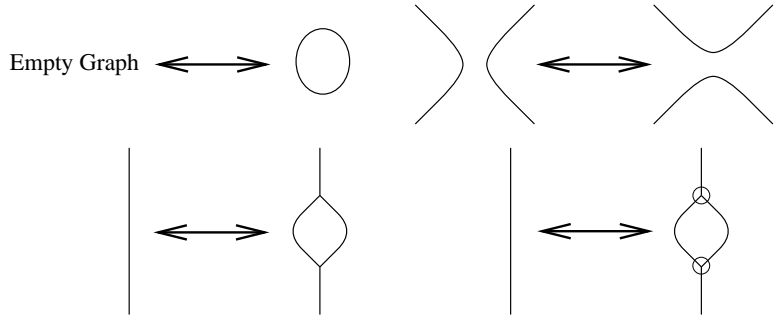}
}
\end{center}
\caption{Some elementary changes of graph movies}
\label{othermovies}
\end{figure}

\begin{enumerate}

\item The change of the movie at a face (a $6$-valent vertex 
of a carrier surface)  is as defined above
(the change of graphs  shown in Fig.~\ref{face}).

\item The change of the movie at a  $4$-valent 
vertex is as depicted in Fig.~\ref{alpha} bottom.

\item The changes of the movie at critical points of edges and faces 
of the carrier surface is generic. 
They are depicted in Fig.~\ref{othermovies}.
 \end{enumerate}

In the graphs we use circled 
 vertices and uncircled 
vertices.
These are cross section of two types of edges. 
In the figures of carrier surfaces
(Fig.~\ref{surfcubeLHS},~\ref{surfcubeRHS}, and \ref{pertcube}),
the edges corresponding to
circled vertices are depicted by thin tubes.
The graph movie defined here includes definitions given above 
(which are clearly equivalent). The graph movie allows us to view the 
perturbed carrier via a sequence of 2-dimensional cross-sections whereas 
the carrier surface itself does not embed in 3-dimensional space.

\begin{sect}{\bf 
Taco moves and graph movies.\/} 
\addcontentsline{toc}{subsection}{Taco moves and graph movies}
{\rm
Herein we directly relate the  graph movies
to the taco move.
 In Figs.~\ref{tacomovieLHS} and \ref{tacomovieRHS} the left-hand-side 
and the right-hand-side of the taco move are depicted, respectively.
 In each figure, the underlying union of tetrahedra remains unchanged
 from frame to frame. Instead the thick lines
change as follows. Consider the $(i,j)\/$th entry of the figure to
 be that illustration in the $i\/$th row $j\/$th column.
 Going from the $(i,1)\/$st entry to the $(i,2)\/$nd entry,
 the graphs change by one of the graph movie changes (either going across 
$Y\times Y$ or going across tetrahedra).
 There is no change from the $(i,2)\/$nd entry to the $(i+1,1)\/$st entry.
In these figures thick lines indicate the graph that was defined 
in Section~\ref{f&d}. 
 The transitions between the two entries on the same row may be visualized
 by means of a cross-eyed stereo-opsis. Place a pen in the center of
 the figure, and move the pen towards your face while keeping it in focus.
 The two images  on the left and right 
should converge into one with the thick lines popping out of
 the plane of the paper. In this way, the   difference between
 the figures can be experienced directly.

Observe that the differences in the graphs are illustrated as well
 in the Fig.~\ref{cocymovie1} which illustrates the graph movies for
 the cocycle conditions and which is obtained by purely algebraic information. 
The time elapsed version of the graph movie for the taco move
is illustrated in Fig.~\ref{surfcubeLHS} and Fig.~\ref{surfcubeRHS}.
Similar diagrams can be drawn for the cone move and the pillow move
 and in this way a direct correspondence can be obtained among the moves,
 the cocycle conditions, and the axioms of a Hopf category
(Section~\ref{Hopfsec}).
The taco, cone, and pillow moves all
correspond to the first coherence cube. 
The correspondence among these moves should not be surprising
 since all of these moves correspond to splitting a tetrahedron open
 (the higher dimensional analogue of the coherence relation between
 multiplication and
 comultiplication).

\begin{figure}
\begin{center}
\mbox{
\epsfxsize=4in
\epsfbox{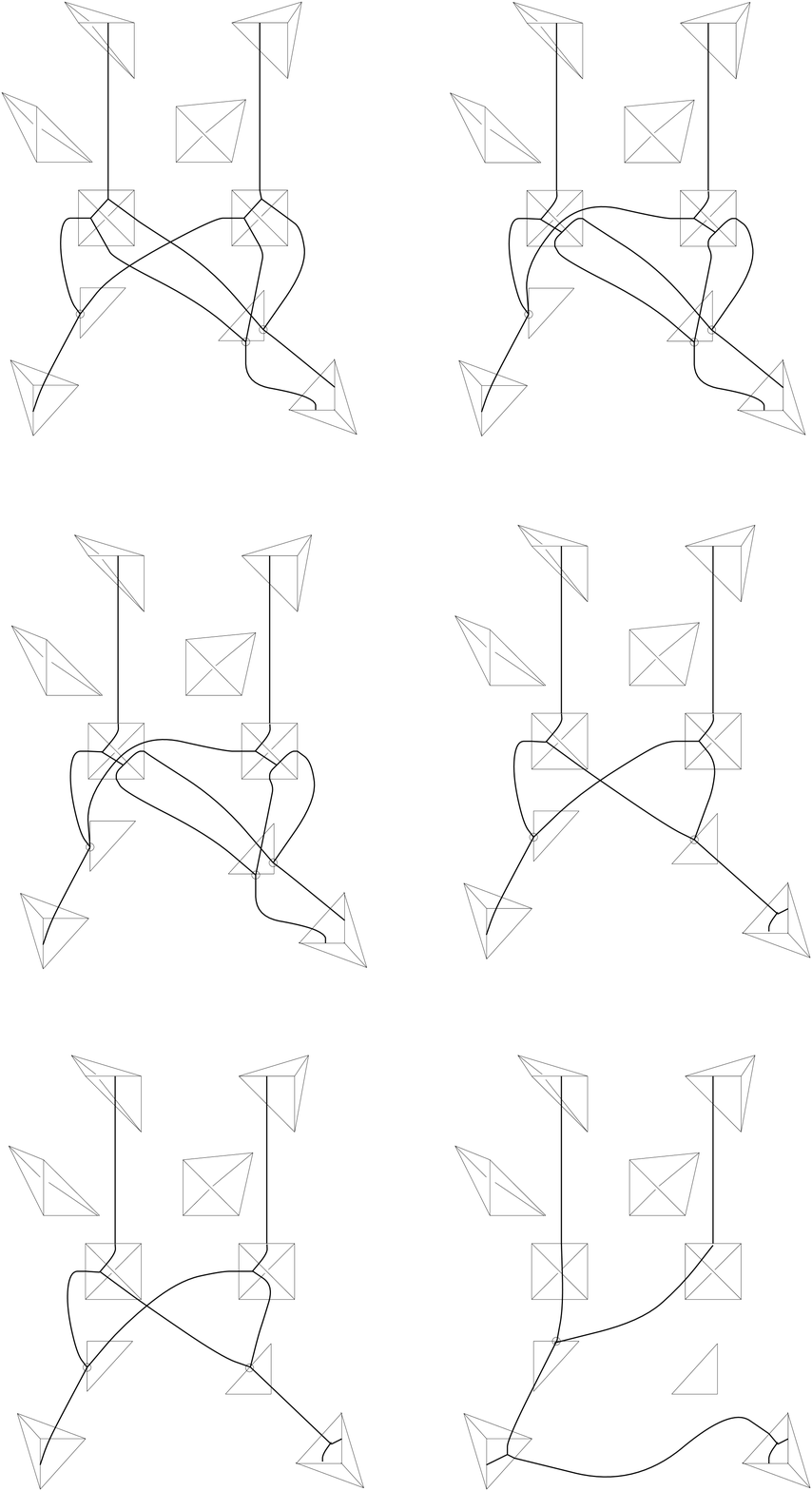}
}
\end{center}
\caption{The
taco move and graph movie, left-hand-side}
\label{tacomovieLHS} 
\end{figure}

\begin{figure}
\begin{center}
\mbox{
\epsfxsize=4in
\epsfbox{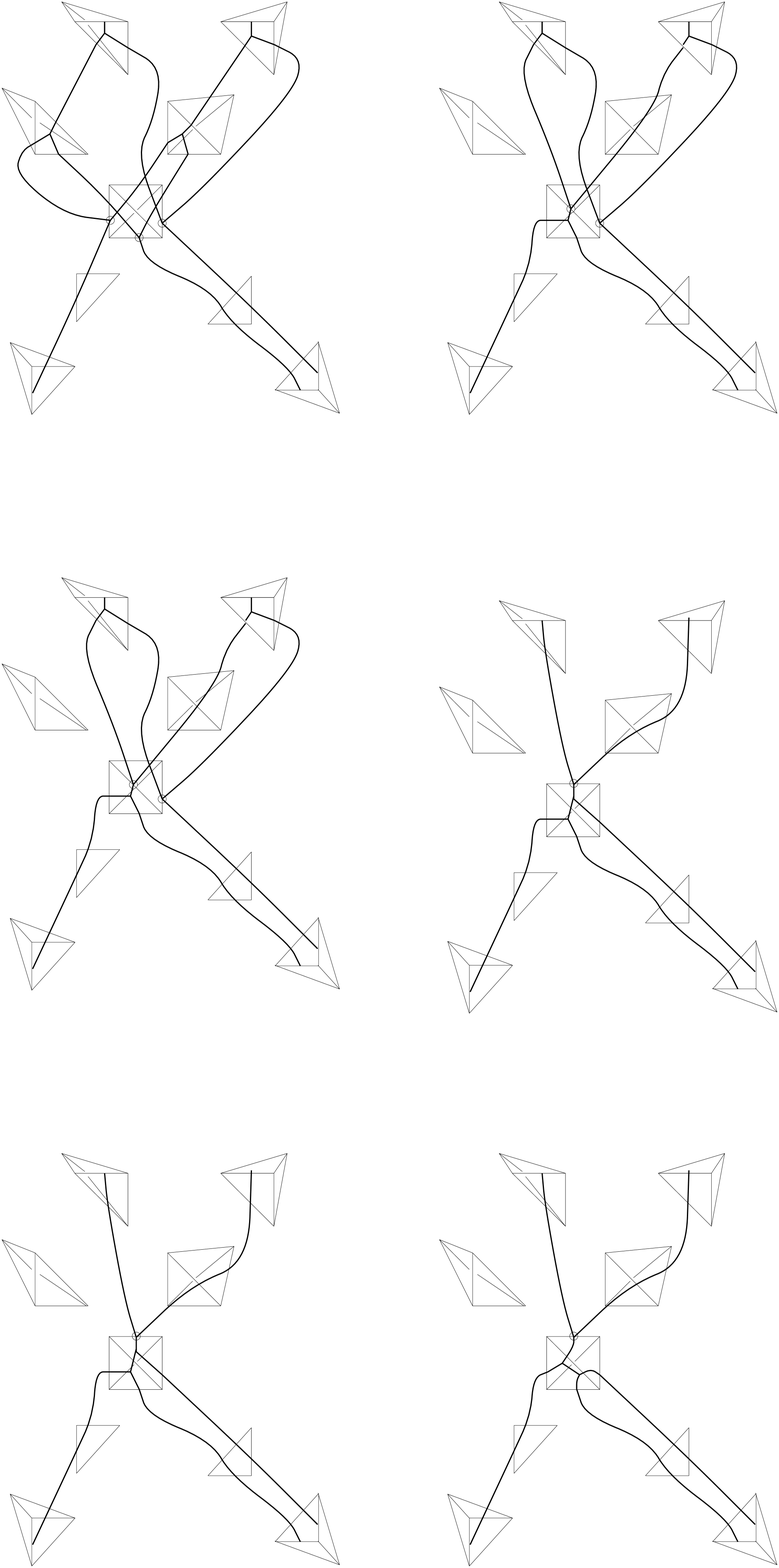}
}
\end{center}
\caption{The 
taco move and graph movie, right-hand-side}
\label{tacomovieRHS} 
\end{figure}

}\end{sect}


\section{Cocycles and cocycle conditions} \label{data}

In this section, we list cocycles and their equalities
that will be used in the following sections.
These cocycles are
given in \cite{CYex} in relation to Hopf categories
(See Section~\ref{Hopfsec}).
Some non-trivial examples are given therein. 
First, we mention that two of the cocycle 
conditions are depicted in 
Figs.~\ref{cocymovie1} and \ref{cocymovie2}
as relations to graph movies 
where the edges of the graph
have been colored with pairs 
of group elements and dual group 
elements. These graph movies
correspond to the dual graphs 
that correspond to the taco move (Fig.~\ref{tacomovieLHS} and \ref{tacomovieRHS}). The coloring will be explained in the subsequent section.

Let $G$ be a finite group and $K^{\times}$ be the multiplicative group
of a field $K$. 
Let $C_{n,m}=C_{n,m} (G, K^{\times})$ 
denote the abelian group of all functions from
$G^n \times \hat{G}^m$ to  $K^{\times}$
where
$$ \hat{G} = \{ \hat{g} : G \rightarrow  K |
g \in G, \hat{g}(h)=1 \;  \mbox{if} \; g=h, 
          \hat{g}(h)=0 \;   \mbox{if} \; g\neq h \} . $$

We need the following functions (called cocycles if they 
satisfy the conditions given in the next section).

\begin{itemize}

\item 
$\alpha ( g, k, m ;  \hat{n} ) \in C_{3,1} $,

\item 
$\beta (g ;  \hat{i}, \hat{j}, \hat{k} ) \in C_{1,3} $,

\item
$\phi (g, k; \hat{m}, \hat{n} ) \in C_{2,2} $. 

\end{itemize}

\begin{sect}{\bf Cocycle conditions.\/} 
\addcontentsline{toc}{subsection}{Cocycle conditions}
{\rm 
The  following are called the cocycle conditions
\cite{CYex}.
} \end{sect}

\begin{itemize}
\item 
$\alpha (k,m,p;\hat{q}) \alpha (g, km, p;\hat{q} ) \alpha (g, k, m; \widehat{p q p^{-1}} )
$ 

$= \alpha (gk, m, p; \hat{q} ) \alpha (g, k, mp; \hat{q} )$

\item 
$\beta(g; \hat{j}, \hat{k}, \hat{\ell}) 
\beta (g; \hat{i}, \widehat{jk}, \hat{ \ell})
\beta ( g; \hat{i}, \hat{j}, \hat{k} )
$

$= \beta (g; \hat{ij}, \hat{k}, \hat{\ell} ) 
\beta (g; \hat{i}, \hat{j}, \widehat{k\ell} ) 
$,

\item
$ \alpha (g,k,m;\hat{p}) \alpha ( g,k,m;\hat{q} )
\phi (k,m;\hat{p}, \hat{q} ) \phi (g,km;\hat{p}, \hat{q} )
$

$ = \phi( g,k; \widehat{mpm^{-1} }, \widehat{mqm^{-1}} )
\phi (gk, m ; \hat{p} \hat{q} )
\alpha (g,k,m;\widehat{pq} )
$,

\item 
$ \phi (g,k;\hat{p}, \hat{r} ) \phi ( g,k; \widehat{pr}, \hat{s} ) 
\beta (gk; \hat{p}, \hat{r}, \hat{s} ) 
$   

$= \beta (g; \widehat{k p k^{-1} } ,   
            \widehat{krk^{-1} } , \widehat{ksk^{-1} } )
\beta (k; \hat{p}, \hat{r} , \hat{s} )
\phi (g,k; \hat{r},  \hat{s} )
\phi (g,k; \hat{p}, \widehat{rs} )
$. 

\end{itemize}

\begin{sect}{\bf Cocycle symmetries.\/} 
\addcontentsline{toc}{subsection}{Cocycle symmetries}
{\rm 
In addition to the above cocycle conditions, we will suppose 
that the cocycles satisfy some equations that correspond to the symmetries
 of tetrahedra and of the space $Y\times Y$.
 The imposition of such conditions will be sufficient to construct
 an invariant. We do not know if the symmetry conditions are necessary.
(They may be satisfied automatically for certain cocycles,
or the invariants may be defined without symmetry conditions.) 
} \end{sect}

\begin{subsect} {\bf Definition.\/}
\label{cocysym} {\rm
The following are called the
{\it cocycle symmetries}.
} \end{subsect}

\begin{itemize}
\item
$\alpha (g,k,m; \hat{n}) = \alpha (g^{-1}, gk, m; \hat{n} )^{-1} $

$= \alpha (gk, k^{-1}, km; \hat{n})^{-1} = 
\alpha ( g, km, m^{-1} ; \hat{\ell} )^{-1},$

where $\ell = mnm^{-1}$.

\item
$\phi (g, k; \hat{m}, \hat{\ell})  
= \phi (g^{-1}, gk; m , \ell  )^{-1} 
= \phi (gk, k^{-1}; kmk^{-1}, k\ell k^{-1})^{-1}$    

$= \phi (g, k ; m \ell , \ell^{-1} )^{-1} 
= \phi (g,k;  m ^{-1}, m \ell ^{-1} )^{-1} .$

\item
$\beta (g ; \hat{h}, \hat{\ell}, \hat{n} )= 
\beta (g; \widehat{h^{-1}} , \widehat{h \ell}, \hat{n} )^{-1} $

$= \beta (g; \widehat{h \ell}, \widehat{ \ell ^{-1}} , \widehat{\ell n} )^{-1}
= \beta (g; \hat{h}, \hat{\ell}, \widehat{n^{-1}} ) ^{-1} .$

\end{itemize}

\section{Labels, weights, and the partition function} \label{bzman}

\begin{sect}{\bf Labeling.\/}
\addcontentsline{toc}{subsection}{Labeling}
 {\rm
Let $\Phi$ denote a triangulation of the $4$-manifold $M$, and let 
 $\Phi^*$     
 denote the dual complex. 
Each $3$-face of $\Phi^*$ is a polytope 
that corrresponds to an edge of $\Phi$.
Choose a triangulation of 
the $3$-skeleton of $\Phi^*$.
There are no 
interior vertices in the $2$-faces of $\Phi^*$.
Thus the restriction on each 
polygonal 
$2$-face consists of $(n-2)$
triangles.
As before, 
 such a choice of triangulation is called a
{\it $3$-face triangulation} 
of $\Phi ^*$.
A  $3$-face triangulation is denoted by $\Phi^!$.

When an order, ${\cal O}$, is fixed for 
the vertex set,
${\cal V}$, 
we define the orientation of dual edges as follows.
A vertex of  $\Phi ^*$ is a $4$-simplex of $\Phi$
whose vertices are ordered. Then $4$-simplices are ordered 
by lexicographic ordering of their vertices.
This gives an order on vertices of  $\Phi ^*$,
giving orientations of edges of  $\Phi ^*$.
Orientations of edges of $\Phi^!$ are 
ones that are compatible with the above orientation. 
} \end{sect}

\begin{figure}
\begin{center}
\mbox{
\epsfxsize=4in
\epsfbox{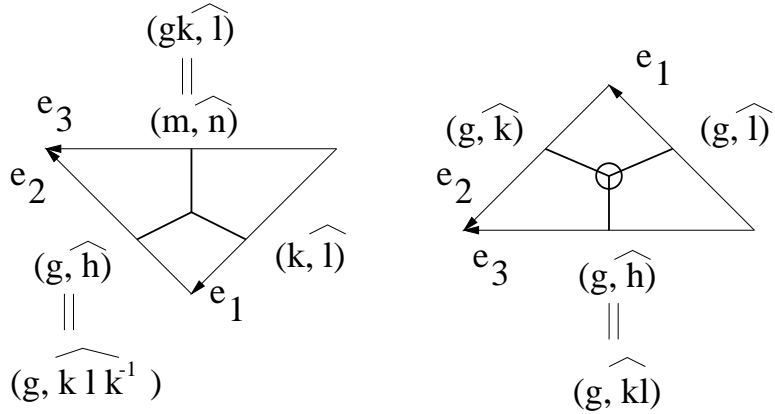}
}
\end{center}
\caption{Rules of cocycle colors}
\label{rule}
\end{figure}

\begin{subsect} {\bf Definition.} 
\label{deflabel}
{\rm
A {\it labeling} (or {\it color}) of $\Phi$ 
with oriented edges 
with respect to  a finite group $G$ is a function
$$ S_0 : {\cal ET} \rightarrow {\cal G} $$
where ${\cal G}= \{ (g, \hat{h}) \in G \times \hat{G}  \} $ and
$$ {\cal ET}= \{ (e, t) \in {\cal E} \times {\cal T}
 | e \subset t \}. $$
Here ${\cal E} $ denotes the set of oriented edges, 
and ${\cal T}$ is the set of tetrahedra. 
We require the following compatibility condition.

If $(e_1, e_2, - e_3)$ forms 
an oriented boundary of a face
of a tetrahedron $t$,   
$S_0(e_1, t) = (k, \hat{\ell}) $, 
and $S_0(e_2, t)=(g,\hat{h})$,
then $S_0(e_3, t)=(m, \hat{n})$,
where  
$m=gk$,  $n=\ell$, and $k^{-1} h k = \ell$.
We call this rule {\it the local rule of colors at a triangle}
(or simply a {\it local rule}). 
The situation is depicted 
on the left of  
Fig.~\ref{rule}.

When an order of vertices is given, the edges are oriented by
ascending order of vertices (if the vertices $v$ and $w$ of an edge 
have 
the order $v<w$, then the edge is oriented from $v$ to $w$).
However in this definition the order on vertices is not required,
although orientations on edges are required.
For an oriented edge $e$, the same edge with the opposite 
orientation is denoted by $-e$.
Consider the edge $e_2$ 
on the left of Fig.~\ref{rule}, 
and reverse the orientation
of $e_2$ to get $- e_2$.
Then the color of $- e_2$ is required to be 
  $S(-e_2,t) = (g^{-1}, \widehat{ gk \ell k^{-1} g^{-1} } )$
where $S(e_2, t)=(g, \widehat{ k \ell k^{-1}  } )$
as depicted in the figure.
In other words, the color for an edge with reversed 
orientation
is defined to satisfy the local requirement of 
the left of Fig.~\ref{rule}.

} \end{subsect}

We often use 
sets of  
non-negative integers to represent
simplices of $\Phi$.
For example, fix a $3$-face  (or tetrahedron)
 $T$ of $\Phi$. 
 Let $0$, $1$, $2$, and $3$ denote the vertices
of $T$.
For a pair of an  oriented  edge $(01)$ and a  tetrahedron
$T= (0123)$ a labeling assigns a pair $(g, \hat{h})$ 
 which we sometimes denote by $S_0(01 | 0123)
= S_0( (01), (0123))$.
When a total order is fixed, the integers are assumed 
to have the compatible order ($0 < 1 < 2 < 3 $).

\begin{figure}
\begin{center}
\mbox{
\epsfxsize=4in
\epsfbox{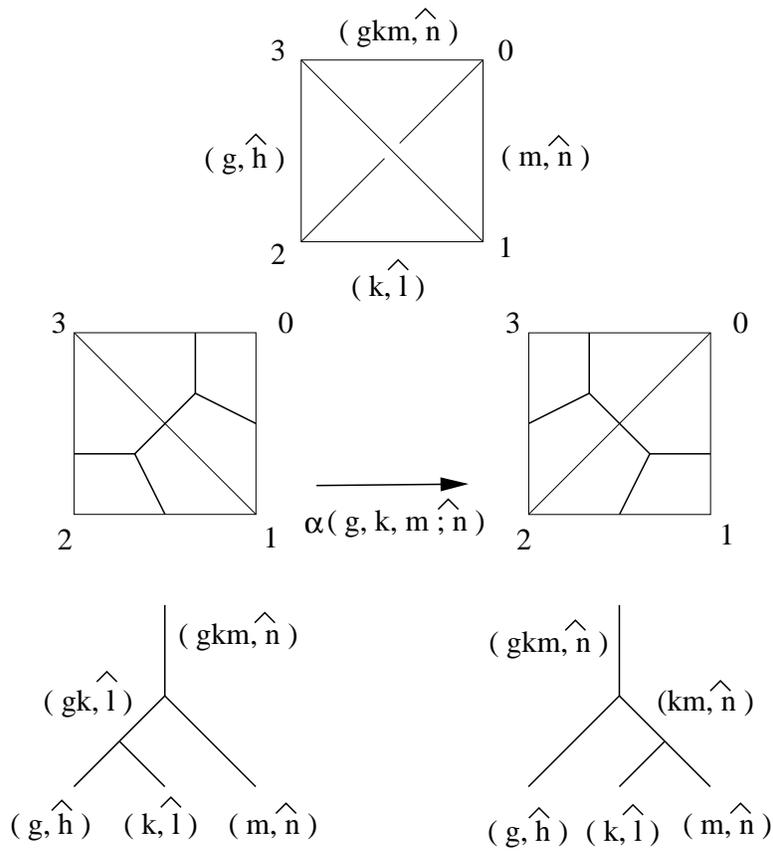}
}
\end{center}
\caption{Weight for tetrahedra}
\label{alpha}
\end{figure}
 
\begin{figure}
\begin{center}
\mbox{
\epsfxsize=4in
\epsfbox{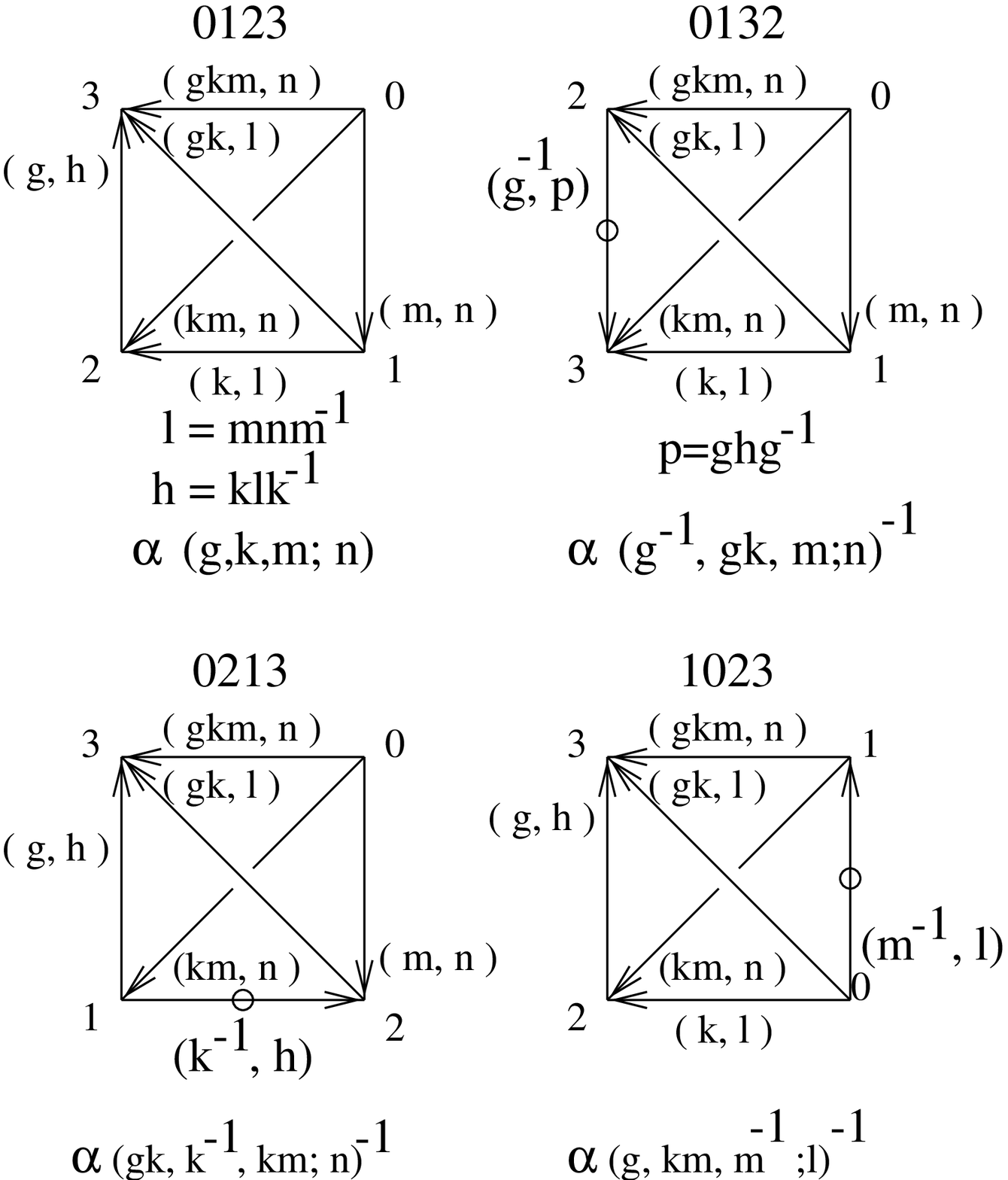}
}
\end{center}
\caption{Symmetry of colors of a tetrahedron}
\label{alphasymmetry}
\end{figure}

We will show (Lemma~\ref{6in1}) that there is a coloring of each
 tetrahedron 
satisfying  
the local rule.
 Furthermore, we will show that changing the orientations of edges of 
a colored tetrahedron results in a unique coloring.

\begin{subsect} {\bf Definition.} {\rm
A {\it labeling} (or {\it color})  of $\Phi^!$
with oriented dual edges       
 is a function
$$ S^! : {\cal EP^!} \rightarrow {\cal G} $$
where
$${\cal EP^!}=\{ (p, e) \in {\cal P^*}\times{\cal E^!}
 | e \subset p \} .$$
Here ${\cal E^!}$ (resp. ${\cal P^*}$ ) denotes
the set of oriented edges (resp. 3-polytopes) of $\Phi^!$
(resp. $\Phi^*$). 
The following compatibility conditions are required.

If $(e_1, e_2, - e_3)$ form an oriented boundary of a face
of a tetrahedron $t$ of $\Phi^!$,  
then the first factors of colors (group elements) 
coincide, and if they are  
$S^! (e_1, t) = (g,\hat{l})$, 
and $S^! (e_2, t)= (g, \hat{k})$,
then $S^! (e_3, t)=(g, \hat{h})$,
where it is required that $h=kl$.

When an order of vertices is given, the edges are oriented by
ascending order of vertices as before.
Consider the edge $e_2$ in the Fig.~\ref{rule} right, and reverse the orientation
of $e_2$ to get $- e_2$.
Then the color of $- e_2$ is required to be 
  $S^!(-e_2,t) = (g, \widehat{  k^{-1}  } )$
where $S(e_2, t)=(g, \widehat{ k   } )$
as depicted in the figure.
In other words, the color for an edge with reversed 
orientation
is defined to satisfy the local requirement of the 
right side   of
Fig.~\ref{rule}.

}\end{subsect}

In the figure, dual graphs in triangles  are also depicted.
We put a small circle around a trivalent vertex for 
the dual faces.
As in the case for tetrahedra, dual tetrahedra can be colored, 
and changing the orientation for colored dual tetrahedra gives a
 unique new  coloring (Lemma~\ref{6in1}).
Note that there is a pair $(p, e')\in {\cal EP^!}$
which is dual to a pair $(e, t)\in {\cal ET} $, in the sense
that $e'$ is dual to 
the tetrahedron
$t$ and $p$ is dual to 
the edge 
$e$.
However there are pairs in ${\cal EP^!}$ that are not to
dual to pairs in ${\cal ET}$.

\begin{subsect} {\bf Definition.} {\rm
A {\it labeling} (or {\it color}) of $ \Psi = \Phi \cup \Phi^!$ is a function
$$  S :   {\cal ET}\cup {\cal EP^!} \rightarrow {\cal G}          $$
such that $ S(p, e')= S(e,t)$
if $(p, e')\in {\cal EP^!}$ is dual to $(e, t)\in {\cal ET} $.
This function is also called a ${\it state}$.
For a particular pair $(e,t) \in {\cal ET}$ (resp. $(p, e') \in {\cal
EP^!}$),
the image $S(e,t)$ (resp. $S(p, e')$) is also called a {\it spin}.
This is sometimes denoted by $S(p|e')$.
} \end{subsect}

\begin{subsect} {\bf Definition.\/} {\rm
We say that two simplices are {\it adjacent} if they intersect.
We say a simplex $\sigma$ and a dual simplex $\tau$ 
are 
{\it adjacent} if 
$\sigma$ intersects 
the polyhedron of the dual complex in which $\tau$ is included.
} \end{subsect}

\begin{figure}
\begin{center}
\mbox{
\epsfxsize=4in
\epsfbox{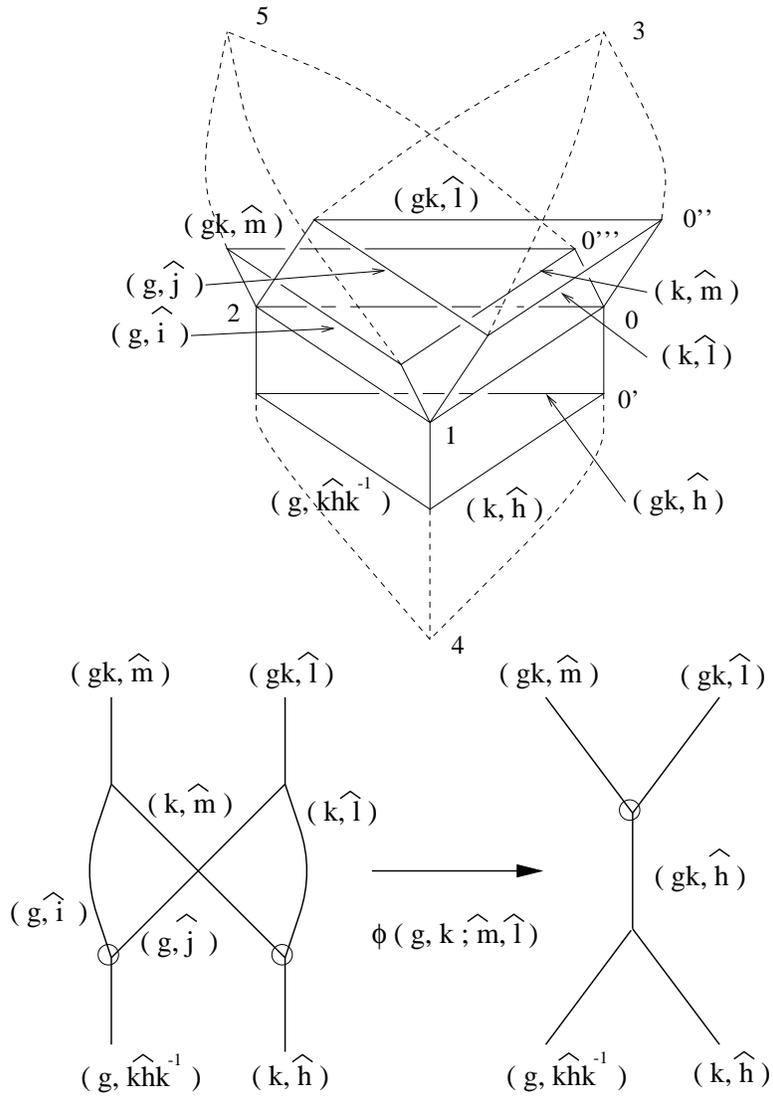}
}
\end{center}
\caption{Weight for faces}
\label{phi}
\end{figure}

\begin{subsect} {\bf Lemma.\/}
\label{6in1}
(1) For a tetrahedron or dual tetrahedron, 
there are colors satisfying the local rule at 
every face or dual face.

(2) There are colors on the edges and dual edges adjacent to
a given face satisfying the local rules.

(3) Let $C$ be a color assigned to the oriented edges (or dual edges) 
of a tetrahedron
(or dual tetrahedron).
Let $C'$ be a color assigned to the same tetrahedron (or dual tetrahedron) with
orientations reversed on some of the edges.
Then $C'$ is 
uniquely determined. If the color $C$ is assigned to oriented edges and dual
edges that are adjacent to a face, then the color $C'$ is uniquely determined
 when some of the edges or dual edges have their orientations reversed. 

\end{subsect}
{\it Proof.}
We prove (1) and (2) in the case of tetrahedra.
 The proof for dual tetrahedra is similar and follows from \cite{Wakui}. 

For tetrahedra, the situation is depicted in Fig.~\ref{alpha}.
In the top of the figure, a 
tetrahedron 
with colors on three edges
is depicted. Other edges receive compatible colors that 
are determined by these three.
In the middle,  pairs of front faces and back faces 
are depicted in the left and right respectively,
together with dual graphs. In the bottom of the figure,
only the dual graphs are depicted, together with colors on
all the edges.
We check that the first factors of colors match
in multiplication convention in Fig.~\ref{rule}.
The second factors are also checked as follows:
in the bottom left figure, $\ell=mnm^{-1}$ and 
$h=k \ell k^{-1}$. 
In the bottom right figure, 
$\ell =mnm^{-1}$ (the same relation as above) 
and $h= ((km) n (km)^{-1} )$
which reduces to 
$kmn m^{-1} k^{-1} = k \ell  k^{-1}$, 
the same relation as above.
Thus the requirements on faces match at a tetrahedron.

To prove (3) for tetrahedra, first consider the case 
where the orientation of the edge $23$ is reversed. 
Then the face $(123)$ forces the change
 $S(- 23 | 0123)= (g^{-1}, \widehat{  g k \ell k^{-1} g^{-1}  } )$.
The other face $(023)$ forces the change
$S(- 23 | 0123)= (g^{-1}, \widehat{ (gkm)n(gkm)^{-1} } ) $.
These are equal since $\ell = mnm^{-1}$.
The situation is depicted in the top right of Fig.~\ref{alphasymmetry} 
where the reversed orientations are depicted by a small circle
on the edge. The top left figure indicates the original colors. 
In the figure, the ``hats'' on the dual group elements 
are abbreviated for simplicity. 

The other cases when the orientation of
a single edge is reversed, are also depicted for the cases
$(0213)$, $(1023)$.
The general case follows
because all cases are obtained by compositions of these 
changes.

Next consider statement (2).
In Fig.~\ref{phi} the colors are depicted using dual graphs
(identify this graph with the graph in Fig.~\ref{face}).
First we check the orientation conventions in the figure.
Identify the circled vertex of the 
right-hand-side graph of bottom of Fig.~\ref{face}
with the right-hand-side of Fig.~\ref{rule}.
Then the orientation conventions of dual edges 
coincide 
where $e_1$ of 
Fig.~\ref{rule} corresponds to the edge $(02) \subset (0123)$,
$e_2$ to $(02) \subset (0125)$, $e_3$ to $(02) \subset (0124)$.
The tetrahedron $(0123)$ is shared by $(01234)$ and $(01235)$,
that are ordered as $(01234) < (01235) < (01245)$ among three
$4$-simplices. Thus this correspondence to $e_1$ matches with
the definition of the orientation of dual edges.

Now we check the constraints. 
In the left bottom of Fig.~\ref{face}
 the following relations must hold:
$i=kmk^{-1}$ (from the top left vertex),
$j=k \ell k^{-1}$ (from the top right vertex),
$h=m \ell $ (from the bottom right vertex),
and $ij=khk^{-1}$ 
(from the bottom left vertex).
The last relation is reduced by substitution to
$k (m \ell ) k^{-1} $ both sides, so that the weight is compatible.
In the right of the figure, we get the relation
$h=m \ell $ from the top vertex, which is the same as above,
and the condition for the bottom vertex 
is already incorporated (by using $khk^{-1}$ in bottom left).
Thus the colors around a face 
are 
compatible.

Now let us check that the orientation conventions are compatible 
in Fig.~\ref{phi}.
The orientations on the edges are the  orientations
from the vertex ordering as seen in the figure.
The orientations on dual edges are checked as follows.
In the figure the face $(012)$ is shared by three $4$-simplices,
$(01234)$, $(01245)$, and $(01235)$. 
The dual edge labeled by $(gk, \hat{\ell})$ is dual to the 
tetrahedron $(0123)$ and oriented from 
$(01234)$ to $(01235)$, corresponding to the edge $e_1$ 
on the right of Fig.~\ref{rule}.  
Respectively, the one labeled by $(gk, \hat{m})$ 
goes from  $(01235)$ to $(01245)$ corresponding to $e_2$, 
the one labeled $(gk, \hat{h})$ 
goes from $(01234)$ to $(01245)$ corresponding to $e_3$.
Thus the orientations defined from the order on vertices match 
the convention in Fig.~\ref{rule}.

To prove part (3), we check the cases of 
interchanging the orientation of some of the edges. The general 
case will follow from the cases depicted in Fig.~\ref{phisymmetry1} 
and \ref{phisymmetry2} since the orientation changes depicted therein generate all the orientation changes.

First consider the case where the edge $(12)$ has reversed orientation.
This corresponds to changing the vertex order from
$(012345)$ in the figure to $(021345)$. (The orientations on dual edges
do not change.)
Then three colors change:
$(g, \hat{i})$ to $(g^{-1}, \widehat{ (gk)m(gk)^{-1} })$,
$(g, \hat{j} )$ to $(g^{-1}, \widehat{ (gk) \ell (gk)^{-1} })$, and 
$(g, \widehat{ khk^{-1}} )$ to $(g^{-1}, \widehat{ (gk) h (gk)^{-1} })$.
These changes are forced by the rules at faces
(uncircled trivalent vertices of the graphs). 
Hence we check the rules at circled trivalent vertices.
The only relevant vertex is the one on the left bottom in the figure.
It must hold that 
$ (gk) h (gk)^{-1} =  (gk)m(gk)^{-1} \cdot  (gk) \ell (gk)^{-1}$.
This indeed follows from $h=m \ell$.
The other cases are similar.
 $\Box$

\begin{subsect} {\bf Lemma.\/} \label{2compcolor} 
The colors define a function 
$\Psi' : {\cal F}C \rightarrow G \times \hat{G}$
where $C$ is 
a perturbed carrier
of $\Phi$ and
$ {\cal F}C $ is the set of $2$-faces of  $S$.
Conversely, a function $\Psi' : {\cal F}C \rightarrow G \times \hat{G}$
defines a color $S$ defined for the triangulation $\Phi$ and 
$\Psi=\Phi \cup \Phi^!$. 
\end{subsect}
{\it Proof.\/} 
This follows from Lemma~\ref{corresplemma} 
 and the definition of $S$. $\Box$

\begin{sect}{\bf Weighting.\/} 
\addcontentsline{toc}{subsection}{Weighting}
{\rm 
A {\it weighting} (also called a {\it Boltzmann 
weight})
is defined for each   tetrahedron, face, edge
of  triangulations $\Phi$ and $\Phi^!$ as follows.
}\end{sect}

\begin{subsect} {\bf Definition (weights for tetrahedra)\/}. {\rm
 Let $T \in \Phi$ be a 
tetrahedron 
with vertices $0$, $1$, $2$, and  $3$.
Suppose $S(01|0123)  =(m, \hat{n})$, 
$S(12|0123)=(k, \hat{l})$, and 
$S(23|0123)  =(g, \hat{h})$. 

The weight of $T$ with respect to the  given labelings of edges is
a number (an element of the ground field) defined by

$$ B(T) = B(0123) = \alpha (g, k, m ; \hat{n})  ^{\epsilon(T)} . $$

Here $\epsilon(T)$ is $\pm1$ and is defined as follows.
Let $T=(a_0, a_1, a_2, a_3)$ be the tetrahedron in consideration
where $a_0 <  a_1 < a_2 < a_3$.
Then $T$ is shared by two $4$-simplices, say, 
 $S=(a_0, a_1, a_2, a_3, v)$ and $(a_0, a_1, a_2, a_3, w)$.
Here we ignore the given 
labels 
of $v$ and $w$, 
and consider the orders written above ($v$ and $w$ coming last).
Then exactly one of these two $4$-simplices, say, $S$, with 
this order  $a_0 <  a_1 < a_2 < a_3 < v$, matches the orientation
of the $4$-manifold, and the other has the opposite orientation.

Consider the label on $v$ induced by the 
ordering on the vertices. If the 
integer index of 
$v$ is such that 
the oriented simplex $(a_0, a_1, a_2, a_3, v)$ is obtained from 
the order induced by labeling by an 
even permutation, then $\epsilon(T)=1$.
Otherwise, $\epsilon(T)=-1$.
(Sometimes we represent the order of vertices by labeling the vertices 
by  integers.)     

} \end{subsect}

\begin{subsect} {\bf Definition (weights for faces).\/} {\rm
Suppose that a face $F= (012)$ is shared by three  tetrahedra
$(0123)$, $(0124)$, and $(0125)$.
Suppose
$S_0(01|0123)=(k, \hat{\ell})$,
$S_0(12|0123)=(g, \widehat{k\ell k^{-1}})$,
$S_0(02|0123)=(gk, \hat{\ell})$,
 $S_0(01|0124)=(k, \hat{h} ) $, and
 $S_0(02|0125)=(gk, \hat{m} )$.
The situation is depicted in  Fig.~\ref{phi}.

Then
the weight for the face $F=(012)$ is defined as the number
$$ B(F) = B(012) = \phi (g, k ; \hat{m}, \hat{l})^{\epsilon(F)} . $$

Here the sign 
$\epsilon(F)= \pm 1$ 
is defined as follows.
If the local orientation defined by $(012)$ in this order
together with the orientation of the link 
$(34)
\rightarrow 
(45)
\rightarrow
(53)$
of $F$ in this order gives the same orientation as that of
the $4$-manifold, then
 $\epsilon(F)= 1$, otherwise  $\epsilon(F)= -1$.

Suppose the face $(012)$ is shared by $n$ (more than three) tetrahedra
$(012k)$, $k=3, \cdots, n$. 
Note that the vertices of these tetrahedra other than $0,1,2$ 
form a link of the face $(012)$.
Assume that these vertices are  
 cyclically
ordered by $3, 4, \cdots n$,
and that the link of $(012)$ and 
$(34)
\rightarrow 
(45)
\rightarrow 
\cdots 
\rightarrow
(n3)$ matches
the orientation of the $4$-manifold.
Suppose $S(01|0124)=(k, \hat{\ell})$,
$S_0(12|0124)=(g, \widehat{ k \ell k^{-1} } )$,
$S_0(02|0124)=(gk, \hat{\ell})$,
$S_0(01|0125)= (k, \widehat{ m_1 })$,
$S_0(01|012k)=(k, \widehat{ m_{k-4} })$, ($k=5, \cdots, n-1$),
$S_0(01|012n)=(k, \widehat{ m_{n-4} ^{-1} } )$.
Then

$ B(F)=B(012)= \phi(g,k ;  \widehat{ m_1} , \hat{ \ell} )
\phi(g, k; \widehat{ m_2 } , \widehat{ m_1 \ell} ) \cdots $

$ \phi(g, k; \widehat{ m_{n-5}} , \widehat{ m_{n-6} \cdots m_1 \ell })
\phi(g, k; \widehat{m_{n-4}},  \widehat{ m_{n-5} \cdots m_1 \ell }). $

If the cyclic order of  vertices is not as above, then it can 
be obtained from the above by transpositions.
When a transposition between $k$-th and $(k+1)$-st vertex occurs,
change the argument of $k$-th weight 
$\phi (g, k; \widehat{ m_k },  \widehat{ m_{k-1} \cdots m_1 \ell })$
to $\phi (g, k; \widehat{ m_k^{-1} },  \widehat{ m_{k-1} \cdots m_1 \ell })$.


} \end{subsect}

Notice that by the conditions in the definition of a triangulation
each face must be shared by at least three tetrahedra.
However 
In the course of computation we may have to deal with
{\it singular triangulations}.
In this case it can happen  
 that only two tetrahedra
share a face.
Let $0123$ and $0124$ be such two terahedra sharing the face $012$.
Then the weight assigned to the face $012$ in this case is
the product of Kronecker's deltas:

{\large
$$ \delta_{ S(01|0123), S(01|0124)}
 \delta_{ S(02|0123), S(02|0124)} \delta_{ S(12|0123), S(12|0124)} $$
}

There are other cases of singular triangulations that appear in our
proofs of well-definedness, and their weights are defined as follows.

Suppose that 
two tetrahedra share 
vertices $0,1,2,3$
and share 
all of their
 faces except  $(013)$. 
Meanwhile, suppose that
the face $(013)$ 
on
each of  
the respective 
tetrahedra is 
shared with 
tetrahedra  
$(0134)$ and $(0135)$. The other faces are shared by 
the tetrahedra $(0126)$, $(0237)$, and $(1238)$.
The situation is depicted in Fig.~\ref{conemove} top.
In the figure, colors and weights are also depicted. 
The signs for each weight $\phi$ are also depicted in the figure in this situation, by indicating the power $-1$ on one of the $\phi\/$s.

Suppose in another situation that two tetrahedra share all 
vertices 
and all 
$2$-faces (triangles) as depicted in Fig.~\ref{pillow}. 
The signs for this situation are also depicted in the figure.

In general if the order of vertices are different,
 then they are obtained from the 
above specific situations by compositions of permutations.
Then the weights and signs are defined by applying Lemma~\ref{6in1}.

\begin{figure}
\begin{center}
\mbox{
\epsfxsize=4in
\epsfbox{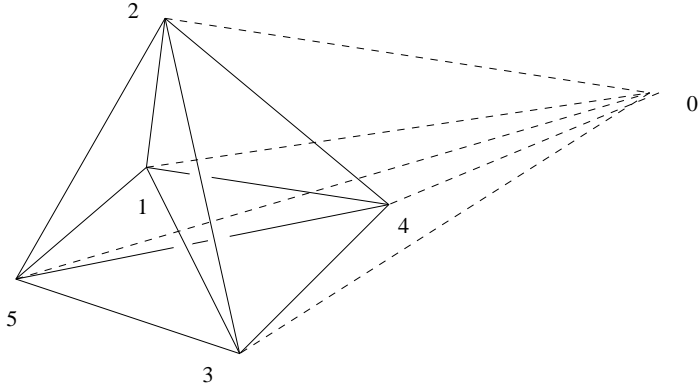}
}
\end{center}
\caption{A triangulation of the 4-ball}
\label{edgecone}
\end{figure}

\begin{subsect} {\bf Definition (weights for edges).} {\rm
Consider the triangulation of the 3-sphere $S^3$ consisting of
5 tetrahedra $(1234)$, $(1235)$, $(1245)$,
$(1345)$, 
 and $(2345)$,
where integers represent the vertices.
This triangulation is depicted in Figure~\ref{edgecone} by
solid lines. Here, the subdivided tetrahedron with vertices
$2,3,4$ and $5$ with the interior vertex $1$
is a triangulation of a 3-ball, and together with the
``outside'' tetrahedron $(2345)$ they form a triangulation of
$S^3$.
Now take a cone of this triangulation with respect to the vertex $0$
to obtain a triangulation of a $4$-ball consisiting of
5 $4$-simplices $(01234)$, $(01235)$, $(01245)$, $(01345)$,
and $(02345)$. This is depicted in Figure~\ref{edgecone} also,
where edges having $0$ as end point are depicted by dotted lines.
(Regard dotted lines as lying in the interior of the $4$-ball.)

Suppose an edge $E= (01)$ has this particular triangulation as the
neighborhood.
Suppose $S(01|0123)=(g, \hat{\ell})$,
$S(01|0124)=(g, \hat{k})$, and 
$S(01|0125)= (g, \hat{j})$.
 Then the weight for the edge $(01)$ is defined by

$$ B(E) = B(01) = \beta (g; \hat{j}, \hat{k}, \hat{\ell} )^{\epsilon (E)}  . $$

The sign $\epsilon (E) = \pm 1$
is defined in the same manner as $B(T)$ simply taking the dual
orientations.

If the neighborhood of an edge has a different triangulation, then
the weight is defined as follows.
Let $H_1, \cdots , H_s \in \Phi^!$ be the set of tetrahedra  of
the polytope $p \in \Phi^*$, and let $h_i^j$ be the set of edges
of $H_i$, $j=1, \cdots, 6$.
Then the weight is defined by

$$ B(E) = 
1/|G|^{2a} 
\sum \prod \beta(g; \widehat{j_f}, \widehat{k_f},\widehat{\ell_f})
$$
where each $\beta$ is assigned to a tetrahedron of the above 
triangulation following the order convention of vertices.
The product of the above expression is taken over all the shared
edges, and the sum is taken over all the possible states on shared
  edges. 
The exponent, $a$, on the normalization factor, $1/|G|^{2a}$, 
is the number of verticies in the interior of the polyhedron 
dual to the given edge. 

} \end{subsect}

\begin{figure}
\begin{center}
\mbox{
\epsfxsize=6in
\epsfbox{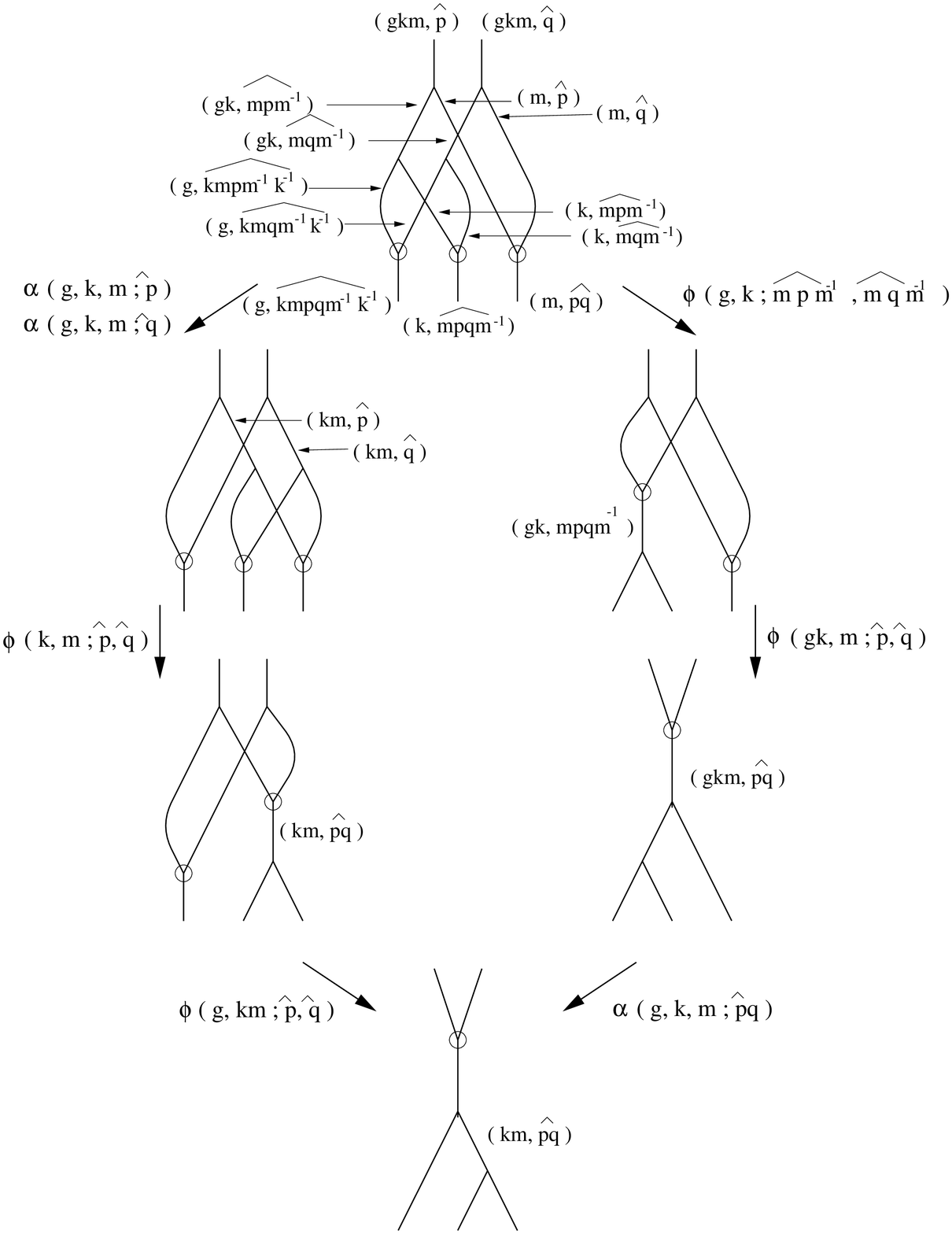}
}
\end{center}
\caption{Movies of cocycle trees, Part I}
\label{cocymovie1}
\end{figure}

\begin{figure}
\begin{center}
\mbox{
\epsfxsize=5in
\epsfbox{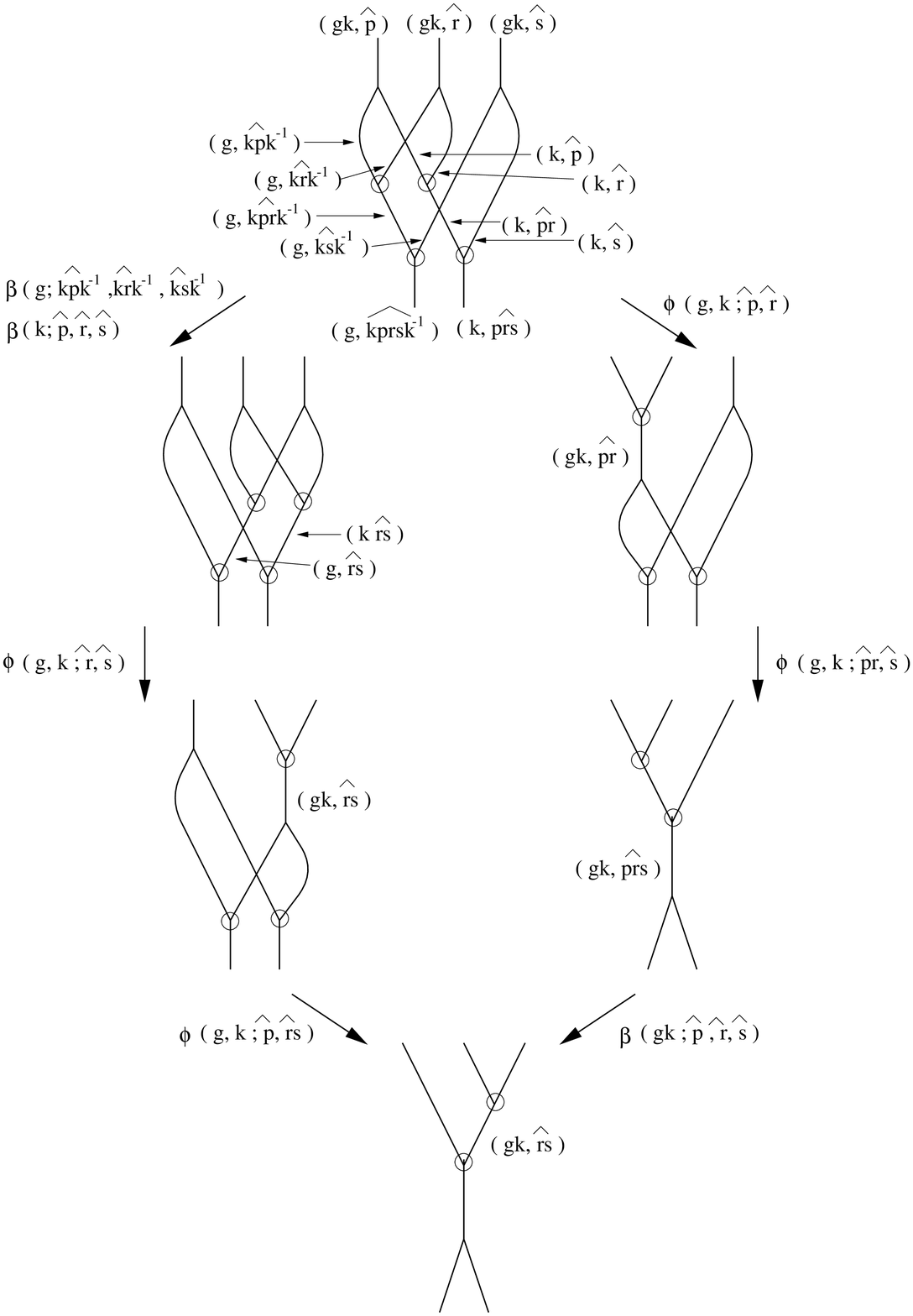}
}
\end{center}
\caption{Movies of cocycle trees, Part II}
\label{cocymovie2}
\end{figure}

\begin{sect}{\bf Partition function.\/}      
\addcontentsline{toc}{subsection}{Partition function}
{\rm
Let $\Phi$ be a triangulation  of a $4$-manifold $M$
with the set of vertices (resp. edges, faces, tetrahedra)
${\cal V}$ (resp. ${\cal E}$, ${\cal F}$, ${\cal T}$).
Fix also a triangulation $\Phi^!$ of the dual $\Phi^*$.

}\end{sect}
\begin{subsect} {\bf Definition.} {\rm
The partition function $\psi (\Psi) $ for a triangulation $\Psi=
\Phi \cup \Phi^!$
with a total order on vertices is
defined by

$$ \psi(\Psi) = 
1/|G|^{2a}
\sum_S \prod_{ {\cal T}, {\cal F}, {\cal E} } B(T) B(F)
B(E) $$
where the product ranges over tetrahedra, faces and edges of the
triangulation $\Phi$, 
the summation ranges over all the possible
states, 
and the exponent $a$ on the normalization 
factor is the number of vertices
in the triangulation.

} \end{subsect}

\begin{subsect} 
{\bf Main Theorem.}
 \label{ftn} 
The partition function $\psi(\Psi)$ defined above
for triangulations $\Psi = \Phi \cup \Phi^!$ of a 4-manifold $M$
is independent of the choice of the triangulation $\Phi$ and $\Phi^!$
and independent of choice of order on vertices.

Therefore the partition function  $\psi$ defines an invariant of
a 4-manifold $M$.
\end{subsect}

Section~\ref{invsec}
is devoted to 
giving the proof of this theorem.

\begin{sect}{\bf Diagrams, cocycles, and triangulations.\/}
\addcontentsline{toc}{subsection}{Diagrams, cocycles and triangulations}
\label{dc&t}
{\rm Here we explain relations among diagrams, cocycles, and triangulations.
Figure~\ref{rule} illustrates the coloring rules at triangles
 and dual triangles. In these triangles and dual triangles graphs
 are embedded; the verticies of the graphs in the dual triangles are
 labeled by small circles. 
The cocycles $\alpha$ are assigned to tetrahedra; 
the cocycles $\beta$ are assigned to dual tetrahedra,
 and the cocycles $\phi$ are assigned to triangular faces.
 Each such figure also corresponds to a graph movie which depicts
 a part of the perturbed carrier surface.
 We can think of the cocycles as being assigned to the vertices
 of the perturbed carrier surface which has a tripartition on its vertex set.
 Indeed the 
many 
scenes that constitute the graph movie are found on the boundary 
of a regular neighborhood 
of the vertices of the carrier surface.
 In this way we can directly visualize the construction of the invariant
 as a colored surface with weighted vertices or as a colored graph movie
 with weights associated to the scenes.

Simliarly, 
the cocycle conditions can be
described 
as relations on movies of tree diagrams.
Figures~\ref{cocymovie1} and \ref{cocymovie2} 
depict  
these relations.
Each change of  a tree diagram 
(scene in  the movie) 
corresponds to a cocycle as indicated. 
When we multiply the left-hand-side and the right-hand-side
 of cocycles in the movies,
we obtain cocycle conditions among $\alpha$, $\phi$, and $\beta$.

The cocycle conditions can also be understood in terms of certain 
singular surfaces that are embedded in the 4-manifold.
 These surfaces are depicted
in Figs.~\ref{surfcubeLHS} and \ref{surfcubeRHS}. 
In these figures the cocycles $\phi$ corresponds to the surface $(Y\times Y)$
 and the cocycle
$\alpha$ corresponds to the surface  that is dual to a tetrahedron. 
The assignment of $\phi$ to $Y\times Y$ is indicated in the weights
 on Fig.~\ref{YxY}. The reasons for these assignments is that the
 cocycle $\phi$ is found when three tetrahedra share a triangular face,
 and the cocycle $\alpha$ is assigned to a tetrahedra. 
In Figs.~\ref{surfcubeLHS} and \ref{surfcubeRHS} some edges are denoted
 as tubes. A tube of the form $\circ \times Y$ corresponds to a triangle
 that is shared  by three tetrahdra as in Fig.~\ref{face}.

\begin{figure}
\begin{center}
\mbox{
\epsfxsize=4in 
\epsfbox{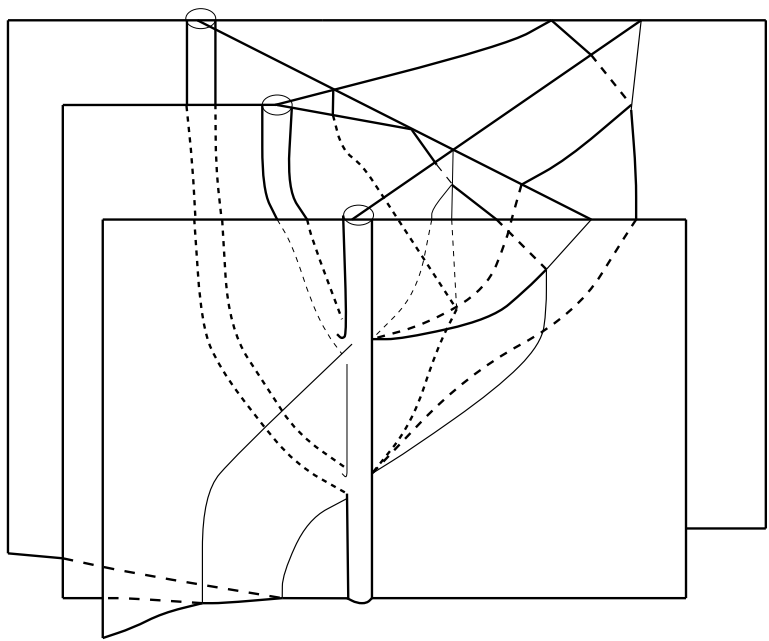}
}
\end{center}
\caption{Surface of cocycle movies: left hand side}
\label{surfcubeLHS}
\end{figure}

\begin{figure}
\begin{center}
\mbox{
\epsfxsize=4in 
\epsfbox{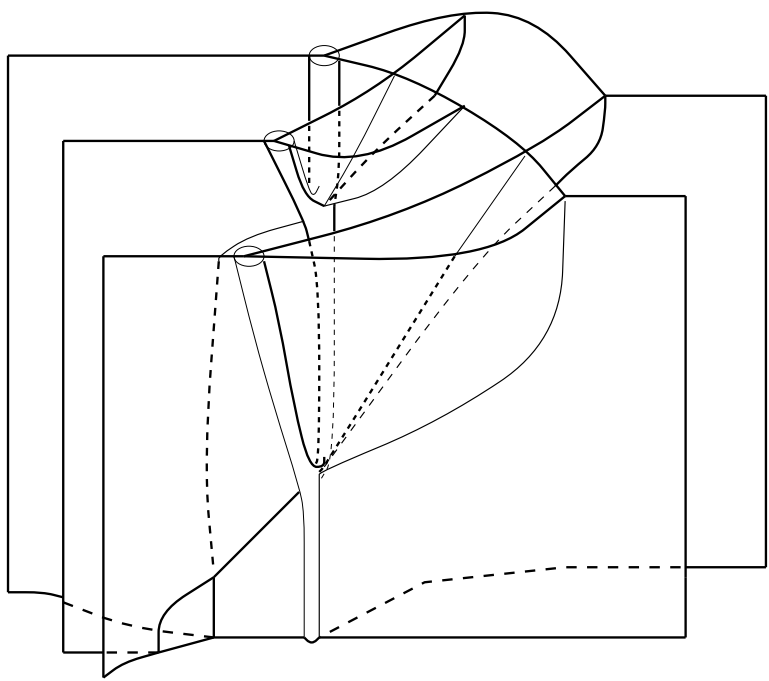}
}
\end{center}
\caption{Surface of cocycle movies: right hand side}
\label{surfcubeRHS}
\end{figure}

}\end{sect}

\section{On invariance of the partition function} \label{invsec}  
 
Recall the notation 
in Section~\ref{bzman}:
$\Phi$ denotes a triangulation of a 4-manifold $M$,
$\Phi^*$ its dual complex,
$\Phi^!$ a 3-face triangulation of $\Phi^*$.
In Section~\ref{ordersec} we show that the partition function defined 
is independent of the order on vertices. In 
Section~\ref{invP},  
 we show that the partition function is independent of the triangulation. In 
Section~\ref{nodual} we show that the partition function is independent
 of the choice of dual triangulation.

\begin{figure}
\begin{center}
\mbox{
\epsfxsize=5in 
\epsfbox{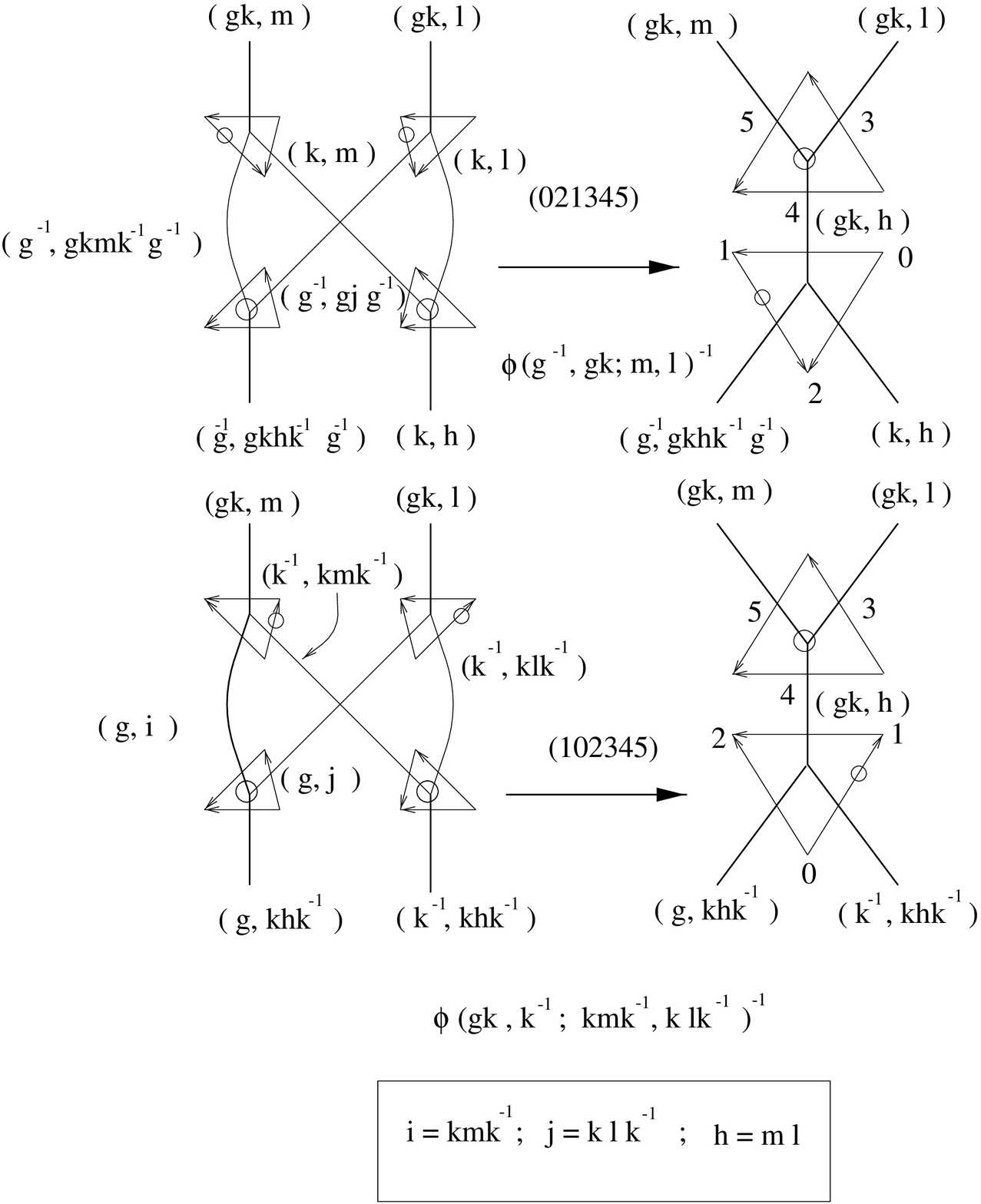}
}
\end{center}
\caption{Symmetries of $\phi$, Part I}
\label{phisymmetry1}
\end{figure}

\begin{figure}
\begin{center}
\mbox{
\epsfxsize=5in 
\epsfbox{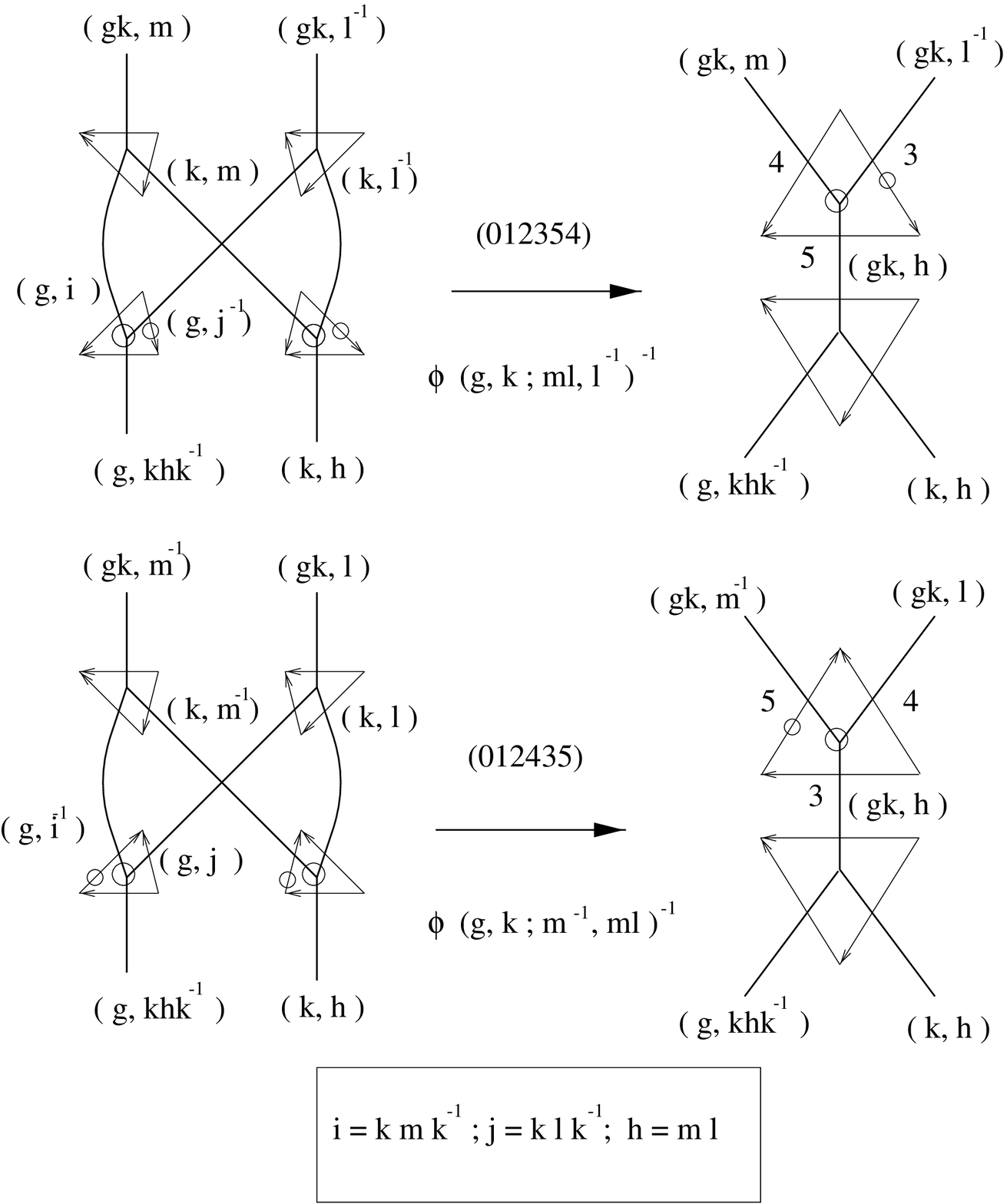}
}
\end{center}
\caption{Symmetries of $\phi$, Part II}
\label{phisymmetry2}
\end{figure}

\begin{sect} {\bf Independence on order of vertices.}
\addcontentsline{toc}{subsection}{Independence on order of vertices}
 \label{ordersec} {\rm 
In this section we prove 
} \end{sect}

\begin{subsect} {\bf Lemma.\/} 
\label{invlemma1} 
The cocycle symmetries imply the
independence of the partition function
on the  order on vertices of the triangulation.
\end{subsect}
{\it Proof.\/}
For tetrahedra and dual tetrahedra, 
the weights are the cocycles $\alpha$ and $\beta$ 
respectively. 
As in \cite{Wakui}, it is sufficient to
check how weights change when the order of vertices 
are changed from $(0123)$ to 
$(0132)$, $(0213)$, and $(1023)$. 
Such changes are 
illustated in Fig.~\ref{alphasymmetry}
for $\alpha$, and 
the corresponding conditions are listed for $\alpha$, $\beta$ and $\phi$
in Section~\ref{cocysym}.

For a face,
 we check as follows.
In Fig.~\ref{phi} an  order of vertices are given, where 
the face is given by  $(012)$, and the
other vertices are given 
labels $3$, $4$, and $5$.
The changes of orders of vertices are generated by 
the changes from $(012345)$ to
$(021345)$, $(102345)$, (for the face)
$(012354)$, $(012435)$, (for the dual face)
since only the relative orders among the vertices of the face  and
those of the  dual faces 
 are in consideration.
For these changes, the colors are listed in Fig.~\ref{phisymmetry1}
and Fig.~\ref{phisymmetry2}.
They are depicted in terms of dual graphs, and on the right hand side,
the orientations of edges of faces/dual faces are shown. 
The small circles indicate reversed orientations.
The corresponding conditions are listed in Section~\ref{cocysym}.

If more than three tetrahedra share a face, 
then a change in order of the vertices can be achieved by such pairwise
  switches. Futhermore, in order to affect such changes,
 we may have to group the vertices in sets of 3. This grouping is achieved by
 a 3-face triangulation.
 So the proof will follow once we have shown invariance
 under the 3-face triangulation.
$\Box$

\begin{figure}
\begin{center}
\mbox{
\epsfxsize=5in
\epsfbox{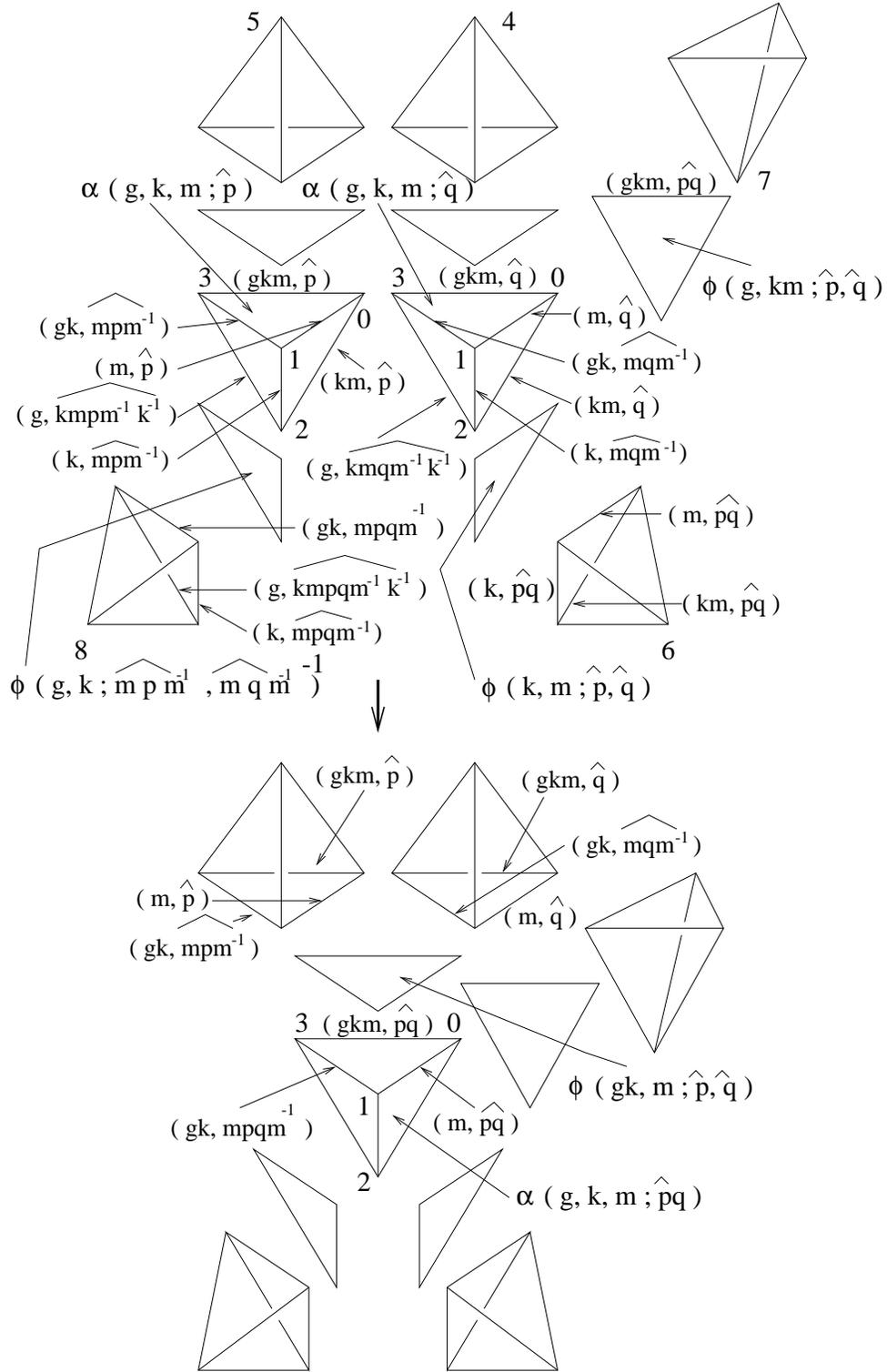}
}
\end{center}
\caption{Colors and cocycles for the cone move}
\label{conemove}
\end{figure}

\begin{sect}{\bf Independence under Pachner moves.}
\addcontentsline{toc}{subsection}{Independence under Pachner moves}
\label{invP}
{\rm
In this section,
we explicitly relate the cone move, taco move
 and pillow move to the cocycle conditions.
 Since these moves and lower dimensional moves generate the Pachner moves,
 we will use the cocycle conditions to show that the partition
 function is invariant under 
the Pachner moves.
}\end{sect}

\begin{subsect}{\bf Lemma.  } \label{coneinv}
The partition function is invariant under
the cone move for a
local triangulation with 
a specific choice of order depicted in Fig.~\ref{conemove}.
 \end{subsect}
    {\it Proof.}
Let $(0123)_1$ and $(0123)_2$ be tetrahedra sharing the same faces
$(012)$, $(013)$ and $(023)$, but having different
faces $(123)_1$ and $(123)_2$, such that
(1) 
the union of the triangles
$(123)_1$
$\cup$
$(123)_2$ bounds 
a 3-ball $B$ in
the 4-manifold,  (2) the union of $B$ , $(0123)_1$
and $(0123)_2$ is diffeomorphic to the $3$-sphere
bounding a $4$-ball $W$ in the 4-manifold.
(See 
Figs.~\ref{4dcone}, \ref{conemove}.)
In these figures, the movies of dual graphs are depicted
where 
each of the faces $(013)$, $(023)$, $(123)$ is shared by another 
tetrahedron 
($(0124)$, $(0125)$, $(0126)$, respectively).
We prove the invariance in this case. 
The general case follows from
such computations together with the pentagon identity of $\beta$.

Fig.~\ref{conemove} shows the colors and cocycles assigned to
this local triangulation (again note the direct relation
between this assignment and those for the top graph in Fig.~\ref{cocymovie1}).
The left hand side of the cone move (top of Fig.~\ref{conemove})
has the local contribution 
$$\phi (g, km; \hat{p}, \hat{q} )
\phi ( k, m; \hat{p}, \hat{q} )
\phi (g, k; \widehat{mpm^{-1}}, \widehat{mqm^{-1}} )^{-1}
\alpha (g, k, m; \hat{p})
\alpha (g, k, m; \hat{q}), $$

(note that the orientation of the face $(9123)$ is opposite),
and the right hand side of the cone move (bottom of the figure)
has the local contribution
$$\phi (gk, m; \hat{p}, \hat{q} )
\alpha (g,k,m;\widehat{pq} ).$$
Thus the partition function is invariant 
under the cone move because the cocycle condition is satisfied.
$\Box$

\begin{figure}
\begin{center}
\mbox{
\epsfxsize=5in 
\epsfbox{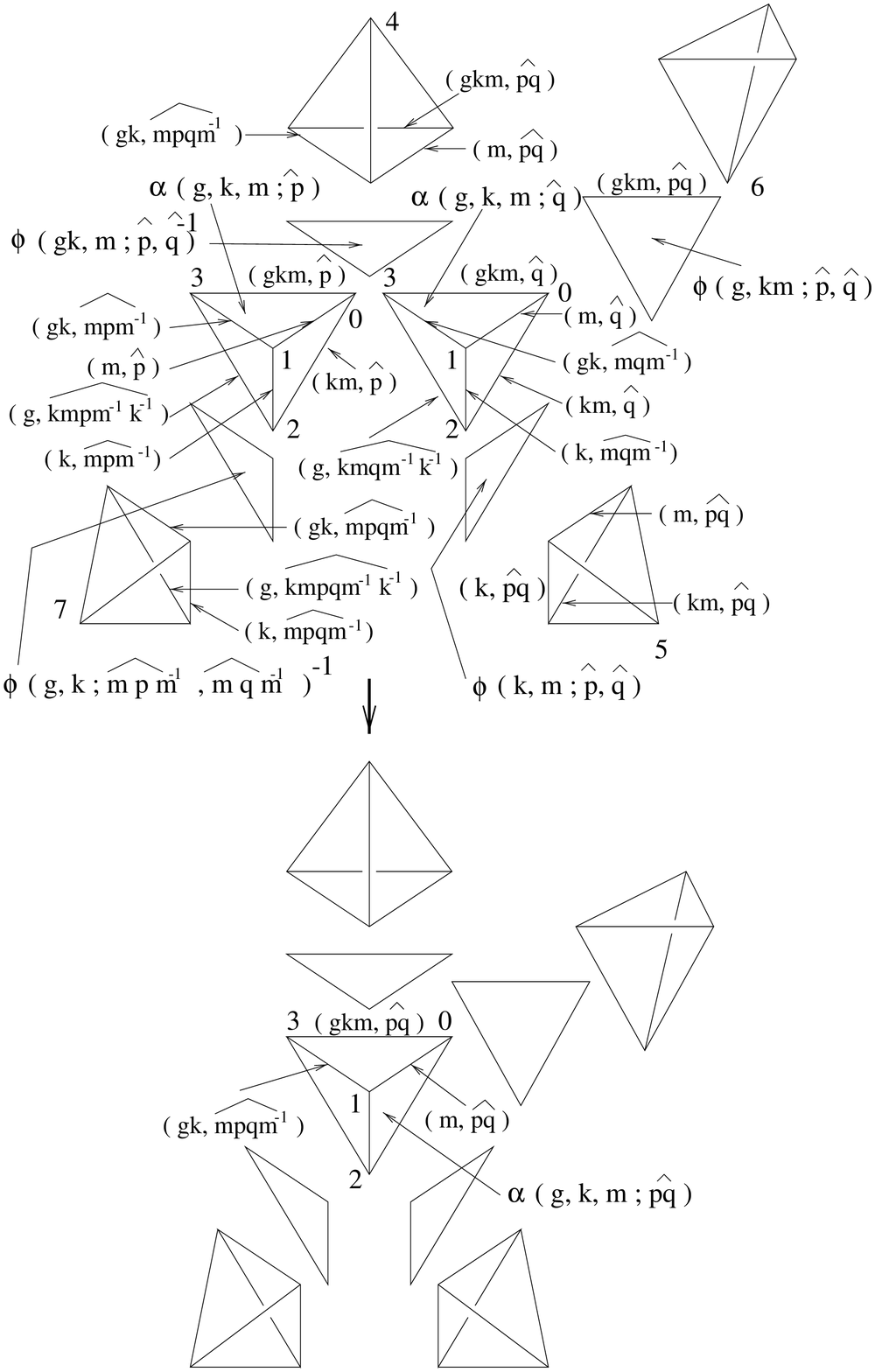}
}
\end{center}
\caption{Colors and cocycles for the pillow move}
\label{pillow}
\end{figure}

\begin{figure}
\begin{center}
\mbox{
\epsfxsize=3.5in 
\epsfbox{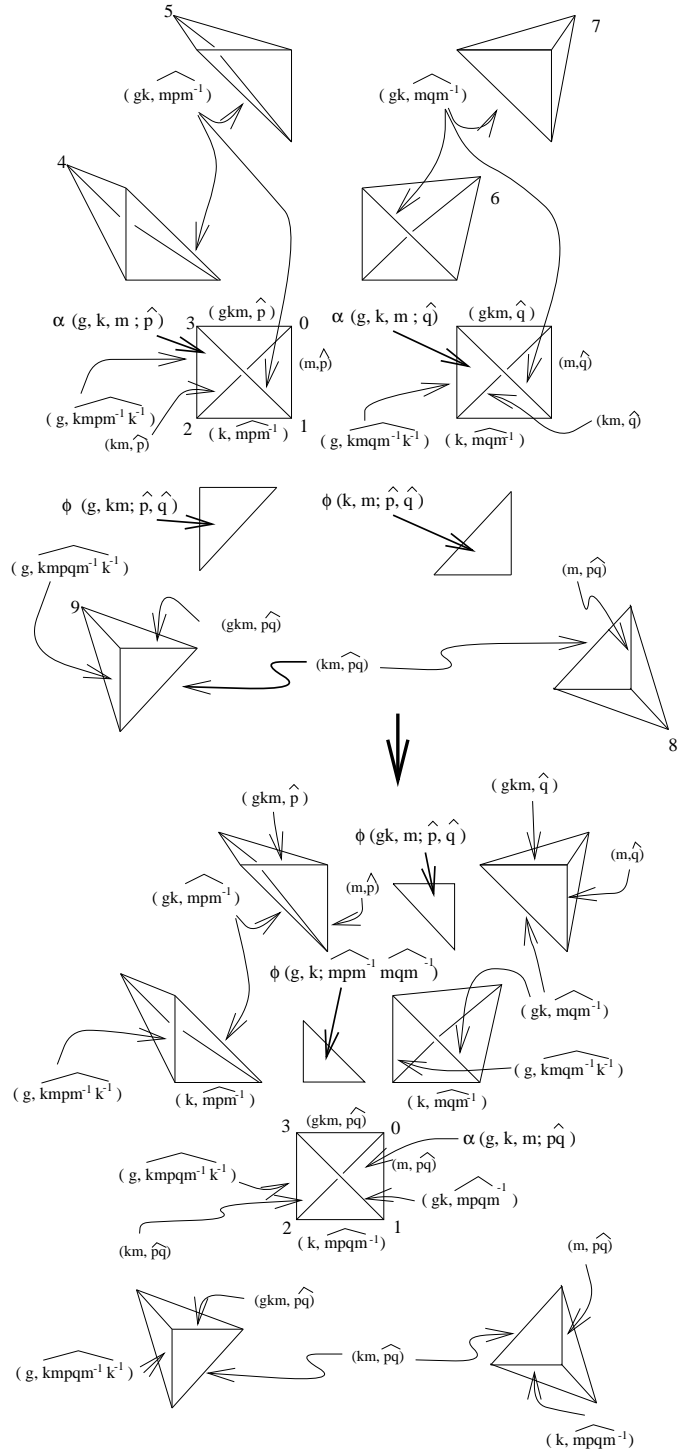} 
}
\end{center}
\caption{Colors and cocycles for the taco move}
\label{taco}
\end{figure}

\begin{subsect} {\bf Lemma.} 
The partition function is independent 
under the pillow move for a specific local triangulation 
with the order depicted in Fig.~\ref{pillow}.
\end{subsect}
    {\it Proof.}
In Fig.~\ref{pillow} 
the 
assignments  
 of colors and cocycles are 
shown. 
The left hand side of the pillow move (top of Fig.~\ref{pillow})
has the local contribution 
$$\phi (g, km; \hat{p}, \hat{q} )
\phi ( k, m; \hat{p}, \hat{q} )
\phi (gk, m; \hat{p}, \hat{q} )^{-1}
\phi (g, k; \widehat{mpm^{-1}}, \widehat{mqm^{-1}} )^{-1}
\alpha (g, k, m; \hat{p})
\alpha (g, k, m; \hat{q}), $$
and the right hand side of the pillow move (bottom of the figure)
has the local contribution
$\alpha (g,k,m;\widehat{pq} ).$
This follows from the cocycle condition used in the above lemma.
$\Box$

\begin{subsect} {\bf Lemma.}
The partition function is independent 
under the taco move for a specific local triangulation 
with the order depicted in Fig.~\ref{taco}.
\end{subsect}
    {\it Proof.}
In Fig.~\ref{taco}  
the 
assignments 
 of colors and cocycles are 
shown. 
The left hand side of the
taco 
move (top of Fig.~\ref{taco})
has the local contribution 
$$\phi (g, km; \hat{p}, \hat{q} ) \phi ( k, m; \hat{p}, \hat{q} )
\alpha (g, k, m; \hat{p})
\alpha (g, k, m; \hat{q}), $$
and the right hand side of the 
taco 
move (bottom of the figure)
has the local contribution
$$\phi (gk, m; \hat{p}, \hat{q} )
\phi (g, k; \widehat{mpm^{-1}}, \widehat{mqm^{-1}} )
\alpha (g,k,m;\widehat{pq} ).$$
This  is exactly one of the cocycle conditions. 
$\Box$

Observe that the diagrammatics of the graph movie move that results from the 
taco move match exactly the graph movie move that represents the cocycle
 condition Fig.~\ref{cocymovie1}. Similar graph movies can be drawn for the 
cone and pillow moves and the correspondence with the move and the coycle
 conditions can be worked out via the graph movies.
 Making such correspondence shows explicitly the method of constructing
 invariants via Hopf categories where, instead of cocycle conditions,
 coherence relations are used. The coherence relations can be expressed
 by such graph movie moves
(See Section~\ref{Hopfsec}).

Since the partition function is invariant under the cone, taco, 
and pillow moves, and since 
$\alpha$ satisfies a pentagon relation, we have the partition function
 is invariant under the Pachner moves.

\begin{sect}{\bf Independence on triangulations of the dual complexes.\/}
\addcontentsline{toc}{subsection}{Independence on triangulations of the dual complexes}
\label{nodual}
 {\rm 
In this section, we complete the proof that the partition function is 
well-defined by showing that the partition function does not depend 
on the 3-face triangulation, $\Phi^!$.   } 
\end{sect}

\begin{subsect}{\bf Lemma.\/}
\label{relpaclemma} 
If $T_1$ and $T_2$ are triangulations
of a $3$-dimensional polytope which is diffeomorphic to
a $3$-ball such that  $T_1$ and $T_2$ restrict to the same
triangulation on the boundary, then they are related by a finite sequence of 
Pachner moves. 
 \end{subsect} 
  {\it Proof.} We first prove the corresponding statement in dimension 2, then use the result in dimension 2 to achieve the result in dimension 3. 
In the proof we use the notation 
$(i \rightleftharpoons j)$-move 
to indicate the move in which $i$ simplices are replaced by 
$j$ simplices. So the  
$(j \rightleftharpoons i)$-move 
is the inverse move, and the order of $i$ and $j$ matter.

In dimension 2, we have two triangulations of the disk that agree on the boundary, and we are to show that they they can be arranged by Pachner moves fixing the boundary to agree on the interior. We prove the result by induction on the number of vertices

on the boundary.

Recall \cite{RS} \cite{BW}, that the {\it star of a $k$-simplex}
(in a simplicial complex) is the union of all the simplices that contain the 
$k$-simplex. The {\it link of a $k$-simplex} is the union of all the simplices in the star that do not contain the $k$-simplex. We will examine the stars and links of vertices on the boundary of a disk (and later on the bounday of a 3-ball). So denote the

 star of $v$ with respect to the boundary by: ${\mbox{\rm st}}_{S} (v)$. Similarly, the link of $v$ with respect to the boundary is ${\mbox{\rm lk}}_{S} (v)$ while these sets with respect to the interior are ${\mbox{\rm st}}_{B} (v)$ and ${\mbox{\rm lk}}_
{B} (v)$, respectively.

In dimension 2, ${\mbox{\rm st}}_{S} (v)$ is a pair of edges that share 
the vertex $v$. Meanwhile, ${\mbox{\rm lk}}_{B} (v)$ is a polygonal path properly embedded in the disk that is the  most proximate to  $v$ among all paths in the interior that join the points of ${\mbox{\rm lk}}_{S} (v)$.

We fix our consideration on one of the triangulations, say $T_1$ of $D^2$.
We want to alter this triangulation so that ${\mbox{\rm lk}}_{B} (v)$
is an edge (so it has no interior vertices). If we can achieve this alteration, then we can perform similar moves to $T_2$. The vertex $v$ on either triangulation then will become the vertex of a triangle that is attached to the disk along a single edge.

We can remove such a triangle 
(or alternatively, work in the interior) and apply induction
on the number of vertices 
on the boundary. 

Consider an interior vertex, $v'$ in $B$. 
If the star of $v'$ in $B$ is the union of three triangles at $v'$, then  
we can remove this vertex from $D$ by means of a 
$ ( 3 \rightleftharpoons 1)$-move. 
Perform such moves until there are no interior vertices of valence 3. 
In this way we may assume that a vertex, 
$v'' \in {\mbox{\rm lk}}_{B} (v)$ has valence larger than 3. If the valence of $v''$  is greater than $3$, then there are a pair of triangles 
in ${\mbox{\rm st}}_{B} (v)$
sharing edge $v,v''$ upon which a Pachner move of type $ ( 2 \rightleftharpoons 2)$ can be performed. 
Such a move removes $v''$ from the link of $v$. 
After such a move, check for interior vertices of valence 3 and remove them by type $ ( 3 \rightleftharpoons 1)$-moves.
In this way we can continue until the link of $v$ is an edge.
If $D$ is a triangle, then the process will reduce the triangulation until there are no interior vertices.

Now we mimic the proof given in dimension 2, to dimension 3. First,  assume that an interior vertex $v'$ in $D^3$ has as its star the union of 4 tetrahedra. Then we may 
eliminate such an interior vertex by means of a type 
$ ( 4 \rightleftharpoons 1)$ Pachner move.

Consider the link, ${\mbox{\rm lk}}_{B} (v)$ of a vertex, $v$, on the boundary. If this link is a triangle, then we may eliminate the vertex from the boundary, as in the 2-dimensional case. For the star of $v$ is a single tetrahedron that is glued to the
 ball along a single face.

More generally, ${\mbox{\rm st}}_{S} (v)$ is a union 
of triangles forming a polygon, 
 so ${\mbox{\rm st}}_{S} (v)$ is the cone on 
the polygon ${\mbox{\rm lk}}_{S} (v)$ where $v$ is the cone point.
Consider the disk properly 
embedded in 
$B$ that is the link of $v$.
This link, ${\mbox{\rm lk}}_{B} (v)$, is a triangulated disk.
There is a sequence of 2-dimensional Pachner moves that change 
${\mbox{\rm lk}}_{B} (v)$ to a triangulation of an $n$-gon,
with no interior vertices. We use these 2-dimensional moves to determine 
$3$-dimensional moves performed in a neighborhood of ${\mbox{\rm st}}_{B} (v)$ as follows.

Suppose that a $ ( 3 \rightleftharpoons 1)$-move is used to simplify the disk that is the link of $v$. Then consider the vertex $v'$ at which such a move is performed. By our first step, 
its star is not the union of 4 tetrahedra. 
Three tetrahedra intersect along the edge, $v,v
'$, and a $ ( 3 \rightleftharpoons 2)$-move can 
be performed in the
star of $v$ to remove the 
vertex $v'$ from the link. After such a move, then
check for vertices in the interior whose 
valence is $4$. Remove these by $ ( 4 \rightleftharpoons 1)$-moves, 
until no such vertices remain. Potentially, some 
vertices from the link of $v$ are removed, and the 
effect of such a removal on the link is to perform a 
$ ( 3 \rightleftharpoons 1)$-move. In general a 
$ ( 3 \rightleftharpoons 1)$-move 
to the link corresponds to a 
$ ( 3 \rightleftharpoons 2)$-move to the star, or a 
$ ( 4 \rightleftharpoons 1)$-move to a 
part of the star and a tetrahedron on the other side of the link.

If a 
$ ( 2 \rightleftharpoons 2)$-move 
is used to simplify 
${\mbox{\rm lk}}_{B} (v)$, then there is either 
a $ ( 3 \rightleftharpoons 2)$-move 
or a $ ( 2 \rightleftharpoons 3)$-move to that 
 ball which induces it. 
Specifcally, if an edge in ${\mbox{\rm lk}}_{B} (v)$ 
has as its star the union of 3 tetrahedra, 
then two of these are found in ${\mbox{\rm st}}_{B} (v)$ 
and the other one is on the other side of the link of $v$. In this case perform a $ ( 3 \rightleftharpoons 2)$-move to $B$. The 
link of $v$ changes 
by a $ ( 2 \rightleftharpoons 2)$-move. If the link of the edge to be changed is more than 3 tetrahedra, then perform a 
$ ( 2 \rightleftharpoons 3)$-move to 
${\mbox{\rm st}}_{B} (v)$. 
In this way, a triangulation always results from these moves.
After each such move, one must go and check for vertices of valence $4$ and remove them by $ ( 4 \rightleftharpoons 1)$-moves. 

Eventually, we can remove all interior 
vertices from 
${\mbox{\rm lk}}_{B} (v)$ and we can 
further make sure that the 
link of $v$ is in some standard position. We can remove $v$ from the boundary of $B$ by removing $v$ 
and the $(n-2)$-tetrahedra in its 
star where $n$ is the valence of $v$ with 
repsect to the boundary. The result follows by induction. $\Box$

\begin{figure}
\begin{center}
\mbox{
\epsfxsize=4in
\epsfbox{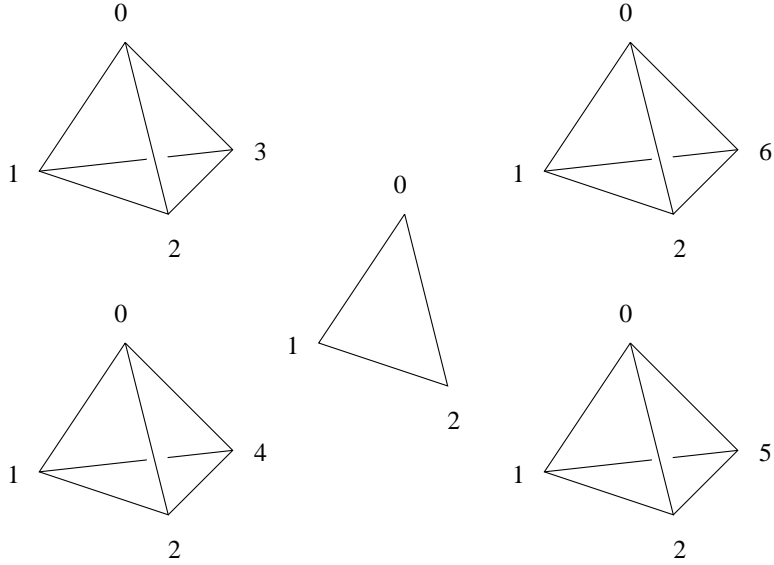}
}
\end{center}
\caption{A face sharing four tetrahedra}
\label{fourface}
\end{figure}

\begin{figure}
\begin{center}
\mbox{
\epsfxsize=4in
\epsfbox{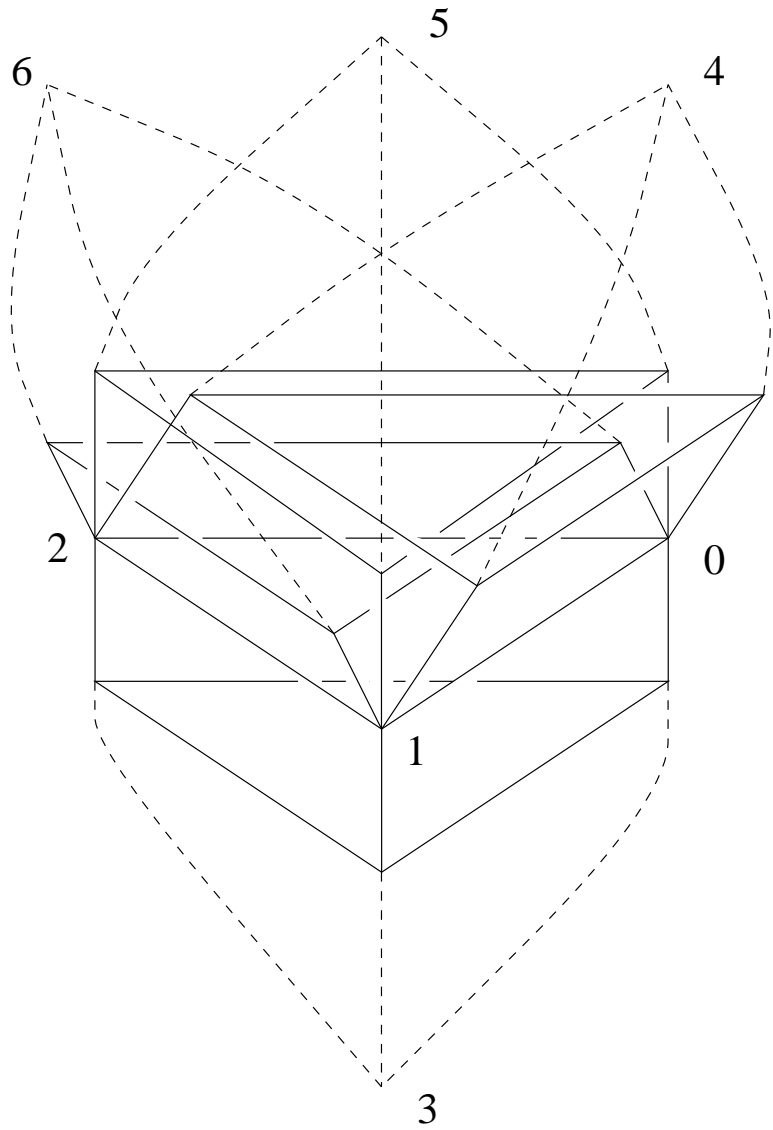}
}
\end{center}
\caption{A face sharing four tetrahedra, another view}
\label{4faces}
\end{figure}

\begin{figure}
\begin{center}
\mbox{
\epsfxsize=3.5in 
\epsfbox{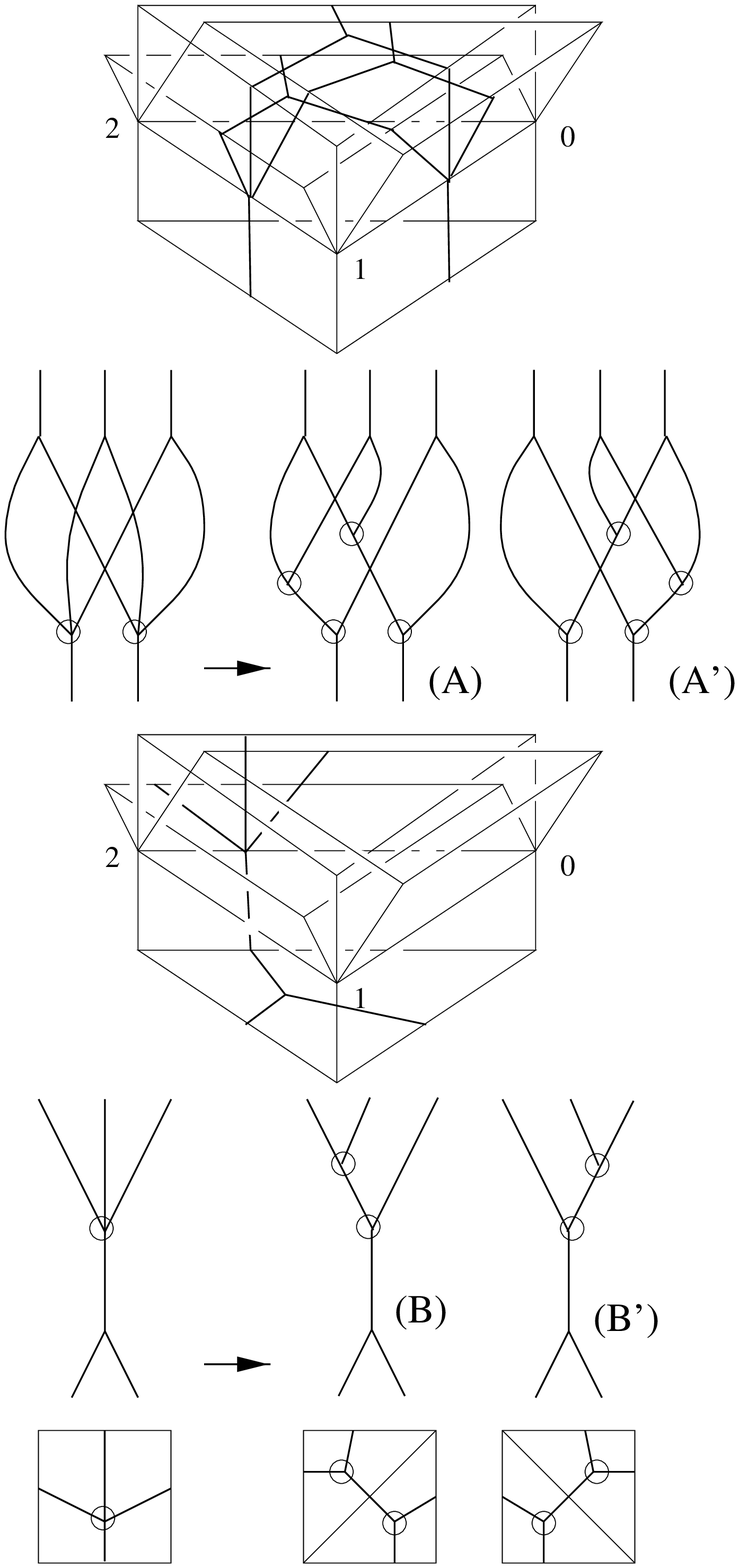}
}
\end{center}
\caption{Dual graphs around a face sharing four tetrahedra}
\label{4faces1}
\end{figure}

\begin{subsect}{\bf Lemma.} \label{inv3d} 
The partition function $\psi$  does not depend on
the choice of $3$-face triangulations of $\Phi^*$.
\end{subsect}
  {\it Proof.}
First let us analyze the case when a face $(012)$ is shared by 
four tetrahedra.
Then we will discuss 
the 
general case.
Figures~\ref{fourface} and \ref{4faces} depict the case where 
a face  $(012)$ is shared by four tetrahedra 
$(0123)$, $(0124)$, $(0125)$ and $(0126)$.


Then the dual complex has a rectangular $2$-face  $(012)^*$ which is dual to 
the face $(012)$. 
There are two triangulations of a rectangle, say $t_1$ and $t_2$,
for  $(012)^*$. 
(These are the triangulations that have no interior vertices).
The $3$-polytopes in $\Phi^*$ that share  $(012)^*$
are duals $(01)^*$, $(02)^*$ and $(12)^*$. 
Let $T_1$ and $T_2$ be 
$3$-face triangulations of $\Phi^*$ that restrict 
to $t_1$ and $t_2$
respectively and restrict to the same triangulation
on all the other $2$-faces of $\Phi^*$.
We show that the partition functions defined from $T_1$ and $T_2$
give the same value.

Recall (Fig.~\ref{3ass}) that two pairs of faces of a tetrahedron give
two triangulation of a rectangle.
We attach a tetrahedron in between 
$(012)^*$ and  $(01)^*$ and change
the triangulations on the face. 
More specifically, attach a pair of adjacent faces of a 
tetrahedron 
onto the  face $(012)^*$ of   $(01)^*$
along the triangulation $T_1$ restricted to $t_1$.
Perform the same attachment for  $(02)^*$ and $(12)^*$. 
Then 
we get 
a new triangulation $T_1 '$ of $\Phi^*$ 
which restricts to $t_2 $ on  $(012)^*$.
Thus $T_1'$ and $T_2$ have the same triangulation on the $2$-skeleton
of $\Phi^*$ by 
the 
assumption 
that all the other faces have the same triangulation.  
Thus  $T_1'$ and $T_2$ are related by a finite sequence of
Pachner moves fixing the boundary triagulation by Lemma~\ref{relpaclemma} 
which does not change the partition function by the pentagon identity
and the orthogonality of the cocycle $\beta$.
Hence it remains to prove that $T_1$ and $T_1'$ 
give the same value of the partition function.

If the dual face $(012)^*$ is in general a polygon of more than four faces
(say $n$-gon),
then triangulations consist of $(n-2)$ triangles
(by the condition of the definition of $3$-face triangulation).
Such trianglations are related by only 
$(2\rightleftharpoons 2)$-moves,
which are realized by attaching a pair of faces of tetrahedra 
once at a time as above. Thus the above argument is applied
to general cases by repeating the argument.

Now we prove that  $T_1$ and $T_1'$ in the case $(012)^*$ is a rectangle 
give the same value of the partition function.
Figure~\ref{4faces1} depicts the graphs for the triangulation.
In Fig.~\ref{4faces1} the perturbations of these
graphs to trivalent graphs are also depicted. These perturbations 
correspond to $\Psi^!$ (triangulations of a rectangle $(012)^*$
in this case) as depicted in the bottom of Fig.~\ref{4faces1}.
Thus the colors assigned near the face $(012)$ with triangulations 
$T_1$ and $T_1'$ are also assigned to edges of the perturbed graphs
(the right pictures of arrows in the figure,
marked (A), (A'), (B), and (B')).

These graphs are identified by the following graphs in Fig.~\ref{cocymovie2}:
(A) corresponds to top graph, (A') to top left, 
(B) to bottom right, (B') to bottom.
The weights $\phi$ assigned to each triangulation are  thus 
$\phi (g,k; \hat{p}, \hat{r}) \phi (g,k ; \widehat{pr}, \hat{s})$
for $T_1$ and
$\phi (g,k; \hat{r}, \hat{s}) \phi (g,k ; \hat{p}, \widehat{rs},)$
for $T_1'$ (or vice versa),
if the group elements assigned are as indicated in Fig.~\ref{cocymovie2}.

Now since $T_1'$ is obtained from $T_1$ by attaching three 
tetrahedra, and they receive the weights
$\beta(g;  \widehat{kpk^{-1}}, \widehat{krk^{-1}}, \widehat{ksk^{-1}})$,
$\beta(k; \hat{p}, \hat{r}, \hat{s})$,
$\beta(gk; \hat{p}, \hat{r}, \hat{s})$, 
therefore the cocycle condition
depicted in Fig.~\ref{cocymovie2} shows that they are equal.
$\Box$

\section{Hopf categories} \label{Hopfsec}

\begin{sect} {\bf Overview of Hopf categories.\/} 
\addcontentsline{toc}{subsection}{Overview of Hopf categories}
{\rm 
We refer the reader to \cite{CF,CYex,Neu} for details 
about Hopf categories.
Here we give a brief overview.
A Hopf category is a categorification of a Hopf algebra.
Roughly speaking, this means that the multiplication between
elements becomes 
the tensor product between two objects, and there is a cotensor product
defined for objects. 

First we recall $2$-categories from \cite{KeSt,KV,BN}.
} \end{sect}

\begin{figure}
\begin{center}
\mbox{
\epsfxsize=5in
\epsfbox{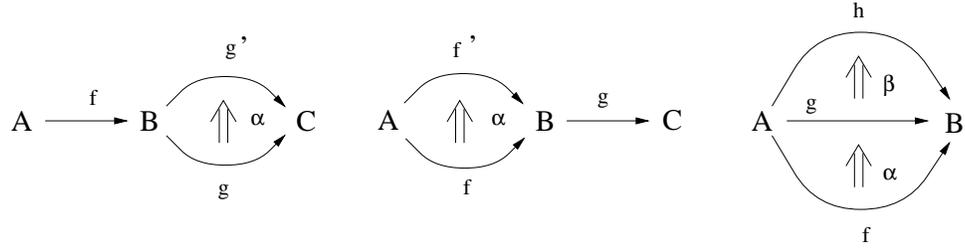}
}
\end{center}
\caption{Compositions of $2$-morphisms}
\label{compos}
\end{figure}

\begin{figure}
\begin{center}
\mbox{
\epsfxsize=3in
\epsfbox{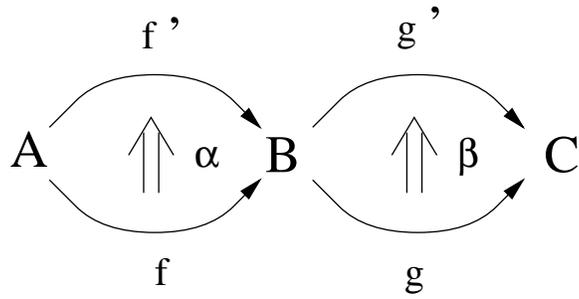}
}
\end{center}
\caption{Horizontal composition}
\label{twosq}
\end{figure}

\begin{figure}
\begin{center}
\mbox{
\epsfxsize=3in
\epsfbox{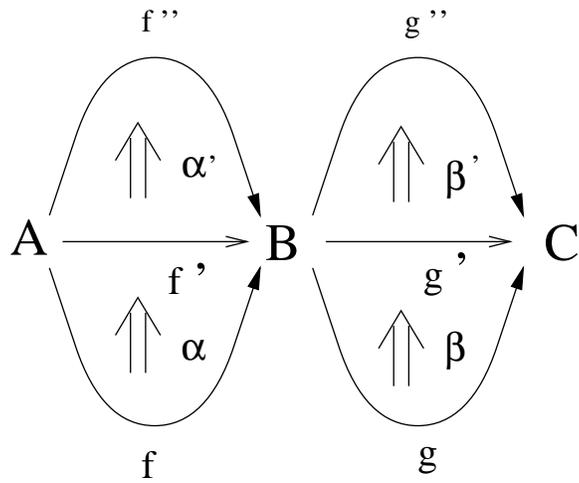}
}
\end{center}
\caption{Unambiguous composition}
\label{4sq4real}
\end{figure}

\begin{subsect}{\bf Review of 2-categories.\/}{\rm
A (small) strict
$2$-category consists of the following data:

(1) a set of {\it objects} {\bf Obj},

(2) a set of {\it $1$-morphisms} {\bf $1$-Mor}, whose elements have
{\it source } and
{\it target } objects,

(3) a set of {\it $2$-morphisms} {\bf $2$-Mor},  whose elements have
{\it source } and {\it target } $1$-morphisms.
Explicitly, given any two objects
$A,B$, there is a set of $1$-morphisms    
{\bf $1$-Mor}$(A, B)$ between them; the
object $A$ is called the source object, $B$ is the target.
Such a $1$-morphism $f$ is also denoted by $f: A \rightarrow B$.
 Given any
two 1-morphisms $f$, $g \in $ {\bf $1$-Mor}$(A, B)$
there is a set of 2-morphisms between
them; $f$ is the {\it source arrow ($1$-morphism)}, and $g$ is
 {\it the target arrow ($1$-morphism).}
The object $A$ is the {\it source object} and the object $B$ is the
{\it target object} (of the $2$-morphism).
Therefore
both the sources and targets
of $f$ and $g$
are required to coincide.
Such a $2$-morphism is also denoted by
$\alpha : f \Rightarrow g$.
(The {\it smallness} is the condition that these are sets.)

The following compositions of morphisms are defined
and related as indicated in items (1) through (5).

(1) For any $1$-morphisms $f : A \rightarrow B$ and
$g : B \rightarrow C$,
a $1$-morphism $g \circ f (=gf) : A \rightarrow C$
is the
composite.

(2) For any $1$-morphism $f  : A \rightarrow B$
and a $2$-morphism $\alpha : g \Rightarrow g' $
between $1$-morphisms $g, g' :  B \rightarrow C$,
there is
a composition (a $2$-morphism)
$\alpha f : gf \Rightarrow g'f$.
Similarly, there is a 2-morphism $g \alpha : gf \Rightarrow gf'$
when $\alpha: f \Rightarrow f'$.
These  compositions are
depicted on the left and middle of Fig.~\ref{compos}.

(3) For any $2$-morphisms $\alpha : f \Rightarrow g$ and
$\beta : g \Rightarrow h$ where $f,g,h : A \rightarrow B$,
there is
a composition  $\beta \cdot \alpha : f \Rightarrow h$.
This is depicted in Fig.~\ref{compos} right.
The composition $\beta \cdot \alpha$ is called the
{\it vertical composition} \index{vertical composition}
of 2-morphisms.

(4)
The composition depicted in
Fig.~\ref{twosq}
is unambiguous
in the sense that
$$ \beta \circ \alpha 
=(\beta f')\cdot (g \alpha) = (g' \alpha )\cdot( \beta f).
$$
The result  $\beta \circ \alpha$ is called the
{\it horizontal composition}
of the  2-morphisms $\alpha$ and $\beta$.

(5) The composition depicted in Fig.~\ref{4sq4real} is also unambiguous in the 
sense that
$$(\beta' \cdot \beta)\circ (\alpha' \cdot \alpha)
= (\beta' \circ \alpha ') \cdot (\beta \circ \alpha)$$
as 2-morphisms from $gf$ to $g'' f''$.

 We assume further that
$$ {\mbox{\rm id}}_f  \circ \alpha = f \circ \alpha $$
and
$$ \alpha \circ {\mbox{\rm id}}_g = \alpha g$$
where $f,g$ are  1-morphisms, $\alpha$ is a 2-morphism
and these composites are defined.

For any object $A$ (resp. $1$-morphism $f$),
 the identity $1$-morphism ${\mbox{\rm id}}_A$
(resp. $2$-morphism ${\mbox{\rm id}}_f$) is defined.
The identity $2$-morphism satisfies $( {\mbox{\rm id}}_f)\cdot \alpha =  \alpha$
,
$ \alpha \cdot ( {\mbox{\rm id}}_f) =  \alpha $ for any $2$-morphism $\alpha$.

For any $2$-morphism $\alpha : f \Rightarrow g $
where $f, g : A \rightarrow B $ , $ \alpha ( {\mbox{\rm id}}_A) = \alpha$
and $( {\mbox{\rm id}}_B) \alpha = \alpha$.

The following are the conditions for the {\it strictness}.

(1) The compositions of $1$- and $2$-morphisms are associative
( $(fg)h=f(gh)$, $(\alpha \beta) \gamma= \alpha (\beta \gamma)$ ).

(2) For any $1$-morphism $f: A \rightarrow B$,
$f ({\mbox{\rm id}}_A) = f = ({\mbox{\rm id}}_B) f $.

This concludes the definition of a small strict $2$-category.

Note that by composing morphisms we can represent $2$-morphisms by
planar polygons. More on 2-categories can be found also in \cite{CS:book}.

We follow  the definitions in \cite{Neu} for the rest of this section.

A {\it  monoidal $2$-category}
consists of a $2$-category ${\cal C} $
together with the data: 
(1) an object $I$, 
(2) for any object $A$ two $2$-functors
${\cal L}_A = A \sqt - : {\cal C}  \rightarrow {\cal C}$ and 
${\cal R}_A = - \sqt A : {\cal C} \rightarrow {\cal C}$
such that ${\cal L}_A (B) = {\cal R}_B (A) $ for any objects $A$, $B$,
and
(3) for any two $1$-morphisms $f: A \rightarrow A'$ and 
$g: B \rightarrow B'$ a $2$-morphism
$$ \sqt_{f,g} : (f \sqt B')(A \sqt g) \Rightarrow (A' \sqt g)(f \sqt B). $$

These data satisfy 8 conditions that we omit here. 

A {\it $2$-vector space}
 is a
$k$-linear additive category $({\cal C}, \oplus )$
that admits a subset $B \subset {\cal C}$ such that 
(1) any object $X \in B$ is simple, (2) for any $X \in B$, 
$dim_k(Hom(X,X))=1$, and (3) for any object $A \in {\cal C}$
there is a unique finite subset $B' \subset B$ 
such that $A \cong \oplus _{X_i \in B'} X_i ^{n_1}$. 
The set $B$ is called a basis and a $2-$vector space is finite if 
$B'$ is finite. 

The $2$-category of finite dimensional $2$-vector spaces is 
denoted by $2-vect$. The set of $1$-morphisms are $k$-linear functors,
and $2$-morphisms are natural transformations.
This is a subcategory of the $2$-category ${\cal C}_k$ of
small $k$-linear additive categories with the $k$-linear functors as 
$1$-morphisms. 

Theorem 2.7 of \cite{Neu} says that  $2-vect$
admits the structure of a strongly involutory 
moniodal 2-category 
in such a way that 
for ${\cal C}, {\cal D} \in 2-$vect, 
Obj$({\cal C} \sqt {\cal D}) = \{ (A_1, \cdots, A_n) | A_i $
 are pairs of objects of ${\cal C}$ and $ {\cal D}$, 
Mor$( (A_1, \cdots, A_n), (B_1, \cdots, B_m) )= 
\{ (f_{ij} ) | f_{ij} : A_j \Rightarrow B_i \} $.

A   {\it comoidal category} is a category ${\cal C}$ in ${\cal C}_k$
with the following data. 
Two $k$-linear functors
 $\diamondsuit : {\cal C} \rightarrow {\cal C} \sqt {\cal C}$ and
$\Gamma :  {\cal C} \rightarrow 
{\mbox{\rm vec}} 
$ where 
vec=the category of vector spaces,
and some  natural isomorphisms, including 
 $\alpha : (1  \sqt \diamondsuit) \circ \diamondsuit \Rightarrow 
	( \diamondsuit \sqt 1)  \circ \diamondsuit$,
satisfying certain relations, including the pentagon relation for $\alpha$.

A {\it $2$-bialgebra} or a {\it bimonoidal category} 
is a category ${\cal C} \in {\cal C}_k $ 
which has both
monoidal and comoidal category structures and with additional data
including a natural isomorphism 
$ \Xi : (\otimes \sqt \otimes )(1 \sqt {\cal T} \sqt 1)( \diamondsuit 
\sqt  \diamondsuit) \Rightarrow  \diamondsuit \otimes$
satisfying certain conditions,
where  ${\cal T}$ denotes
permutations,   including figures depicted in
 \ref{sjcube}, \ref{coh1}, and \ref{coh2}.
The faces in the figures 
are $2$-morphisms, and these cubes are relations among $2$-morphisms 
in the Hopf category.
}\end{subsect}

\begin{sect} {\bf Cocycle conditions  and Hopf categories.\/} 
\addcontentsline{toc}{subsection}{Cocycle conditions and Hopf categories}
{\rm 
The cocycles and their equations given in the previous
sections to construct partition functions
are given by Crane and Yetter \cite{CYex} 
to construct examples of Hopf categories. 
We refer the reader to \cite{CYex} and \cite{Neu} for the detailed 
relation between cocycles and Hopf categories.

Here we give a sketch description.
The Hopf category constructed is a categorification
of the Drinfeld quantum double of a finite group.
Thus first we recall $D(G)$, the quantum double of a finite group $G$.
It is an algebra generated by pairs of gropus and dual group elements,
$\{ (g, \hat{h}) \}$.
The multiplication is defined by 
$(g, \hat{h}) \cdot (k, \hat{\ell}) = \delta_{k^{-1} h k , \ell} (gk, \hat{\ell})$, comultiplication by
 $\Delta ( (g, \hat{h}) ) = \sum_k (g, \hat{k})
 \otimes (g, \widehat{hk^{-1}})$. 
The unit is $I=\sum_h (g, \hat{h})$,
 the counit is  $\Gamma (g, \hat{h})=\delta_{h, e}$,
 the antipode is $S(g, \hat{h}) = (g^{-1} , \widehat{ g h^{-1} g^{-1} } ) $.

The Hopf category ${\cal C} = {\cal D}(G)$, a categorification of $D(G)$, 
has the one dimensional vector space over a field $k$ generated by
pairs $\{ (g, \hat{h}) \}$ as objects. They are also simply denoted by pairs.
The functor
${\cal C} \sqt {\cal C} \rightarrow  {\cal C}$
is a family of isomorphisms 
$ (g, \hat{h}) \otimes   (k, \hat{\ell}) \rightarrow  (gk, \hat{\ell})$
where $h=k \ell k^{-1}$, 
and the functor $\diamondsuit : {\cal C} \rightarrow  {\cal C} \sqt {\cal C} $
is given by a family $  (g, \hat{h}) \rightarrow \sum_k (g, \hat{k}) 
\sqt  (g, \widehat{xk^{-1}})$.

The transformation ($2$-morphism) 
$$\alpha _{ (g, \hat{h}),  (k, \hat{\ell}), (m, \hat{n}) } :
(  (g, \hat{h})\sqt (k, \hat{\ell}) ) \sqt  (m, \hat{n})
\rightarrow  (g, \hat{h})\sqt (  (k, \hat{\ell})\sqt  (m, \hat{n}) )$$
is a family of scalars 
$\alpha (g,k,m; \hat{n})$ since the objects are nontrivial 
only when 
$\ell = k^{-1} hk$ and $n= \ell ^{-1} m \ell$.
Similar for $\beta(g; \hat{i}, \hat{j}, \hat{k})$, 
and for the transformation $\Xi$, 
one has scalars 
$$\phi(g,k; \hat{m}, \hat{n}):
\oplus _{mn=\ell } (gk, \hat{m}) \sqt (gk, \hat{n} )
\rightarrow \oplus _{ab=h } (gk, \widehat{k^{-1} ak}) \sqt (gk, \widehat{k^{-1} bk} ) .$$

Thus we get these cocycles as transformations (or $2$-morphisms)
in a Hopf category.
The cocycle conditions are obtained from the axioms of the Hopf categories.

} \end{sect}

\begin{figure}
\begin{center}
\mbox{
\epsfxsize=4in 
\epsfbox{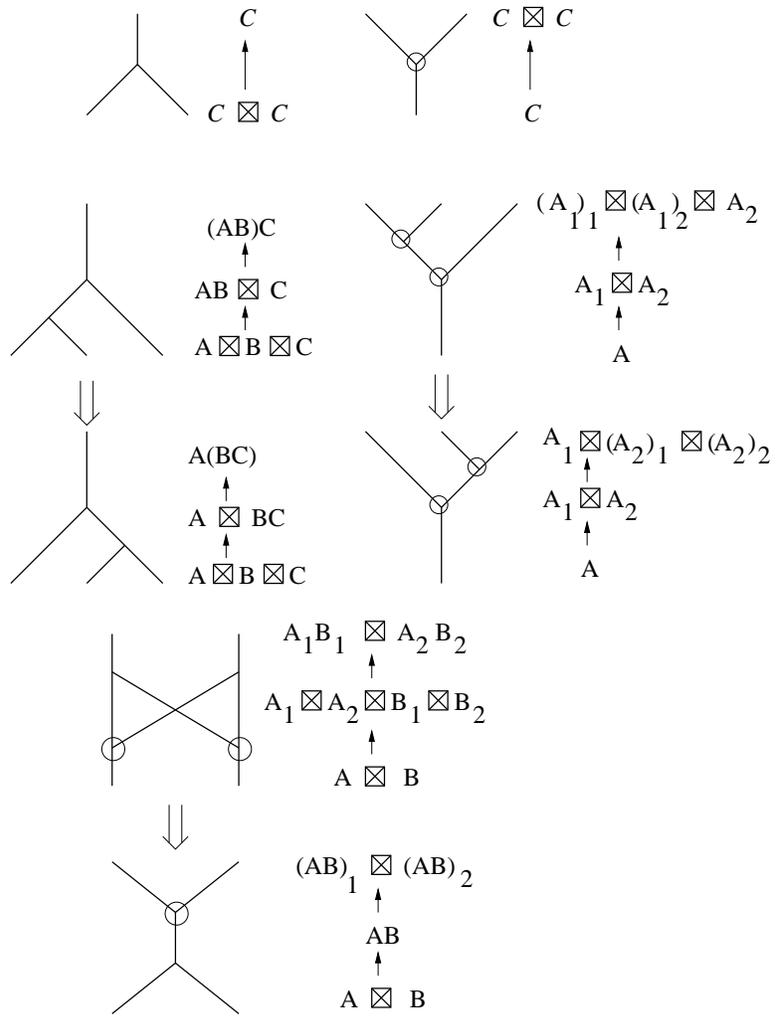}
}
\end{center}
\caption{Diagrams for morphisms}
\label{why}
\end{figure}

\begin{sect} {\bf Diagrams for morphisms.\/} 
\addcontentsline{toc}{subsection}{Diagrams for morphisms}
{\rm 
We use the diagrammatic convention depicted in Fig.~\ref{why}.
The  
top two figures represent tensor and cotensor functors ($1$-morphisms).
We read the diagram from bottom to top,  the segments
correspond to categories, and trivalent vertices correspond to 
tensor (uncircled vertices) functors and cotensor (circles ones) functors.
Thus a tree diagrams represent compositions of such functors.
The bottom three figures represent $2$-morphisms.
For each of these, two compositions of functors are related by a double arrow
representing a transformation ($2$-morphism) between $1$-morphisms that are compositions
of tensor and cotensor functors.

For the categorification of the quantum double, the objects are 
represented by pairs $\{ (g, \hat{h}) \}$,
therefore the edges of the graphs are labled by these pairs.
The situation is seen in Fig.~\ref{rule}.
The functors, then, are represented by trivalent vertices of the graphs.
These are depicted in Figs.~\ref{alpha}, \ref{phi}.
In Fig.~\ref{alpha}, the functor corresponding to $\alpha$ is depicted, 
and in Fig.~\ref{phi}, the functor corresponding to $\phi$
(or $\Xi$) is depicted. The diagrams for $\beta$ is similar to $\alpha$,
except that the vertices have small circles.

} \end{sect}

\clearpage

\begin{figure}
\begin{center}
\mbox{
\epsfxsize=5in
\epsfbox{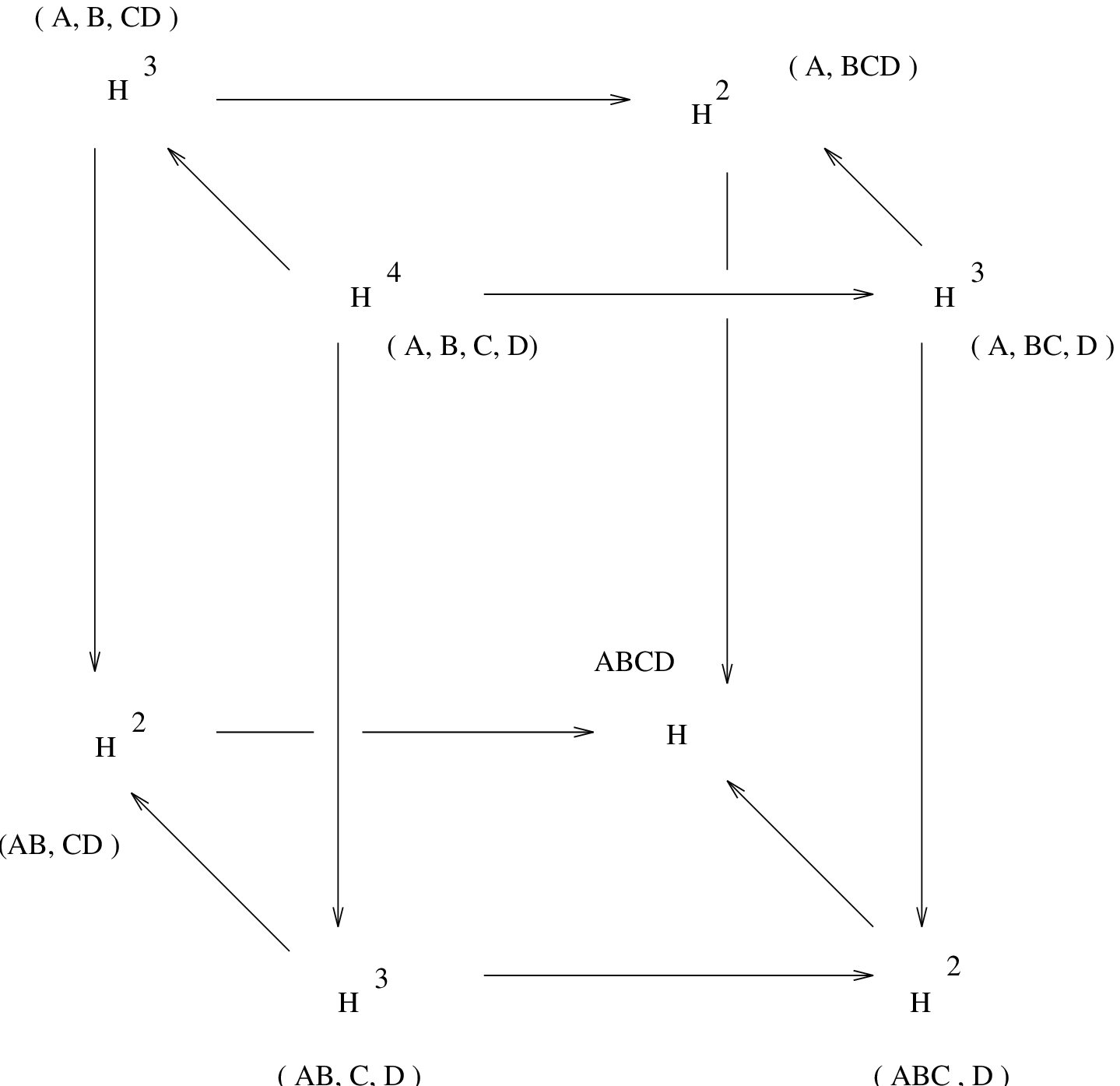}
}
\end{center}
\caption{The coherence cube for tensor operators}
\label{sjcube}
\end{figure}

\newpage

\begin{figure}
\begin{center}
\mbox{
\epsfxsize=5in 
\epsfbox{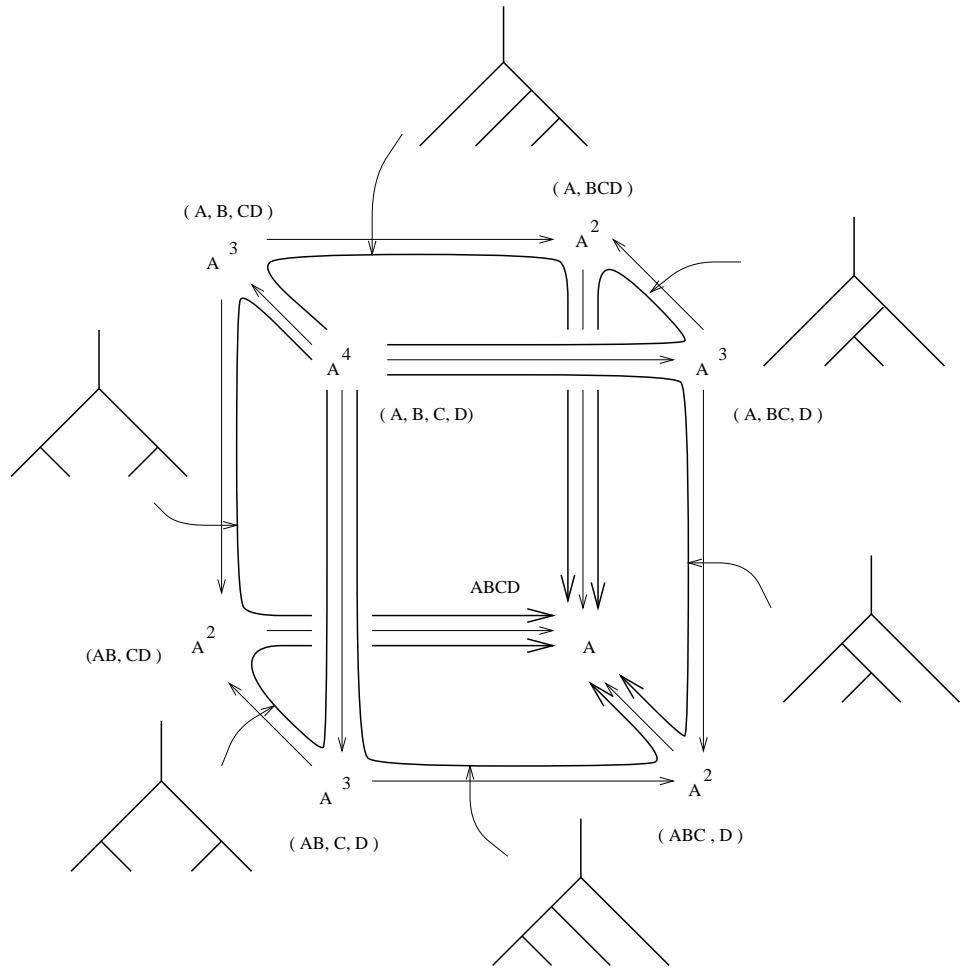}
}
\end{center}
\caption{Networks for tensor operators }
\label{netsj}
\end{figure}

\newpage 

\begin{figure}
\begin{center}
\mbox{
\epsfxsize=5in 
\epsfbox{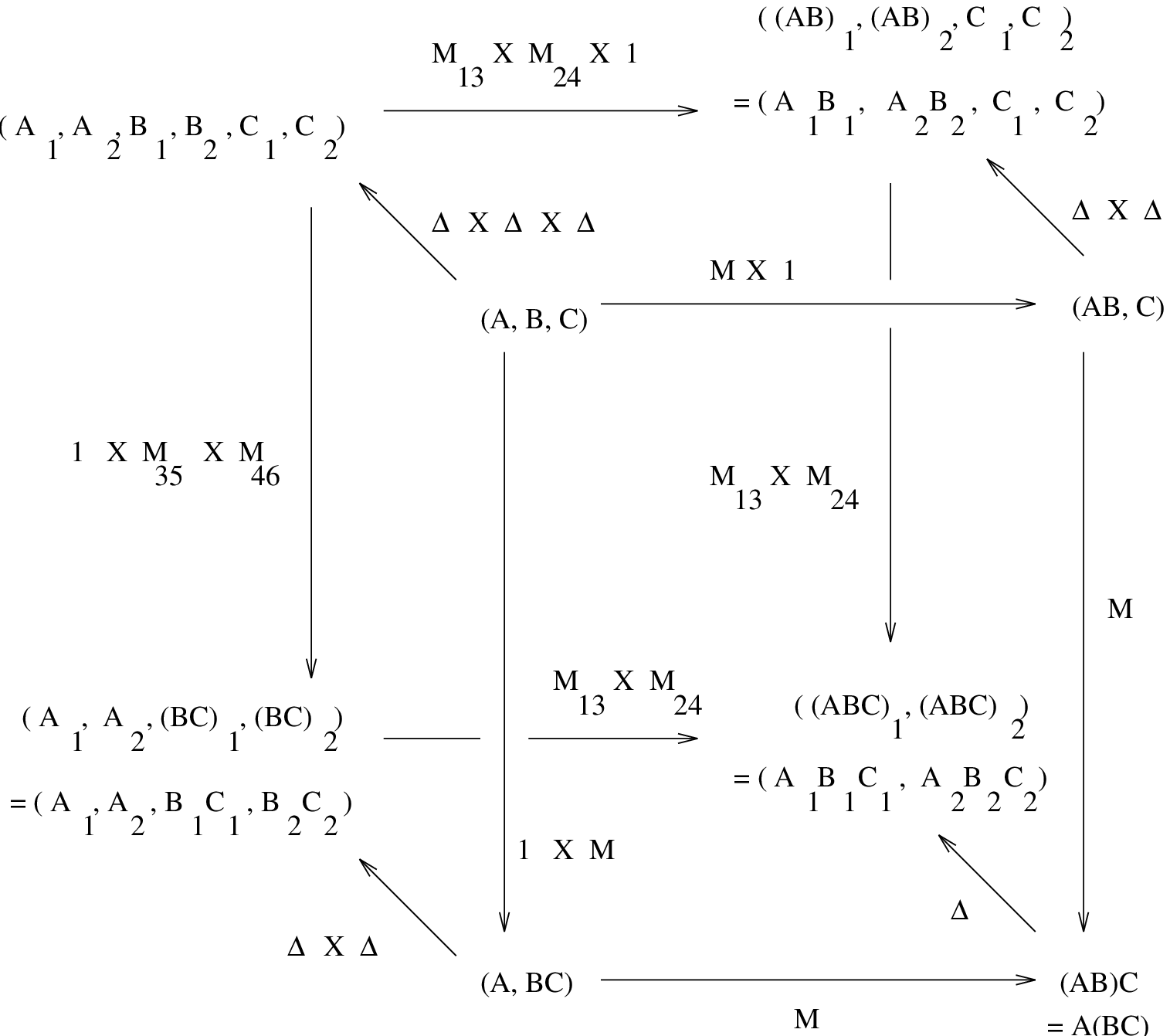}
}
\end{center}
\caption{The coherence cube, type I}
\label{coh1}
\end{figure}

\newpage

\begin{figure}
\begin{center}
\mbox{
\epsfxsize=6in 
\epsfbox{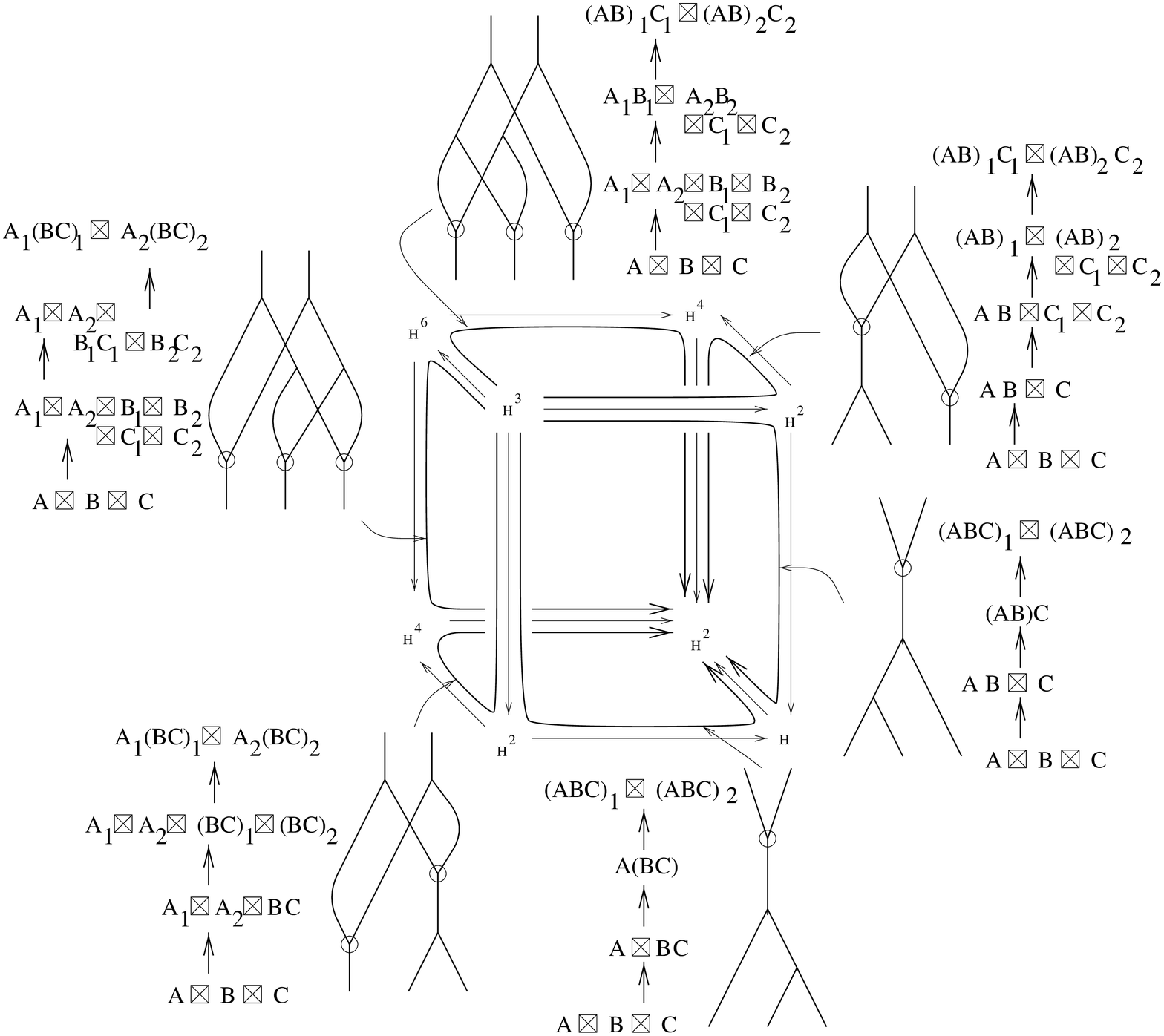}
}
\end{center}
\caption{Networks for the coherence cube I}
\label{coh11a}
\end{figure}

\newpage 

\begin{figure}
\begin{center}
\mbox{
\epsfxsize=4in
\epsfbox{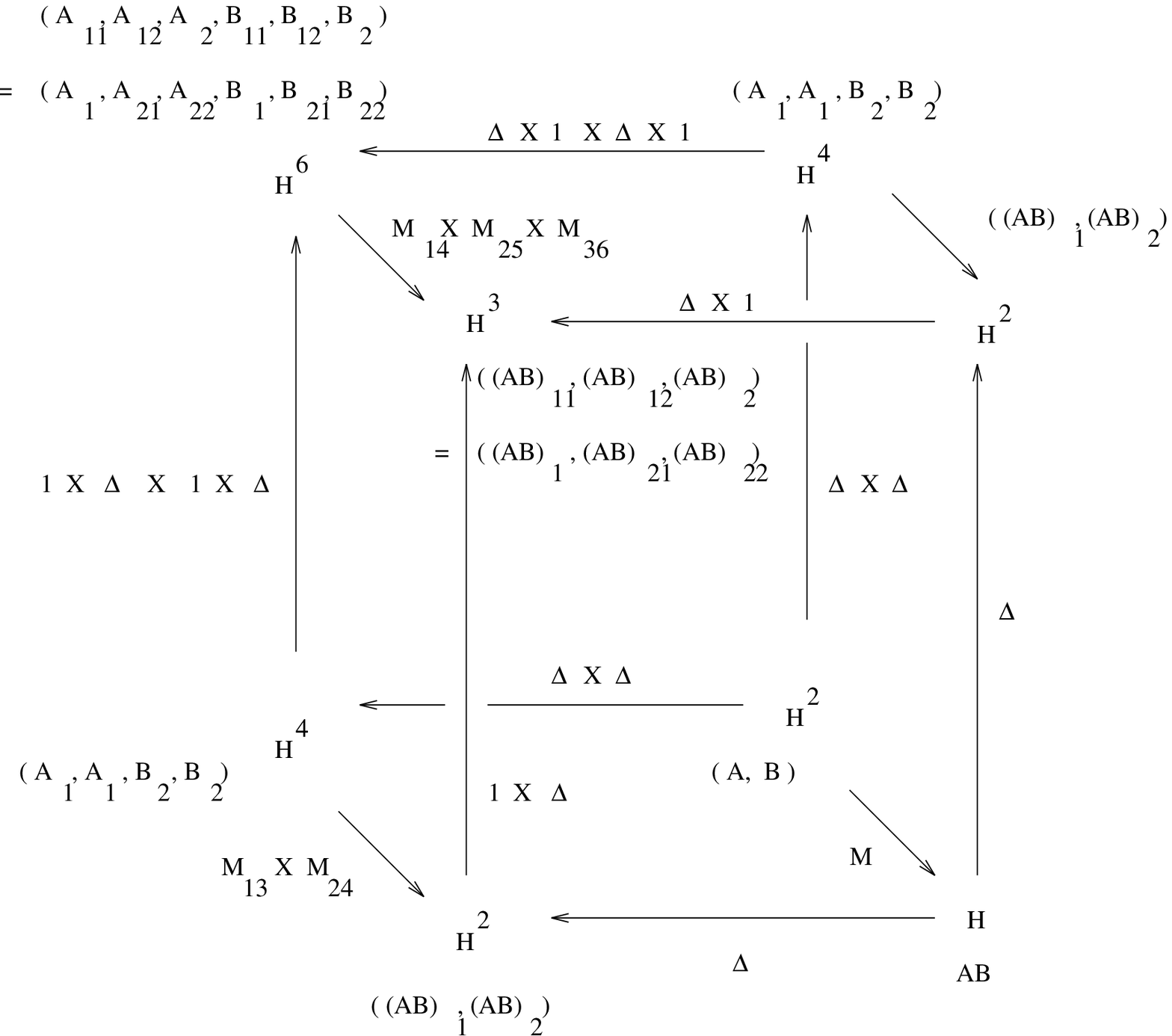}
}
\end{center}
\caption{The coherence cube, type II}
\label{coh2}
\end{figure}
 
\newpage 

\begin{figure}
\begin{center}
\mbox{
\epsfxsize=4in 
\epsfbox{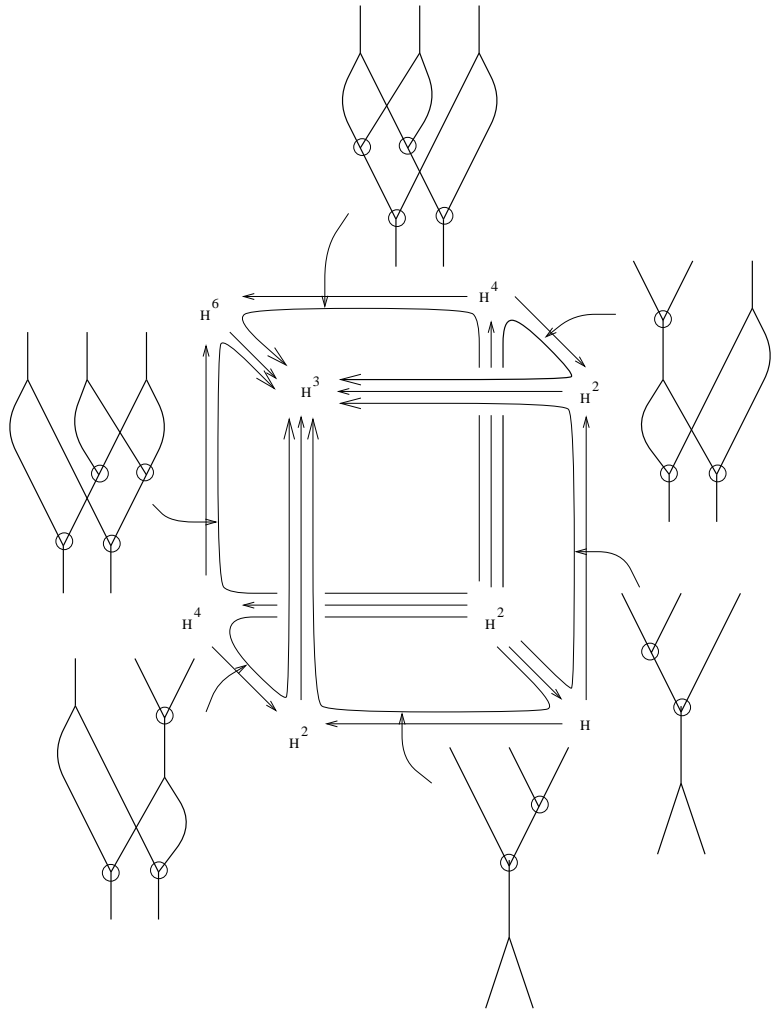}
}
\end{center}
\caption{Networks for the coherence cube II}
\label{coh22}
\end{figure}

\clearpage

\begin{sect}{\bf Coherence cubes and diagrams.\/}
\addcontentsline{toc}{subsection}{Coherence cubes and diagrams}
{\rm 
Certain compositions of $2$-morphisms in Hopf categories 
are required to satisfy relations, such as those called
{\it coherence cubes} in \cite{CF}.
There are relations between the  coherence cubes in
\cite{CF} and the diagrams (graphs) we have defined in this paper. 
Here we include Figures~\ref{netsj}, \ref{coh11a}, 
\ref{coh22} that indicate such relations. 
Figures~\ref{sjcube}, \ref{coh1}, \ref{coh2} are the coherence cubes
given in \cite{CF}. 
The labels in the parentheses are in fact (co-)tensor products
and equalities are morphisms.
There is another cube in \cite{CF} which is dual to Figure~\ref{sjcube}
and is ommited.

In each of these figures, there is an initial vertex and a terminal
vertex in the cube.
These cubes appear in \cite{CF} in the definition of Hopf categories.
 Each arrow from the initial vertex to 
the terminal vertex represents a morphism, which is represented
by a diagram with trivalent vertices. 

Through each rectangular face, such a path is homotoped to
another such path.  
This homotopy causes a change in diagrams.
These changes in diagrams are the scenes from our movies, and in the 
case of the Hopf category associated to a Drinfeld double, these scenes
are group cocycles.

Note that relations are represented by diagrams (graphs)
and the relations take forms of equating two sequences (graph movies)
of graphs. 
These equations are in fact identical that appeared in cocycle conditions
(Figs.~\ref{cocymovie1}, \ref{cocymovie2}) 
since these cocycles came from constructions of 
Hopf categories in \cite{CY}. 
On the other hand these diagrams are related to the 
taco move via the dual graphs. Thus we see that the the formalism proposed 
in \cite{CF} can          
be used to construct invariants of 4-manifolds via a 
state-sum  provided the propsed set of states is finite.
}\end{sect}

%
%
%

\section{Concluding remarks}

In this paper we established diagrammatic machinery 
for the study of $4$-manifold invariants using triangulations
and graphs.
In particular, invariance under Pachner moves  of Crane-Frenkel invariants 
for cocycles constructed by Crane-Yetter is proved 
by using graphs. 
This 
strongly suggests generalizations of
 the Dikgraaf-Witten invariant to $4$-manifolds
using cocycles defined in \cite{CY}.
We have shown direct relations among algebraic
structures (Hopf categories), triangulations,
 and (graph) diagrams in dimension $4$,
generalizing spin network theory in $3$-dimensions.

Further study on higher dimensional TQFTs and 
higher algebraic structures are 
anticipated. 
We expect that our diagrammatic machinery established in this paper 
serve as tools for further developments in the area.

Open questions remain. Which finite groups contain cocycles 
 that satisfy the symmetry conditions? Can other examples of Hopf Categories be constructed. What do  the invariants from this construction  measure?
Can these invariants be related to invariants that arise from Donaldson 
Theory?

\section*{Acknowledgements}

We are grateful to J Barrett, L. Crane,  
 and D. Yetter for valuable conversations.
The first named author 
was supported by the NSA while some of the research for this paper
was being conducted.
The second named author is partially supported by NSF DMS-2528707. 
The third named author is partially  supported by 
the University of South Florida
Research and Creative Scholarship Grant Program
under Grant Number 1249932R0.

\end{document}